\documentclass[8pt]{article}
\usepackage{amssymb}
\usepackage{amsmath}
\usepackage{amsfonts}
\usepackage{amscd}
\usepackage{amsthm}
\usepackage{graphicx}
\usepackage{enumerate}
\usepackage{empheq}
\usepackage{xcolor}
\usepackage{caption}
\usepackage{subcaption}
\newtheorem{theorem}{Theorem}[section]

\newtheorem{lemma}{Lemma}[section]
\newtheorem{example}{Example}[section]

\newcommand{\se}{\setcounter{equation}{0}}

\newcommand{\R}{\mathbb{R}}

\newcommand{\bu}{\textbf{u}}
\newcommand{\bU}{\textbf{U}}
\newcommand{\bW}{\textbf{W}}
\newcommand{\bX}{\textbf{V}}

\newcommand{\bV}{\textbf{V}}

\newcommand{\bS}{{\mathcal{S}}}
\newcommand{\bR}{\textbf{R}}
\newcommand{\bP}{\mathcal{P}}
\newcommand{\bH}{\textbf{H}}
\newcommand{\bL}{\textbf{L}}
\newcommand{\bv}{\mathbf{v}}
\newcommand{\bw}{\textbf{w}}
\newcommand{\de}{\textbf{e}}
\newcommand{\bE}{\textbf{E}}
\newcommand{\bn}{\textbf{n}}
\newcommand{\bz}{\textbf{z}}

\newcommand{\bJ}{\textbf{G}}
\newcommand{\brf}{\mathbf{f}}
\newcommand{\al}{\alpha}
\newcommand{\Lam}{\Lambda}
\newcommand{\be}{\mbox{\boldmath$\eta$}}
\newcommand{\bphi}{\mbox{\boldmath$\phi$}}
\newcommand{\bxi}{\mbox{\boldmath$\omega$}}
\newcommand{\bpsi}{\mbox{\boldmath$\psi$}}

\newcommand{\bzeta}{\mbox{\boldmath$\zeta$}}
\newcommand{\bchi}{{\mbox{\boldmath$\chi$}}}
\newcommand{\brho}{\mbox{\boldmath$\rho$}}
\newcommand{\bT}{\mbox{\boldmath$\Theta$}}

\newcommand{\biota}{{\mbox{\boldmath$\iota$}}}
\newcommand{\bp}{{\bf {p}}}
\newcommand{\vphi}
{{\underline{\boldsymbol{\varphi}}}}
\newcommand{\e}{{\bf e}}
\newcommand{\su}{{\underline{\boldsymbol{\sigma}}}}
\newcommand{\sib}{{{\boldsymbol{\sigma}}}}
\newcommand{\tu}{{\underline{\boldsymbol{\tau}}}}

\newcommand{\mo}{\mbox{$\mathcal{F}$}}

\newcommand{\Su}{{\underline{\boldsymbol{\Sigma}}}}
\newcommand{\betau}{{\underline{\boldsymbol{\beta}}}}
\newcommand{\Y}{{\underline{\textbf{V}}}}
\newcommand{\dbt}{{\bar\partial_t}}

\newcommand{\bx}{{\textbf{x}}}
\newcommand{\bsq}{\boldsymbol{[}}
\newcommand{\bqs}{\boldsymbol{]}}
\newcommand{\bZ}{\textbf{Z}}
\newcommand{\by}{\textbf{y}}
\newcommand{\vtheta}
{{{\boldsymbol{\vartheta}}}}
\DeclareMathOperator*{\esssup}{ess\: sup}
\newcommand{\vertiii}[1]{{\left\vert\kern-0.25ex\left\vert\kern-0.25ex\left\vert #1 
		\right\vert\kern-0.25ex\right\vert\kern-0.25ex\right\vert}}

\begin{document}
\author{Debendra Kumar Swain\footnotemark[1] \thanks{Tata Institute of Fundamental Research, Center for Applicable Mathematics, Bengaluru--560065, India. Email: debendra26@tifrbng.res.in,},~ Saumya Bajpai\footnotemark[2]\thanks{School of Mathematics and Computer Science, Indian Institute of Technology Goa, Ponda, Goa-403401, India. Email: saumya@iitgoa.ac.in, 
}, ~and Deepjyoti Goswami\footnotemark[2]\thanks{Department of Mathematical Sciences, Tezpur University, Tezpur, Sonitpur, Assam--784028, India. Email: deepjyoti@tezu.ernet.in}}
\date{}
\title
	{On global error estimates of a local discontinuous Galerkin method of the Boussinesq system of equations}
	\maketitle
	\begin{abstract}
		In this paper, we develop and analyze a local discontinuous Galerkin (LDG) finite element method for the Boussinesq system of equations, and  employ an upwind scheme to effectively manage the nonlinear convective  terms, thereby ensuring stability and accuracy. Our analysis establishes new \emph{a priori} regularity estimates for the semidiscrete LDG approximations.  We establish optimal \emph{a priori} error estimates for the velocity approximation in  $\bL^2$--norm, the pressure approximation in  $L^2$--norm  and for temperature in $L^2$--norm for all $t>0$.
Our proof relies on the standard elliptic  and parabolic duality arguments and  modified
Stokes  operators defined on appropriate broken Sobolev spaces.  Under suitable smallness assumptions on the given data, we further establish error estimates that remain uniform in time. Furthermore, the backward Euler
scheme is considered for a full discretization, and optimal fully discrete error estimates
are derived. Finally, we carry out numerical experiments to substantiate our theoretical
findings. 				
\end{abstract}
	
	\vspace{1em} 

\section{Introduction}
Heat transfer in a viscous incompressible heat--conducting ﬂuid is represented by the Boussinesq system of equations, which is described by the
following system of partial diﬀerential equations:
\noindent
 \begin{eqnarray}\label{1.1}
 \left.
 \begin{array}{rcl}
&& { \bu}_{ t}\,-\,\nu \,\Delta \bu\,+ \bu\cdot \nabla \bu\,+\nabla p + \,\boldsymbol{\mu } \theta
\,={\bf f} \hspace{1cm} {\bf x}\in \,\Omega, \;\;t\,>\,0,\\
&&\theta_t -k\Delta\theta + \textbf{u} . \nabla\theta = g \hspace{1.4cm} {\bf x}\in \,\Omega, \;\;t\,>\,0,
\end{array}
 \right\}
\end{eqnarray}
and incompressibility condition
\begin{eqnarray}
\label{1.2}
\nabla \cdot \bu\,=0\hspace{1cm}{\bf x}\in \Omega,\,t\,>\,0,
\end{eqnarray}
with initial and boundary conditions
\begin{eqnarray}
\label{1.3}
\bu(\,{\bf x},\,0)= \bu_0, \; \; \theta(\textbf{x},\, 0)\, =\, \theta_0 \;\;\;\mbox {in} \;\Omega,\;\;\;\;\; \bu\,=0,\; \theta=0\; \;
\mbox {on}\; \partial \Omega,\;\, t\,\ge 0,
\end{eqnarray}
where, \(\bu_t = \frac{\partial \bu}{\partial t}\) and \(\theta_t = \frac{\partial \theta}{\partial t}\). The domain \(\Omega \subset \mathbb{R}^2\) is assumed to be a bounded, convex polygonal region with smooth boundary \(\partial \Omega\). The vector  \(\mathbf{u} = (u_1, u_2)\) denotes the fluid velocity, \(\theta\) represents the temperature, and \(p\) corresponds to the pressure. The term \(\mathbf{f}\) describes the external body force, while \(g\) denotes the prescribed heat source. The constants \(\nu > 0\) and \(\kappa > 0\) stand for the kinematic viscosity and thermal conductivity, respectively.\\
 We further suppose that \(\boldsymbol{\mu} \in \mathbb{R}^2\) and enforce a zero--mean condition on the pressure, 
 that is, $\displaystyle{\int_\Omega} p\, d\Omega=0$.  For a detailed discussion on the mathematical modeling and applications of this system, we refer to \cite{Benard, BHH, KO} and the references cited therein.\\
Much of the earlier work related to the system (\ref{1.1})--(\ref{1.3}) is based on continuous Galerkin methods, while studies involving discontinuous Galerkin (DG) finite element analysis remain relatively limited. Recently, DG methods are increasingly applied to the numerical treatment of partial differential equations in many areas of science and engineering. They offer several benefits, including geometric flexibility, local conservation properties, variable polynomial approximation, and high accuracy. A comprehensive discussion of DG methods can be found in \cite{CKS00}. Extensive studies on DG approximations for both steady and time-dependent incompressible Navier--Stokes equations are reported in \cite{GRW005, GRW05, KR05, MLR22}. Various DG discretization strategies for incompressible NSEs are further discussed in \cite{CCQ17, CKS09, KR19, WDYW17} and the references therein. More recently, DG methods have also been applied to system (\ref{1.1})--(\ref{1.3}) in \cite{BD24B}, where optimal error estimates for velocity, pressure, and temperature have been established.\\
\noindent
The LDG technique, a specific development of the DG approach, distinguishes itself by locally enforcing penalty terms, particularly along element interfaces, to bolster both stability and conservation. In contrast to conventional DG methods, which may apply penalties more globally, the LDG technique prioritizes local conservation by maintaining stability and continuity within each individual element. The selective use of penalty terms together with the associated numerical flux evaluation is a characteristic feature of the LDG framework, which enhances its stability and robustness on an element--by--element level. Because of these properties, LDG schemes are widely applied to problems arising in physics, fluid dynamics, and engineering, particularly for models governed by elliptic and parabolic equations, as well as the incompressible steady Stokes and Navier--Stokes systems.\\
\noindent
Castillo  {\it et al.} \cite{CCPS2000} have presented the first \textit{a priori} error analysis for  LDG method applied to a model elliptic problem and have shown that, when polynomials of total degree at least $k$ are used, the $L^2$-norms of the gradient and the potential converge with orders $k$ and $k+\tfrac{1}{2}$, respectively, on arbitrary meshes; moreover, if stabilization parameters of order $h^{-1}$ are employed, the convergence order of the potential increases to $k+1$.  Cockburn {\it et al.} \cite{cockburn00} have also established a superconvergence result on Cartesian grids using tensor–product polynomials of degree at most $k$, where a specially chosen numerical flux yields convergence orders $k+\tfrac{1}{2}$ for the gradient and $k+1$ for the potential.
 Furthermore, the authors in \cite{cockburn02}  have analyzed LDG method for the Stokes system and have derived optimal $L^2$ error estimates for velocity and pressure on shape-regular meshes with hanging nodes. Optimal convergence is obtained using degree $k$ polynomials for velocity and degree $k-1$ for pressure $(k\ge1)$, while using degree $k$ for all variables preserves the same orders with better efficiency. 
  Guzman \cite{guzman08} has established local and weighted pointwise error estimates for the LDG scheme applied to the Stokes problem in two and three dimensions, obtaining optimal pointwise convergence results for both the velocity and pressure on smooth domains. Cockburn {\it et al.} \cite{cockburn03}   have proposed and analyzed LDG formulation for the Oseen equations governing incompressible fluid flow. For shape--regular meshes allowing hanging nodes, they have established optimal \textit{a priori} error bounds for both velocity and pressure in the $L^2$ as well as in negative--order norms. 
  Similarly, in \cite{cockburn04}, optimal error estimates have been obtained for the velocity in the $\bL^2$--norm and for the stress tensor in the $L^2$--norm for a locally conservative LDG formulation of the incompressible Navier--Stokes equations. Furthermore, in \cite{WDYW17}, the characteristic method has been combined with the LDG framework to construct variational formulations for the time--dependent incompressible Navier--Stokes equations in $\mathbb{R}^2$. Nonlinear stability for a symmetric formulation has also been established, and on general triangulations. Suboptimal {\it a priori} estimates for the velocity in the $\bL^2$--norm together with optimal error estimates for the pressure in the $L^2$--norm have been derived.\\
To the best of our knowledge, although significant efforts have been made in this area, the use of the LDG method for unsteady incompressible fluid problems has not been thoroughly investigated. The present work aims to address this gap. In particular, we have established optimal error estimates for both semidiscrete and fully discrete schemes by incorporating  backward Euler scheme within the LDG framework for system (\ref{1.1})--(\ref{1.3}).
\noindent
The main analytical tools and techniques utilized in this work are described as follows:
\begin{enumerate}
    \item[$\bullet$] Initially, we derive \emph{a priori}  regularity bounds for the semidiscrete LDG approximations of the stress tensor, velocity, temperature and temperature gradient. 
    \item[$\bullet$] Next, we introduce the approximation operator $\bR_h$ along with the $L^2$-{ projections} $\underline{\bP}_h^\sib$, $\bP_h^\bu$, $\mathcal{P}_h^\bchi$  and $\mathcal{P}_h^\theta$ (see Section \ref{ldgsection3}), and discuss their approximation properties. We then define the modified modified Stokes { operators} $\underline{\boldsymbol{\bS}}_h^\sib$, $\boldsymbol{\bS}_h^\bu$, $\boldsymbol{\bS}_h^\bchi$  and $\mathcal{S}_h^\theta$ (see Section \ref{ldgsection3}), which are tailored for broken Sobolev spaces within the LDG context and play a pivotal role in deriving semidiscrete error estimates.
    \item[$\bullet$] By using the modified Stokes operators along with duality arguments, we subsequently establish optimal \emph{a priori} error estimates for the semidiscrete LDG approximations. Under the smallness condition on the data, we achieve uniform in time estimates. 
    \item[$\bullet$] We further extend our analysis to the fully discrete setting by applying  backward Euler method to the semidiscrete LDG Boussinesq system of  equations, thereby deriving optimal \emph{a priori} error estimates for the fully discrete approximations of velocity, pressure and temperature.
    \item[$\bullet$] To conclude, numerical experiments are performed to validate the theoretical error estimates.
\end{enumerate}
\noindent
The organization of this chapter is described as follows:
Section \ref{ldgsection2} introduces the functional setting required for the analysis, including the relevant spaces, notation, and fundamental assumptions. In addition, the LDG weak formulation of the problem is presented.
Section \ref{ldgsection3} is devoted to the semidiscrete LDG scheme. We first establish \emph{a priori} regularity estimates for the semidiscrete solution and discuss the necessary trace inequalities. The analysis further examines the properties of the modified Stokes operators in broken Sobolev spaces within the LDG framework.
Section \ref{ldgsection4} derives optimal \emph{a priori} error estimates for the velocity and temperature using new approximation properties. Under suitable smallness conditions on the data, uniform--in--time $(t>0)$ error bounds are obtained. Section \ref{ldgsection5} establishes optimal error estimates for the pressure.
Section \ref{ldgsection6} develops the fully discrete scheme based on the backward Euler time discretization and derives the fully discrete error bounds.
Section \ref{ldgsection7} provides numerical experiments that illustrate and validate the theoretical convergence results.
Finally, Section \ref{ldgsection8} summarizes the key findings and contributions of this work.
\section{Preliminaries and  Weak Formulation}\label{ldgsection2}
\se
Throughout this article, bold symbols are used to denote vector--valued function spaces.
In particular, $\mathbf{H}_0^1(\Omega) := (H_0^1(\Omega))^2$ and
$\mathbf{L}^2(\Omega) := (L^2(\Omega))^2$ represent $\mathbb{R}^2$--valued spaces.
Similarly, underline notation is employed for $\mathbb{R}^4$--valued spaces, for example,
$\underline{\bH}_0^1(\Omega) := (H_0^1(\Omega))^4$ and
$\underline{\bL}^2(\Omega) := (L^2(\Omega))^4$.
For the analytical framework of this work, we make use of the standard Sobolev spaces
(see \cite{A75}). For a nonnegative integer $n$ and $q \geq 1$, these spaces are defined by
\[
\mathcal{V}^{n,q}(\eta) = \left\{ \eta \in L^q(\Omega) : \forall \: |{\bf p}| \leq n, \, \partial^{\bf p} \eta \in L^q(\Omega) \right\}.
\]
Here, ${\bf p}=(p_1,p_2)$ denotes a multi--index with $|{\bf p}|=p_1+p_2$, and $\partial^{\bf p} \eta$
represents the partial derivative of order $|{\bf p}|$.
The norm and seminorm on $\mathcal{V}^{n,q}(\Omega)$ are denoted by
$\|\cdot\|_{n,q,\Omega}$ and $|\cdot|_{n,q,\Omega}$, respectively.
The inner product in $L^2(\Omega)$ is written as $(\cdot,\cdot)$.
For the Hilbert space $H^n(\Omega) := W^{n,2}(\Omega)$, we use the notation
$\|\cdot\|_{n,\Omega}$ or simply $\|\cdot\|_n$.
The vector--valued space $\mathbf{H}_0^1(\Omega)$ is equipped with the norm
\[
\|\nabla \boldsymbol{\eta}\|
=
\left(
\sum_{m,n=1}^{2} (\partial_n \eta_m, \partial_n \eta_m)
\right)^{1/2}
=
\left(
\sum_{m=1}^{2} (\nabla \eta_m, \nabla \eta_m)
\right)^{1/2}.
\]
\noindent
We additionally introduce the quotient space $H^n(\Omega)/\mathbb{R}$, which is equipped with the norm
\(
\|\varphi\|_{H^n(\Omega)/\mathbb{R}}
:=
\inf_{c \in \mathbb{R}} \|\varphi + c\|_m .
\)
In particular, for $n=0$, this construction yields the space $L^2(\Omega)/\mathbb{R}$.
Let $\mathcal{Y}$ be a Banach space and let $1 \leq q \leq \infty$.
The Bochner space $L^q(0,T;\mathcal{Y})$ is defined as the collection of $\mathcal{Y}$--valued measurable functions $\boldsymbol{\zeta}$ on the interval $(0,T)$ satisfying
\begin{align*}
\int_0^T \|\boldsymbol{\zeta}(t)\|_{\mathcal{Y}}^q \, dt < \infty
\quad \text{for } 1 \leq q < \infty,\quad 
\text{and}\quad 
\esssup_{0<t<T} \|\boldsymbol{\zeta}(t)\|_{\mathcal{Y}} < \infty
\quad \text{for } q=\infty.
\end{align*}
The dual space of $H^n(\Omega)$ is denoted by $H^{n}(\Omega)$ and is endowed with the norm
\[
\|f\|_{-n}
:=
\sup \left\{
\frac{(f,\chi)}{\|\chi\|_n}
\,:\,
\chi \in H^n(\Omega),\; \chi \neq 0
\right\}.
\]
Finally, we define the subspace of solenoidal vector fields by
\[
\bJ_1
=
\left\{
\bv \in \bH_0^1(\Omega)
\,:\,
\nabla \cdot \bv = 0
\right\}.
\]
The analysis in this paper is carried out under the following assumptions:
\\
({\bf A1}). 
For a positive constant \( M_0 \),
the initial velocity \( \bu_0 \in \bH^2(\Omega)\cap \bJ_1 \), the initial temperature \( \theta_0 \in H^2(\Omega) \) and the external force with its time derivatives \( \brf,\, \brf_t \in L^\infty(0,\,\infty:\,\bL^2(\Omega))\), the given heat source with its time derivatives \( g,\, g_t \in L^\infty(0,\,\infty:\,L^2(\Omega))\) satisfy 
\[\|\theta_0\|,\,\|\bu_0\|_2\leq M_0,\quad  \mbox{and} \quad\displaystyle\sup_{0<t< \infty}\{\|\brf(\cdot,\,t)\|,\|\brf_t(\cdot,\,t)\|, \|g(\cdot,\,t)\|,\|g_t(\cdot,\,t)\|\}\leq M_0.\]
To introduce the LDG formulation, we consider vector  
$\mathbf{u}, \mathbf{z} \in \bL^{2}(\Omega)$ and tensor  
$\underline{\boldsymbol{\eta}},\, \underline{\boldsymbol{\omega}} \in \underline{\bL}^2(\Omega)$, where
$\mathbf{u} = (u_{1}, u_{2})$, $\mathbf{z} = (z_{1}, z_{2})$,
$\underline{\boldsymbol{\eta}} = [\eta_{ij}]_{1 \le i,j \le 2}$, and
$\underline{\boldsymbol{\omega}} = [\omega_{ij}]_{1 \le i,j \le 2}$.
We employ the conventional definitions
$(\nabla \mathbf{v})_{pq} = \partial_{q} v_{p}$ and
\(
(\nabla \cdot \underline{\boldsymbol{\eta}})_{p}
= \sum_{q=1}^{2} \partial_{q} \eta_{pq}.
\)
Moreover, for vectors $\mathbf{u}$ and $\mathbf{z}$, we denote by
$\mathbf{u} \otimes \mathbf{z}$ the matrix whose $(p,q)$--th entry is
$u_p z_q$, and define
\begin{align*}
\underline{\boldsymbol{\eta}} : \underline{\boldsymbol{\omega}}
:= \sum_{p,q=1}^{2} \eta_{pq}\omega_{pq},
\quad \text{and} \quad 
\mathbf{u} \cdot \underline{\boldsymbol{\eta}} \cdot \mathbf{z}
:= \sum_{p,q=1}^{2} u_p \eta_{pq} z_q
= \underline{\boldsymbol{\eta}} : \mathbf{u} \otimes \mathbf{z}.
\end{align*}
With these notations in place, the Boussinesq system
\eqref{1.1}--\eqref{1.3} can be rewritten as 
\begin{align}
\su &= \nabla\bu \:\,\,\quad \bx\in \Omega, \;t>0,\label{1.1l}\\
\bchi &= \nabla\theta \:\,\,\quad \bx\in \Omega ,\;t>0,\label{1.1lth}\\
	\bu_t 
		-\nu \nabla \cdot \su  + \bu\cdot \nabla \bu&+ \nabla p+\boldsymbol{\mu}\theta
		=\brf\,\,\,\quad \bx\in \Omega, \;t>0,	\label{1.2l}\\
		\theta_t -k\nabla\cdot\bchi &+ \textbf{u} . \nabla\theta = g\,\,\quad \bx\in \Omega ,\;t>0,\label{1.2lth}\\
	\nabla \cdot \bu&=0\,\,\,\quad\bx\in \Omega,\;t>0.\label{1.3l}
\end{align}	
\noindent 
We now define the variational formulation of (\ref{1.1l})-(\ref{1.3l}) as: For  $ t > 0 $, \\seek $ (\su ,\, \bu,\, \bchi,\, \theta,\, p) \in \underline{\bH}^1_0(\Omega)\times \bH_0^1(\Omega)\times \bH_0^1(\Omega)\times H_0^1(\Omega) \times L^2(\Omega) / \mathbb{R} $,  such that
	\begin{align}
 (\su,\:\underline{\boldsymbol{\upsilon}}) =& (\nabla\bu,\:\underline{\boldsymbol{\upsilon}} ),\label{w0}\\
(\bchi,\, \biota) =& (\nabla\theta,\,\biota ),\label{w0th}\\
	(\bu_t,\,\bphi)+\nu\,(\su,\,\nabla\bphi)+(\bu\cdot\nabla\bu,\,\bphi)-&(p,\,\nabla\cdot\bphi)+(\boldsymbol{\mu}\theta,\,\bphi)=(\brf,\,\bphi),\label{ldgw1}\\
	(\theta_t,\,\psi)+\kappa\,(\bchi,\,\nabla\psi)&+(\bu\cdot\nabla\theta,\,\psi)=(g,\,\psi),\label{w1th}\\
		(\nabla\cdot\bu,\,q)&=0,\label{ldgw2}
	\end{align} 
for all $ (\underline{\boldsymbol{\upsilon}} ,\, \bphi,\, \biota,\, \psi,\, q) \in \underline{\bH}^1_0(\Omega)\times \bH_0^1(\Omega)\times \bH_0^1(\Omega)\times H_0^1(\Omega) \times L^2(\Omega) / \mathbb{R}.$\\	
We recall a set of Sobolev inequalities needed for the nonlinear estimates (see \cite[pp.~291]{T77} and \cite[pp.~11]{T95}).
\begin{lemma}\label{ldgsobolev2.1} For $\mathbf{z} \in \bH_0^1(\Omega)$, there exists a positive constant $K>0$, independent of $h$, such that
\[
\|\mathbf{z}\|_{L^4(\Omega)} \leq 2^{1/4} \|\mathbf{z}\|^{1/2} \|\nabla \mathbf{z}\|^{1/2}.
\]
Additionaly, if $\mathbf{z} \in \bH^2(\Omega)$, then
\[
\|\mathbf{z}\|_{L^\infty(\Omega)} \leq K \|\mathbf{z}\|^{1/2} \|\Delta \mathbf{z}\|^{1/2}.
\] 
\end{lemma}
\noindent
%
\noindent
\noindent
To construct the LDG formulation corresponding to
(\ref{w0})--(\ref{ldgw2}), we consider a quasi-uniform triangulation
$\mathcal{T}_h$ of the closure $\bar{\Omega}$, composed of triangular
elements $T$ with maximal edge length $h$. For each element $T \in
\mathcal{T}_h$, we denote its diameter by $h_T$. Let $\mathcal{N}_h$ denote
the collection of all edges associated with the triangulation
$\mathcal{T}_h$, and for each edge $e \in \mathcal{N}_h$, let $\bn_e$ be a
unit normal vector. Furthermore, let $\mathcal{N}_D \subset \partial\Omega$ be the set of edges
lying on the boundary. On each boundary edge $e \in \mathcal{N}_D$, the
vector $\bn$ is chosen to be the outward unit normal to $\partial\Omega$.
If an interior edge $e$ is shared by two neighboring elements $T_m$ and
$T_n$, we denote by $\bn_m$ and $\bn_n$ the outward unit normal vectors
on $e$ corresponding to $T_m$ and $T_n$, respectively.\\
Let $(\su,\,\bu,\,\bchi,\,\theta,\,q) \in
\Y \times \bX \times \bZ \times Z \times M$ be functions that are smooth
within each element $T_m$ and $T_n$. The average and jump operators for
the scalar variables $p$ and $\theta$, the vector--valued functions
$\bu$ and $\bchi$, and the tensor--valued quantity $\su$ across an edge
$e \in \mathcal{N}_h$ are defined as follows.

\begin{align*}
 \{q\} := \frac{1}{2} (q|_{T_m})|_e + \frac{1}{2} (q|_{T_n})|_e,\quad \{\theta\} := \frac{1}{2} (\theta|_{T_m})|_e + \frac{1}{2} (\theta|_{T_n})|_e,\\ \{\bchi\} := \frac{1}{2} (\bchi|_{T_m})|_e + \frac{1}{2} (\bchi|_{T_n})|_e,\quad \{\bu\} := \frac{1}{2} (\bu|_{T_m})|_e + \frac{1}{2} (\bu|_{T_n})|_e,\\
  \{\su\} := \frac{1}{2} (\su|_{T_m})|_e + \frac{1}{2} (\su|_{T_n})|_e,\quad
\bsq q\bqs := (q|_{T_m})|_e \bn_m+ (q|_{T_n})|_e\bn_n, \\
\bsq \theta \bqs := (\theta|_{T_m})|_e \bn_m+ (\theta|_{T_n})|_e\bn_n,\quad
  [\bchi] := (\bchi|_{T_m})|_e \cdot \bn_m+ (\bchi|_{T_n})|_e\cdot \bn_n,\\
  [\bu] := (\bu|_{T_m})|_e \cdot \bn_m+ (\bu|_{T_n})|_e\cdot \bn_n, \quad \bsq \su \bqs := (\su|_{T_m})|_e \cdot \bn_m+ (\su|_{T_n})|_e\cdot \bn_n. 
\end{align*}
For edges located on the boundary $\partial\Omega$, the average operators
reduce to the trace values of $q$, $\theta$, $\bchi$, $\bu$, and $\su$
on the corresponding edge. The associated jumps are defined by
\(
\bsq q \bqs = q\,\bn, \quad
\bsq \theta \bqs = \theta\,\bn, \quad
[\bchi] = \bchi \cdot \bn, \quad
\bsq \su \bqs = \su \cdot \bn .
\)
We remark that the jumps $\bsq q \bqs$, $\bsq \theta \bqs$, and
$\bsq \su \bqs$ are vector--valued quantities, whereas the jumps
$[\bu]$ and $[\bchi]$ are scalar--valued. In addition, since the velocity  $\bu$ is vector--valued, its jump
across an interior edge is a matrix which is defined defined as
\(
\underline{[\bu]}
:= (\bu|_{T_m})|_e \otimes \bn_m
+ (\bu|_{T_n})|_e \otimes \bn_n .
\)
For boundary edges $e \subset \partial\Omega$, the above definition
reduces to
\(
\underline{[\bu]} = \bu \otimes \bn .
\)\\
To facilitate the analysis that follows, we now specify the functional spaces employed in the LDG formulation
\begin{align*}
    \Y&=\{\su\in \underline{\bL}^2(\Omega):\,\,\, \forall T\in \mathcal{T}_h, \:\:\su|_T\in \underline{\bH}^1(T)\},\\
   \bX &=\{\bphi\in \bL^2(\Omega):\,\,\, \forall T
   \in \mathcal{T}_h, \,\,\, \bphi|_{E}\in \textbf{H}^{2}(T)\},\\
   \bZ &=\{\boldsymbol{\psi}\in \bL^2(\Omega):\,\,\, \forall T\in \mathcal{T}_h, \,\,\, \boldsymbol{\psi}|_{T}\in \textbf{H}^{1}(T)\},\\
   Z &=\{\psi\in L^2(\Omega):\,\,\, \forall T\in \mathcal{T}_h, \,\,\, \psi|_{T}\in {H}^{2}(T)\},\\
   M &= \{q\in L^2(\Omega)/ \R:\,\,\,\forall	T\in \mathcal{T}_h,\,\,\, q|_T\in H^1(T)\}
\end{align*}
      as well as the bilinear forms that are defined over these broken spaces, given by
 $a^l_1 : \Y \times \Y\rightarrow {\R}$, $a^l_2 : \bZ \times \bZ\rightarrow { \R},\,$
$b^l_1 : \bX \times \Y  \rightarrow { \R},\:b^l_2 : Z \times \bZ  \rightarrow { \R},\:
 c^l_1:\bX \times \bX \rightarrow { \R},\, c^l_2:Z \times Z \rightarrow { \R}, 
\,\, d: \bX \times M \rightarrow { \R}  
$  (\cite[pp. 326]{cockburn02}),
\begin{align*}
	 a^l_1(\su,\,\tu)&:=\int\limits_\Omega\su:\tu \,d\bx,\,\quad a^l_2(\bchi,\,\biota):=\int\limits_\Omega\bchi\cdot\biota \,d\bx,\\
b^l_1(\bu,\tu)&:=\sum_{T\in\mathcal{T}_h}\int\limits_T\bu\cdot\nabla\cdot\tu\,d\bx-\sum_{e\in \mathcal{N}_h\backslash \mathcal{N}_D}\int\limits_e(\{\bu\}+\underline{[\bu]}\cdot\textbf{C}_{12})\cdot \bsq\tu\bqs\,ds,\\
b^l_2(\theta,\biota)&:=\sum_{T\in\mathcal{T}_h}\int\limits_T\theta\nabla\cdot\biota\,d\bx-\sum_{e\in \mathcal{N}_h\backslash \mathcal{N}_D}\int\limits_e(\{\theta\}+{[\theta]}\cdot\textbf{E}_{12}) [\biota]\,ds,\\
 c^l_1(\bu,\bphi)&:=\sum_{e\in \mathcal{N}_h}\int\limits_e\frac{c_{11}}{|e|}\underline{[\bu]}:\underline{[\bphi]}\,ds,\quad  c^l_2(\theta,\psi):=\sum_{e\in \mathcal{N}_h}\int\limits_e\frac{c_{22}}{|e|}\bsq\theta\bqs\cdot\bsq\psi\bqs\,ds,\\
d(\bphi,p)&:=-\sum_{T\in\mathcal{T}_h}\int\limits_T p\nabla\cdot\bphi \,d\bx+\sum_{e\in \mathcal{N}_h \backslash \mathcal{N}_D}\int\limits_e(\{p\}-\textbf{D}_{12}\cdot \bsq p \bqs)[\bphi]\,ds+\sum_{e\in \mathcal{N}_D}\int\limits_ep\bphi\cdot \bn\,ds.
	\end{align*}
     Here, $|e|$ denotes the length of the edge $e$, and the positive constants $c_{11}, c_{22} > 0$ are penalty parameters specified for each edge and  
     $\textbf{C}_{12}, \, \textbf{D}_{12}, \, \textbf{E}_{12} \in \R^2$.  Using the integration--by--parts formula, the bilinear forms 
     $b_1^l(\cdot,\cdot)$, $b_2^l(\cdot,\cdot)$, and $d(\cdot,\cdot)$
     can be equivalently expressed in the following alternative forms:
\begin{align*}
	b^l_1(\bu,\tu)&=-\sum_{T\in\mathcal{T}_h}\int\limits_T\nabla\bu:\tu \,d\bx +\sum_{e\in \mathcal{N}_h\backslash \mathcal{N}_D}\int\limits_e(\{\tu\}-\bsq \tu \bqs \otimes\textbf{C}_{12}):\underline{[\bu]} \,ds+\sum_{e\in \mathcal{N}_D}\int\limits_e \tu:(\bu\otimes\bn)\,ds,\\
	b^l_2(\theta,\biota)&=-\sum_{T\in\mathcal{T}_h}\int\limits_T\nabla\theta\cdot\biota \,d\bx +\sum_{e\in \mathcal{N}_h\backslash \mathcal{N}_D}\int\limits_e(\{\biota\}-[\biota] { \textbf{E}_{12}})\cdot{\bsq\theta\bqs} \,ds+\sum_{e\in \mathcal{N}_D}\int\limits_e \biota\cdot\theta\bn\,ds,\\
	d(\bphi,p)&=\sum_{T\in\mathcal{T}_h}\bphi\cdot\nabla p\, d\bx - \sum_{e\in \mathcal{N}_h\backslash \mathcal{N}_D}\int\limits_e(\{\bphi\}+\textbf{D}_{12}[\bphi])\cdot \bsq p \bqs\, ds. 
\end{align*}
\\
   We now define the seminorm (see \cite[pp.~328]{cockburn02}) and the associated energy norm for the spaces $\Y$ and $\bX$, as well as for $\bZ$ and $Z$. In addition, we introduce the norm on the space $M$, which are specified below
\begin{align*}
		|(\tu, \bv)|^2_\mathcal{A} 
		&:=a_1^l(\tu,\tu)+c_1^l(\bv,\bv) \quad \forall \tu\in \Y, \quad \forall \bv\in \bX,\\
		|(\bchi, {\psi})|^2_\mathcal{A} 
		&:=a^l_2(\bchi,\bchi)+c_2^l({\psi},{\psi}) \quad \forall \bchi\in \bZ, \quad \forall {\psi}\in Z,\\
\|\bphi\|_\varepsilon 
		&=\bigg(\sum_{T\in \mathcal{T}_h} \|\nabla\bphi\|^2_{L^2(E)}+c_1^l(\bphi,\bphi)\bigg)^{1/2} \quad \forall \bphi\in \bX,\\
		\|r\|_M &= \|r\|_{L^2(\Omega)/\R} \quad  \forall r\in M,\\
		\|\psi\|_\varepsilon 
		&=\bigg(\sum_{T\in \mathcal{T}_h} \|\nabla\psi\|^2_{L^2(E)}+c_2^l(\psi,\psi)\bigg)^{1/2} \quad \forall \psi\in Z.
		\end{align*}
For the purposes of the subsequent analysis, we state the classical $L^p$--estimate for elements of the broken spaces $\bX$ and $Z$ (\cite[Lemma~6.2]{GRW005}):
	\begin{align}
\|\bz\|_{L^p(\Omega)}\leq \mathcal{N}_\bu \|\bz\|_\varepsilon\quad  \forall \bz \in \bX,\: 2\le p < \infty,\label{ldgLp}\\
\|z\|_{L^p(\Omega)}\leq \mathcal{N}_\theta \|z\|_\varepsilon\quad  \forall z \in Z,\: 2\le p < \infty, \label{ldgLpth}
	\end{align} 
    where $\mathcal{N}_\bu$ and $\mathcal{N}_\theta=C(p)>0$, depends on $p$.\\
We now recall several trace inequalities that will be utilized throughout the analysis for the spaces $\bX$ and $\Y$ (see \cite[Section~2.1.3]{R08} and \cite[Section~7]{OK98}).
\begin{lemma}
\label{ldgtrace}
For every element \(T \in \mathcal{T}_h\), 
\begin{align*}
		\|\bw\|_{0,e} &  \leq C(h_T^{-1/2}\|\bw\|_{0,T}+h_T^{1/2}\|\nabla\bw\|_{0,T})\quad \forall e\in \partial T,\,\, \forall\bw\in \bX,\\
		\|\nabla\bw\|_{0,e} & \leq C(h_T^{-1/2}\|\nabla\bw\|_{0,T}+h_T^{1/2}\|\nabla^2 \bw\|_{0,T})\quad \forall e\in \partial T,\,\,\forall \bw\in \bX,\\
		\|\bw\|_{L^4(e)} & \leq Ch_T^{-3/4}(\|\bw\|_{0,T}+h_T\|\nabla \bw\|_{0,T})\quad \forall e\in \partial T,\,\, \bw\in \bX,\\
        \|\su\|_{0,e} & \leq C(h_T^{-1/2}\|\su\|_{0,T}+h_T^{1/2}\|\nabla \su\|_{0,T})\quad \forall e\in \partial T,\,\,\forall \su\in \Y.
	\end{align*}
\end{lemma}
\noindent
It is worth noting that Lemma \ref{ldgtrace} also applies when the vector--valued function $\bw$ is substituted with a scalar--valued function $\psi$.\\
To formulate the LDG discretization of the nonlinear terms in (\ref{1.1})--(\ref{1.3}), we express the trilinear form $\mo^\by(\cdot, \cdot, \cdot)$ as follows:
\begin{align*}
		\mo_1^{\bw}&(\bxi,\boldsymbol{\eta}, \boldsymbol{\varphi})=\sum_{T\in\mathcal{T}_h}  \left(\int\limits_T(\bxi\cdot \nabla \boldsymbol{\eta}) \cdot \boldsymbol{\eta}\,d\bx+{\int\limits_{\partial T_-}|\{\bxi\}\cdot \bn_T|(\boldsymbol{\eta}^{int}-\boldsymbol{\eta}^{ext})\cdot \boldsymbol{\eta}^{int}} d\bx\,\right)\nonumber\\
		& +\frac{1}{2}\sum_{T\in \mathcal{T}_h}\int\limits_T(\nabla\cdot \bxi)\boldsymbol{\eta}\cdot \boldsymbol{\eta} d\bx \,)-\frac{1}{2}\sum_{e\in \mathcal{N}_h}\int\limits_e [ \bxi] \{\boldsymbol{\eta}\cdot \boldsymbol{\eta}\} ds \, \forall \by,\,\bxi,\,\boldsymbol{\eta},\,\boldsymbol{\varphi}\in \textbf{V},
	\end{align*} 
	where 
	\[\partial T_-=\{\textbf{x}\in \partial T:\{\bw\}\cdot \bn_T<0\}.\]
 The dependence of $\partial T_-$ on the direction $\bw$ is indicated by the parameter $\bw$. The notation \textit{int} and \textit{ext} refer, respectively, to the trace of a function taken from the interior and exterior of the element $T$ on a given side. For sides of $T$ that lie on the boundary $\partial \Omega$, the exterior trace is taken to be zero. In the definition of $\mo_1(\cdot; \cdot, \cdot, \cdot)$, the first two terms were originally introduced in \cite{LR74} in the context of transport equations, while the third term is included to guarantee the positivity property of $\mo_1$ (see (\ref{ldg5.5}) below). The form $\mo_1(\cdot; \cdot, \cdot, \cdot)$ is motivated by the Lesaint--Raviart upwinding scheme \cite{LR74} and was adapted in \cite{GRW005} for the treatment of the Navier--Stokes equations.


    \noindent
By employing the integration--by--parts formula, the trilinear form can be rewritten as follows:
	{ \begin{align}
		&\mo_1^{\bxi}(\bxi,\,\boldsymbol{\eta},\,\boldsymbol{\varphi})= -\sum_{T\in \mathcal{T}_h}\left(\int\limits_T (\bxi \cdot \nabla \boldsymbol{\varphi})\cdot \boldsymbol{\eta}\,d\bx+\frac{1}{2}\int\limits_T(\nabla\cdot \bxi)\boldsymbol{\eta}\cdot\boldsymbol{\varphi}\,d\bx\right)+\int\limits_{\mathcal{N}_+}|\bxi\cdot \bn|\boldsymbol{\eta}\cdot \boldsymbol{\varphi}\,ds
		\nonumber\\
		&\hspace{1cm} +\frac{1}{2}\sum_{e\in \mathcal{N}_h}\int\limits_e [\bxi] \{\boldsymbol{\eta}\cdot\boldsymbol{\varphi}\}\,ds -\sum_{T\in \mathcal{T}_h}\int\limits_{\partial T_-}|\{\bxi\}\cdot \bn_T|\boldsymbol{\eta}^{ext}\cdot (\boldsymbol{\varphi}^{int}-\boldsymbol{\varphi}^{ext})\,ds,\label{ldgintegration}
	\end{align}}
	where $\mathcal{N}_+$ is the subset of $\partial\Omega$, such that $\bxi \cdot \bn > 0$ and $\bxi, \boldsymbol{\eta}\in \mathbf{V}$. 
\\
Selecting $\boldsymbol{\varphi} = \boldsymbol{\eta}$ in equation~(\ref{ldgintegration}) yields
\begin{equation}
		\mo_1^{\bxi}(\bxi,\,\boldsymbol{\eta},\,\boldsymbol{\eta})\geq 0, \,\,\, \bxi,\,\boldsymbol{\eta} \in\bX.\label{ldg5.5}
	\end{equation}
\noindent
We now formulate the LDG weak formulation corresponding to 
\eqref{w0}--\eqref{ldgw2}. 
For any fixed time $t>0$, determine 
$(\su,\, \bu,\, \bchi,\, \theta,\, p) \in \Y \times \bX \times \bZ \times Z \times M$
such that
\begin{align}
	a^l_1(\su,\tu)+b^l_1(\bu,\tu)&=0, \label{2.we1}\\
	a^l_2(\bchi,\biota)+b^l_2(\theta,\biota)&=0, \label{2.we1th}\\
	(\bu_t,\bz)-\nu(b_1^l(\bz,\su)-c_1^l(\bu,\bz))
		+\mo_1^\bu(\bu,\,\bu,\,&\bz)+d(\bz,p)+(\boldsymbol{\mu}\theta,\,\bz)= (\textbf{f},\,\bz), \label{2.we2}\\
		(\theta_t,\psi)-\kappa(b^l_2(\psi_,\bchi)-c^l_2(\theta,\psi))
		&+\mo_2^\theta(\bu,\,\theta,\,\psi)= (g,\,\psi), \label{2.we2th}\\
		 d(\bu,q)= 0,& \label{2.we3}
	\end{align}
	for all $ (\tu,\, \bz,\,\biota,\,\psi,\, q) \in \Y \times \bX\times\bZ\times Z \times M.$\\
To facilitate the error estimates in the later sections, 
we assume that the solution 
$(\su,\, \bu,\, \bchi,\, \theta,\, p)$ of 
\eqref{2.we1}--\eqref{2.we3} satisfies the following condition 
throughout the paper:
\\
{\bf (A2)}.  The solution  $(\su,\, \bu,\,\bchi,\,\theta, p)$ of (\ref{2.we1})-(\ref{2.we3}) satisfy 
\noindent	
{ \[
\begin{array}{ll}
\text{(i)} &
\left\{
\begin{aligned}
\su &\in L^\infty(0,T;\underline{\bH}^{r+1}(\Omega)) \cap L^2(0,T;\underline{\bH}^{r+1}(\Omega)), \\
\su_t &\in L^2(0,T;\underline{\bH}^{r+1}(\Omega)) \cap L^\infty(0,T;\underline{\bH}^{r}(\Omega)),
\end{aligned}
\right. \\[1.5ex]

\text{(ii)} &
\left\{
\begin{aligned}
\bu &\in L^\infty(0,T;\bH^{r+1}(\Omega)) \cap L^2(0,T;\bH^{r+1}(\Omega)), \\
\bu_t &\in L^2(0,T;\bH^{r+1}(\Omega)) \cap L^\infty(0,T;\bH^{r}(\Omega)),
\end{aligned}
\right. \\[1.5ex]

\text{(iii)} &
\left\{
\begin{aligned}
\bchi &\in L^\infty(0,T;\bH^{r+1}(\Omega)) \cap L^2(0,T;\bH^{r+1}(\Omega)), \\
\bchi_t &\in L^2(0,T;\bH^{r+1}(\Omega)) \cap L^\infty(0,T;\bH^{r}(\Omega)),
\end{aligned}
\right. \\[1.5ex]

\text{(iv)} &
\left\{
\begin{aligned}
\theta &\in L^\infty(0,T;H^{r+1}(\Omega)) \cap L^2(0,T;H^{r+1}(\Omega)), \\
\theta_t &\in L^2(0,T;H^{r+1}(\Omega)) \cap L^\infty(0,T;H^{r}(\Omega)),
\end{aligned}
\right. \\[1.5ex]

\text{(v)} &
\left\{
\begin{aligned}
p &\in L^\infty(0,T;H^{r}(\Omega)) \cap L^2(0,T;H^{r}(\Omega)), \\
p_t &\in L^2(0,T;H^{r}(\Omega)).
\end{aligned}
\right.
\end{array}
\]
}
\noindent	
   We proceed to analyze  {\it a priori} error estimates associated with the semidiscrete LDG scheme.
	\section{LDG Semidiscrete  Formulation}\label{ldgsection3}
	\se
We begin by introducing the finite--dimensional discontinuous spaces  $ \Y_h \subset \Y $, $ \bX_h \subset \bX $, $ \bZ_h \subset \bZ $, $ Z_h \subset Z $, and $ M_h \subset M $, which are used to approximate the semidiscrete solutions. We define them as follows:
\begin{align*}
\Y_h&=\{ \su_h\in  \underline{\bL}^2(\Omega)\, : \,\,\, \su_h|_T \in (\mathbb{P}_{r-1}(T))^4 \,\,\, \forall T \in \mathcal{T}_h\},\\
\bX_h  &= \{ \textbf{w}_h\in  \bL^2(\Omega)\, : \,\,\, \textbf{w}_h|_T \in (\mathbb{P}_r(T))^2 \,\,\, \forall T \in \mathcal{T}_h\},\\
\bZ_h  &= \{ \boldsymbol{\psi}_h\in  \bL^2(\Omega)\, : \,\,\, \boldsymbol{\psi}_h|_T \in (\mathbb{P}_{r-1}(T))^2 \,\,\, \forall T \in \mathcal{T}_h\},\\
Z_h  &= \{ {\psi}_h\in  \bL^2(\Omega)\, : \,\,\, {\psi}_h|_T \in \mathbb{P}_{r}(T) \,\,\, \forall T \in \mathcal{T}_h\},\\
	M_h & = \{q_h\in  L^2(\Omega)/ \R\, :\,\,\, q_h|_T\in \mathbb{P}_{r-1}(T) \,\,\, \forall T\in \mathcal{T}_h\},
\end{align*}
where $ r \ge 1 $ is an integer, and $ \mathbb{P}_r(E) $ denotes the space of polynomials of degree at most $ r $ on $ T $.
\begin{lemma}\cite[pp. 61]{GRW005}\label{Rhp}
	 For every \( T \in \mathcal{T}_h \), there exists a projection operator \( \textbf{R}_h \in \mathcal{R}(\bH^1(\Omega); \bX_h) \), such that
	\begin{align*}
\forall s \in [1, r+1],~\forall \bz\in  \bH^s(\Omega)\cap\bH_0^1(\Omega),~~ \|\bR_h(\bz)-\bv\|_\varepsilon&\leq C h^{s-1}|\bz|_{s,\Omega},\\
		 \forall \bz\in \bH^s(\Omega),~~ \|\bz-\textbf{R}_h(\bz)\|_{L^2(T)}&\leq C h_T^{s}|\bz|_{s,\,\Delta_T}
	\end{align*}
	where, \( \Delta_T \) denotes a macro--element that contains \( T \).
	\end{lemma}
	\noindent
Moreover, the operator \( \bR_h \) (see \cite{GRW005}) satisfies    
\begin{equation*}
	 d(\textbf{R}_h(\bz) - \bz,\, p_h) = 0 \qquad \forall \bz \in \bH_0^1(\Omega)^2 \quad \forall p_h \in M_h.
\end{equation*}  
In addition, $\mathbf{R}_h$ is stable with respect to the energy norm
(see \cite{GRW05}); that is, there exists a constant $K>0$,
independent of  $h$, such that 
\begin{equation*}
	\|\bR_h \bz\|_\varepsilon \leq C |\bz|_{1,\,\Omega} \quad  \forall \bz \in \bH_0^1(\Omega).  
\end{equation*}
Furthermore, for every \( T \in \mathcal{T}_h \), there exists an operator \( r_h \in \mathcal{R}(L^2(\Omega)/ \mathbb{R}; M_h) \) (see \cite{GRW005}) such that 
 \begin{align}
		 \int\limits_T  y_h(r_h(p)-p) \,d\bx\,=0\quad \forall p\in L^2(\Omega)/ \R,~~\forall y_h\in \mathbb{P}_{r-1}(T), \label{ldg2}\\
	\|p-r_h(p)\|_{0,T}\leq C h_T^s|p|_{s,T}\quad \forall s\in[0,r],\quad	\forall p\in H^s(\Omega)\cap L^2(\Omega)/ \R. \label{ldg2.1}
	\end{align}
Now, we state the semidiscrete LDG formulation for  (\ref{2.we1})--(\ref{2.we3}) as follows: \\
Seek $(\su_h,\, \bu_h,\,\bchi_h,\,\theta_h,\,p_h)\,\in\, \Y_h\times\bX_h\times\bZ_h\times Z_h \times M_h$, $t>0$ such that 
\begin{align}
	a^l_1(\su_h,\tu_h)&+b^l_1(\bu_h,\tu_h)=0, \label{3.w1}\\
	a^l_2(\bchi_h,\biota_h)&+b^l_2(\theta_h,\biota_h)=0, \label{3.w1th}\\
	(\bu_{ht},\bphi_h)-\nu(b^l_1(\bphi_h,\su_h)-c^l_1(\bu_h,\bphi_h))
		&+\mo_1^{\bu_h}(\bu_h,\bu_h,\,\bphi_h)+d\,(\bphi_h,\,p_h)
		+(\boldsymbol{\mu} \theta_h,\,\bphi_h)= (\textbf{f},\bphi_h), \label{3.w2}\\
        \theta_{ht},\psi_h)-\kappa(b^l_2(\psi_h,\bchi_h)-c^l_2(\theta_h,\psi_h))
		&+\mo_2^{\theta_h}(\bu_h,\theta_h,\,\psi_h)
        = (g,\,\psi_h), \label{3.w2th}\\
		 d\,(\bu_h,\,q_h)&= 0, \label{3.w3}
\end{align}
for all $(\tu_h,\,\bphi_h,\,\biota_h,\,\psi_h,\,\,q_h)\in (\Y_h,\,\bX_h,\,\bZ_h,\,Z_h,\,M_h)$.\\
We also define the discontinuous space $\bX^l_h$, which is analogous to $\bJ_1$, as follows:
 \[\bX^l_h\,=\,\{\bv_h\in \bX_h:\,\, d\,(\bv_h,\,q_h)=0\,\, \forall q_h\in M_h\}.
	\]
Equivalently, LDG scheme corresponding to  (\ref{3.w1})--(\ref{3.w3}) on the spaces $\Y_h$, $\bX^l_h$, $\bZ_h$ and $Z_h$ is given as follows:
Seek $(\su_h, \bu_h,\,\bchi_h,\,\theta_h)\in (\Y_h,\,\bX^l_h,\,\bZ_h,\,Z_h)$, $t>0$, such that
	\begin{align}
	a_1^l(\su_h,\tu_h)&+b_1^l(\bu_h,\tu_h)=0, \label{3.we1}\\
	a^l_2(\bchi_h,\biota_h)&+b^l_2(\theta_h,\biota_h)=0, \label{3.we1th}\\
	(\bu_{ht},\,\bphi_h)-\nu(b_1^l(\bphi_h,\su_h)-c_1^l(\bu_h,\bphi_h))
		&+\mo_1^{\bu_h}(\bu_h,\,\bu_h,\,\bphi_h)+(\boldsymbol{\mu} \theta_h,\,\bphi_h)
		= (\textbf{f},\,\bphi_h), \label{3.we2}\\
	 (\theta_{ht},\,\psi_h)-\kappa(b^l_2(\psi_h,\bchi_h)-c^l_2(\theta_h,\psi_h))
		&+\mo_2^{\theta_h}(\bu_h,\theta_h,\,\psi_h)
        = (g,\,\psi_h), \label{3.we2th}
	\end{align}
for all $(\tu_h,\,\bphi_h,\,\biota_h,\,\psi_h)\in (\Y_h,\,\bX_h^l,\, \bZ_h,\,Z_h)$.\\
\noindent
To establish optimal error estimates for the semidiscrete LDG approximations of the stress tensor, velocity, pressure, temperature gradient, and temperature, the analysis is restricted to weakly divergence--free spaces. Accordingly, an approximation operator defined on $\bJ_1 \cap \bH^{r+1}(\Omega)$ is introduced.
\begin{lemma}\label{ldgih}   
		For each \( \bz \in \bJ_1 \cap \bH^{r+1}(\Omega) \),  there is an approximation $i_h:\bJ_1 \cap \bH^{r+1}(\Omega)\rightarrow \bX_h^l$ such that
\begin{equation*}
		\|\bz-i_h \bz\|+h	\|\bz-i_h \bz\|_\varepsilon\leq Ch^{r+1}\|\bz\|_{r+1},
		\end{equation*}
		 where $i_h$ is the restriction of $\bR_h$ on $\bJ_1 \cap \bH^{r+1}(\Omega)$.
	\end{lemma}
	\noindent
	\noindent
Now we state   boundedness property of the bilinear forms \(a_1^l(\cdot, \cdot)+ c_1^l(\cdot, \cdot)\) and \(a_2^l(\cdot, \cdot)+ c_2^l(\cdot, \cdot)\).
The result follows by applying the Cauchy–Schwarz inequality to the bilinear forms and estimating each term in terms of the associated energy norms. Squaring the resulting bound and using Young’s inequality yields the desired continuity estimates with a constant independent of $h$.
	\begin{lemma}
    \label{ldgcont}
		There exists a constant \( C > 0 \), independent of \( h \), such that  
\begin{align*}
			|a_1^l(\su_h,\tu_h) +c_1^l(\bv_h,\bw_h )|\leq & C |(\su_h,\bv_h)|_\mathcal{A}|(\tu_h,\bw_h)|_\mathcal{A},\quad
			|a_2^l(\bchi_h,\bpsi_h) +c_2^l(\theta_h,\psi_h )|\leq & C |(\bchi_h,\,\theta_h)|_\mathcal{A}|(\bpsi_h,\psi_h)|_\mathcal{A},
		\end{align*}
		for all \( \su_h,\, \tu_h \in \Y_h\),   \(\bv_h,\, \bw_h \in \bX_h \),  $\bchi_h,\,\bpsi_h \in \bZ_h$, $\theta_h\,\text{and}\,\psi_h\in Z$.
	\end{lemma}
	\noindent
Analogous to the continuous framework, 
$(\widetilde{\bV}_h, M_h)$ satisfies a discrete inf--sup stability condition as 
\[ \widetilde{\textbf{V}}_h=\bigg\{\bw_h\in \bX_h: \forall e \in \mathcal{N}_h, \quad \int\limits_e \brho_h\cdot [\bw_h]\bn_e\,=0\,ds\quad \forall \brho_h\in (\mathbb{P}_{k-1}(e))^2\bigg\}, \]
	which plays a crucial role in the derivation the pressure estimates.
	\begin{lemma}\cite[Theorem 4.5]{GRW005}\label{ldginf-sup}
		There exists a constant \( \beta^* > 0 \), independent of \( h \) such that \begin{equation*}
			\inf_{p_h\in M_h}\sup_{\bv_h\in \widetilde{\textbf{V}}_h}\frac{d(\bv_h,p_h)}{\|\bv_h\|_\varepsilon \|p_h\|}\geq \beta ^*.
		\end{equation*}
	\end{lemma}
    \noindent
In analogy with the initial bound presented in Lemma \ref{ldgsobolev2.1}, 
an analogous $L^4$--estimate can also be established for functions belonging to $\bX_h$.\begin{lemma}\cite[Theorem 3.8]{GRW05}\label{ldgl4}
Let $\Omega \subset \mathbb{R}^2$ be convex. 
Then a constant $K>0$, independent of $h$, such that
\begin{align*}
\|\bw_h\|_{L^4(\Omega)}\leq C \|\bw_h\|^{1/2}\|\bw_h\|^{1/2}_\varepsilon \quad \forall\bw_h\in \bX_h.
\end{align*}
\end{lemma}

	\noindent
We state in Lemma \ref{ldgdistrace} several trace and inverse 
estimates for the finite element spaces $\bX_h$ and $\Y_h$, analogous to the inequalities given in Lemma \ref{ldgtrace} 
for $\bX$ and $\Y$.
\begin{lemma}\cite[Section 1.4.3]{DE12}\label{ldgdistrace}
	There exists a  positive constant $C$, independent of $h_T$, $\bw_h$ and $\su_h$, such that 
\begin{align*}
		\|\bw_h\|_{0,e}&\leq Ch_T^{-1/2}\|\bw_h\|_{0,T}\quad \forall e\in \partial T, \,\,\forall\: \bw_h\in\bX_h,\nonumber\\
		\|\nabla\bw_h\|_{0,e}&\leq Ch_T^{-1/2}\|\nabla\bw_h\|_{0,T}\quad \forall e\in \partial T, \,\,\forall \bw_h\in\bX_h, \nonumber\\
	\|\nabla\bw_h\|_{0,T}&\leq Ch_T^{-1}\|\bw_h\|_{0,T}\quad \forall \bw_h\in\bX_h,\\
		\|\bw_h\|_{L^4(T)}&\leq Ch_T^{-1/2}\|\bw_h\|_{0,T}\quad \forall \bw_h\in\bX_h, \nonumber\\
        \|\su_h\|_{0,e}&\leq Ch_T^{-1/2}\|\su_h\|_{0,T}\quad \forall e\in \partial T, \,\,\forall \su_h\in \Y_h,
	\end{align*}
	
\end{lemma}
\noindent
The next  lemma characterizes the relationship between $\|\cdot\|_\varepsilon$--norm and $|\cdot|_\mathcal{A}$--seminorm.
The estimates follow from (\ref{3.we1}) by applying the Cauchy–Schwarz inequality and trace bounds to control $\| \bu_h\|_\varepsilon$ in terms of $|(\su_h,\bu_h)|_{\mathcal A}$ and $\| \theta_h\|_\varepsilon$ in terms of $|(\bchi_h,\theta_h)|_{\mathcal A}$. 
The remaining inequalities are derived analogously by suitable choices of test functions and similar arguments.
\begin{lemma}\label{Enrscript}
    Let $(\su_h, \bu_h)$ and $(\bchi_h,\,\theta_h)$ be the solutions of (\ref{3.we1}) and (\ref{3.we1th}), respectively. Then, 
    \begin{align*}
        \|\bu_h\|^2_\varepsilon\leq C_\bu^1 |(\su_h,\bu_h)|^2_\mathcal{A},\quad \|\theta_h\|^2_\varepsilon\leq C_\theta^1 |(\bchi_h,\theta_h)|^2_\mathcal{A},\quad
        |(\su_h,\bu_h)|^2_\mathcal{A} \leq C_\bu^2\|\bu_h\|^2_\varepsilon,\quad  |(\bchi_h,\theta_h)|^2_\mathcal{A} \leq C_\theta^2 \|\theta_h\|^2_\varepsilon, 
    \end{align*}
    where $C_\bu^1$, $C_\bu^2$  and  $C_\theta^1$, $C_\theta^2$ are positive   constants, independent of $h$.
\end{lemma}
\noindent
We now state bounds for the trilinear form ${\mo_1}(\cdot;\cdot,\cdot,\cdot)$. 
The estimates given in (\ref{ldgtri1}) and (\ref{ldgtri2}) follow from 
\cite[Proposition 4.1]{GRW05} together with Lemma \ref{ldgsobolev2.1}. 
The inequality (\ref{ldgtri20}) is obtained by combining the Cauchy–Schwarz 
and Young inequalities with the results of Lemmas \ref{ldgtrace} and \ref{ldgdistrace}. 
Finally, estimate (\ref{ldgtri200}) can be derived by arguments analogous 
to those used in \cite[Proposition 4.10]{GR09} and Lemma \ref{ldgl4}.
	\begin{lemma} \label{ldgTrilinear}
	(i) For all $\bphi \in \bW^{1,4}(\Omega)$ and $ \bpsi_h,\, \bv_h,\,\bw_h\in\bX_h,$\,
     \begin{align}
|\mo_1^{\bpsi_h}(\bv_h,\bphi,\bw_h)|\leq C\|\bv_h\|\|\bphi\|^{1/2}\|\|\bphi\|^{1/2}_\varepsilon\|\bw_h\|_\varepsilon \label{ldgtri1}.
		\end{align}
		(ii) Let $\bphi \in \bX$, then 
   \begin{align}
			|\mo_1^{\bphi}(\bv_h,\bw_h,\bu_h)|\leq C\|\bu_h\|_\varepsilon\|\bv_h\|_\varepsilon\|\bw_h\|_\varepsilon\quad \forall \bv_h, \bw_h, \bu_h \in \bX_h.\label{ldgtri2}
		\end{align} 
		(iii) Let $\bphi, \bv, \bu \in \bX$ and $\brho_h \in \bX_h$, then
     \begin{align}
			|\mo_1^{\bphi}(\bv,\bu,\brho_h)|\leq C\|\bv\|_\varepsilon\|\bu\|_\varepsilon\|\brho_h\|_\varepsilon.\label{ldgtri20}
		\end{align}
        (iv)  Let $ \brho_h, \bpsi_h, \bv_h, \bw_h \in \bX_h$, then
     \begin{eqnarray}
			|\mo_1^{\brho_h}(\bpsi_h,\bv_h,\bw_h)-\mo_1^{\bphi_h}(\bpsi_h,\bv_h,\bw_h)|\leq C\|\brho_h-\bpsi_h\|^{1/2}\|\brho_h-\bpsi_h\|^{1/2}_\varepsilon\|\bv_h\|_\varepsilon\|\bw_h\|_\varepsilon.\label{ldgtri200}
		\end{eqnarray}
		\end{lemma}
	\noindent
\noindent
We next establish regularity estimates for $\su_h$, $\bu_h$, $\bchi_h$, and $\theta_h$, which play a crucial role in the derivation of the error bounds.
\begin{lemma}\label{3.uap0}
Assume  $\textbf{(A2)}$ is satisfied. Let $0<\al \leq \min \{\frac{\kappa}{2(\mathcal{N}_\theta C_\theta^1)},\,\frac{\nu}{2(\mathcal{N}_\bu C_\bu^1)}\}$. Then, for  $0<t\leq T$
\begin{align*}
\sup_{0<t<\infty}\{\|\bu_h(t)\|+\|\theta_h(t)\|\}+e^{-2\alpha t}\int\limits_0^te^{2\al s}(|(\su_h(s),\,\bu_h(s))|^2_\mathcal{A}+|(\bchi_h(s),\,\theta_h(s))|^2_\mathcal{A}){ ds}\leq C(t),\\
\|\bu_{ht}(t)\|^2+\|\theta_{ht}(t)\|^2+e^{-2\alpha t}\int\limits_0^te^{2\al s}(|(\su_{hs}(s),\,\bu_{hs}(s)|^2_\mathcal{A}+|(\bchi_{hs}(s),\,\theta_{hs}(s)|^2_\mathcal{A})ds\leq C(t),
\end{align*}
furthermore, 
\begin{align*}
\limsup_{t\rightarrow \infty}|(\bchi_h(t),\,\theta_h(t)|_\mathcal{A} \leq \frac{\mathcal{N}_\theta C^1_\theta \|g\|_{L^\infty(\bL^2(\Omega))}}{\kappa},
\end{align*}
where  $C(t)>0$ is a positive constant.
\end{lemma}
\begin{proof}
Choose $\biota_h=\kappa \bchi_h$  in (\ref{3.we1th}) and $\bpsi_h=\theta_h$ in (\ref{3.we2th}). By summing these equations and applying  \eqref{ldg5.5} together with 
Lemma~\ref{Enrscript}, we obtain 
\begin{align}\label{3.uap5th}
\frac{1}{2}\frac{d}{dt}\|\theta_h\|^2+\kappa |(\bchi_h,\,\theta_h)|^2_\mathcal{A}\leq \|g\|\|\theta_h\|\leq \frac{\kappa}{2}|(\bchi_h,\,\theta_h)|^2_\mathcal{A}+ \frac{\mathcal{N}_\theta^2(C_\theta^1)^2}{2\kappa}\|g\|^2.
\end{align}
Multiplying \eqref{3.uap5th} by $e^{2\alpha t}$, 
then using \eqref{ldgLpth} together with Lemma~\ref{Enrscript} 
and integrating from  $0$ to $t$, we arrive at
\begin{align}
 \label{3.uap6}
e^{2\al t}\|\theta_h\|^2\,+(\kappa-2\al (\mathcal{N}_\theta C_\theta^1))\int\limits_0^te^{2\al s} |(\bchi_h(s),\,\theta_h(s))|^2_\mathcal{A} ds\leq \|\theta_h(0)\|^2
+\frac{\mathcal{N}^2(C_\theta^1)^2}{2\nu}\int\limits_0^t e^{2\al s} \|g(s)\|^2ds.
\end{align}
\noindent
We multiply \eqref{3.uap6} by $e^{-2\alpha t}$ 
and employ   
\begin{align}\label{3.uap60}
e^{-2 \al t}\int\limits_0^t e^{-2 \al s}ds
=\frac{1}{2\al}(1-e^{-2\al t})
\end{align}
and choose $0<\al <\frac{\kappa}{2\mathcal{N}_\theta C_\theta^1}$ to obtain
\begin{align}\label{3.uap7}
\|\theta_h(t)\|^2+e^{-2\al t}\int\limits_0^te^{2\al s} |(\bchi_h(s),\,\theta_h(s))|^2_\mathcal{A} ds\leq C.
\end{align}
Multiplying $e^{2\alpha t}$ to 
\eqref{3.uap5th} and integrating  
 from $0$ to $t$, we deduce
 \begin{align*}
e^{2\al t}\|\theta_h(t)\|^2&+\kappa \int\limits_0^te^{2\al s} |(\bchi_h(s),\,\theta_h(s))|^2_\mathcal{A} ds\leq \|\theta_h(0)\|^2+2\al \int\limits_0^te^{2\al s}(\|\theta_h(s)\|^2
+(e^{2\al t}-1)\frac{\mathcal{N}_\theta^2(C_\theta^1)^2\|g\|^2_{L^\infty(L^2(\Omega))}}{2 \al \kappa}.
\end{align*}
We multiply the above equation by $e^{-2\alpha t}$, 
consider the $\limsup$ as $t \to \infty$, and use 
L'H\^opital's rule. Consequently, we arrive at
\begin{align}\label{3.uap9}
\frac{\kappa}{2\al} \limsup_{t\rightarrow \infty}|(\bchi_h(t),\,\theta_h(t))|^2_\mathcal{A}\leq \frac{\mathcal{N}_\theta^2(C_\theta^1)^2\|g\|^2_{L^\infty(L^2(\Omega))}}{2 \al \kappa}.
\end{align}
Furthermore, 
substituting $\tu_h=\nu \su_h$  in (\ref{3.we1}) and $\bphi_h=\bu_h$ in (\ref{3.we2}), applying the similar analysis as used to obtain (\ref{3.uap7}) and (\ref{3.uap9}), we arrive at
\begin{align}\label{3.uap7vel}
&\hspace{2cm}\|\bu_h(t)\|^2+e^{-2\al t}\int\limits_0^te^{2\al s} |(\su_h(s),\,\bu_h(s))|^2_\mathcal{A} ds\leq C,\\
&\frac{\nu}{2\al} \limsup_{t\rightarrow \infty}|(\su_h(t),\,\bu_h(t))|^2_\mathcal{A}\leq \frac{\mathcal{N}_\bu^2 (C_\bu^1)^2\|\brf\|^2_{L^\infty(\bL^2(\Omega))}}{2 \al \nu}+ C\frac{\mathcal{N}_\theta^2(C_\theta^1)^2\|g\|^2_{L^\infty(L^2(\Omega))}}{2 \al \kappa}.\nonumber
\end{align}

\noindent
To derive the subsequent estimates, we differentiate 
\eqref{3.we1th} and \eqref{3.we2th} with respect to time, 
which yields
\begin{align}
a_2^l(\bchi_{ht},\,\biota_h)&+b_2^l(\theta_{ht},\,\biota_h)=0 \quad \quad \quad \forall \biota_h\in  \bZ_h, \label{3.uap10th}\\
(\theta_{htt},\,\psi_h)-\kappa(b_2^l(\psi_h,\,\bchi_{ht})-c_2^l(\theta_{ht},\,\psi_h))&+\mo_2^{\theta_h}(\bu_{ht},\,\theta_h,\,\psi_h)+\mo_2^{\theta_h}(\bu_h,\,\theta_{ht},\,\psi_h)= (g_t,\,\psi_h)\; \forall\psi_h \in Z_h^l. \label{3.uap11th}
	\end{align}
	Substitute  $\tu_h=\kappa\bchi_{ht}$  in (\ref{3.uap10th}),  and $\psi_h=\theta_{ht}$ in (\ref{3.uap11th}). Then, summing the resulting equations and employing the definition 
of the $|\cdot|_{\mathcal{A}}$--seminorm, together with 
\eqref{ldgLpth}, \eqref{ldg5.5}, and Lemma~\ref{Enrscript}, 
we obtain
\begin{align}\label{3.uap12}
\frac{1}{2}\frac{d}{dt}\|\theta_{ht}\|^2+\kappa|(\bchi_{ht},\,\theta_{ht})|^2_\mathcal{A} \leq --\mo_2^{\theta_h}(\bu_{ht},\,\theta_h,\,\theta_{ht})+C\|g_t\|^2.
\end{align}
 The nonlinear term in \eqref{3.uap12} is expanded as follows: 
	\begin{align*}
		&\mo_2^{\theta_h}(\bu_{ht},\theta_h,\theta_{ht})=\sum_{T\in\mathcal{T}_h}  \int\limits_T(\bu_{ht}\cdot \nabla \theta_h)\cdot\theta_{ht}\,d\bx +\sum_{T\in\mathcal{T}_h}\int\limits_{\partial T_-}|\{\bu_{ht}\}\cdot \bn_T|(\theta_h^{int}-\theta_h^{ext})\cdot \theta_{ht}^{int}\,ds \\&+\frac{1}{2}\sum_{T\in \mathcal{T}_h}\int\limits_T(\nabla\cdot \bu_{ht})\theta_h\cdot \theta_{ht}\,d\bx
		 -\frac{1}{2}\sum_{e\in \mathcal{N}_h}\int\limits_e[\bu_{ht}]\{\theta_h\cdot \theta_{ht}\}\,ds
		= B_1+B_2+B_3+B_4,\quad { \text{(say)}}.
	\end{align*}
	 Using H\"older, Young's inequality, and 
Lemmas~\ref{ldgl4} and \ref{Enrscript}, we deduce  
  \begin{align}
		|B_1|
        &\leq \epsilon|(\bchi_{ht},\,\theta_{ht})|_\mathcal{A}^2+C|(\bchi_h,\theta_{h})|_\mathcal{A}^2\|\bu_{ht}\|^2+\epsilon|(\su_{ht},\,\bu_{ht})|_\mathcal{A}^2+C|(\bchi_h,\theta_{h})|_\mathcal{A}^2\|\theta_{ht}\|^2.\label{pr1*}
	\end{align}
	Applying H\"older, Young's  and trace inequalities and Lemmas \ref{ldgl4}, \ref{Enrscript}, we arrive at
\begin{align}
		|B_2|
		&\leq \epsilon|(\bchi_{ht},\,\theta_{ht})|_\mathcal{A}^2+C|(\bchi_h,\theta_{h})|_\mathcal{A}^2\|\bu_{ht}\|^2+\epsilon|(\su_{ht},\,\bu_{ht})|_\mathcal{A}^2+C|(\bchi_h,\theta_{h})|_\mathcal{A}^2\|\theta_{ht}\|^2.\label{pr2*}
	\end{align}
	To estimate the term $B_3$, we employ H\"older and 
Young's inequalities along with 
Lemmas~\ref{ldgl4} and \ref{Enrscript}, which yields
	\begin{align}
		|B_3|
		\leq & \epsilon|(\su_{ht},\,\bu_{ht})|_\mathcal{A}^2+\epsilon|(\bchi_{ht},\,\theta_{ht})|_\mathcal{A}^2+C\|\theta_{ht}\|^2 |(\bchi_h,\theta_{h})|_\mathcal{A}^2.\label{pr3*}
	\end{align}
	Using arguments analogous to those employed in the 
estimates of $B_2$ and $B_3$, we  bound 
$B_4$ as
	\begin{align}
		|B_4|\leq \epsilon|(\su_{ht},\,\bu_{ht})|_\mathcal{A}^2+\epsilon|(\bchi_{ht},\,\theta_{ht})|_\mathcal{A}^2+C\|\theta_{ht}\|^2 |(\bchi_h,\theta_{h})|_\mathcal{A}^2.\label{pr4*}
	\end{align} 
\noindent
Utilizing \eqref{pr1*}--\eqref{pr4*} in \eqref{3.uap12}, 
multiplying derived equation by $e^{2\alpha t}$, and using \eqref{ldgLpth} 
together with Lemma~\ref{Enrscript}, we integrate with respect time from $0$ to $t$ to obtain
\begin{align}\label{3.9ap1}
e^{2\al t}(\|\theta_{ht}\|^2&)+(\kappa-2\al \mathcal{N}_\theta C_\theta^2)\int\limits_0^te^{2\al s} |(\bchi_{hs}(s),\,\theta_{hs}(s))|^2_\mathcal{A} ds\leq \|\theta_{ht}(0)\|^2+\epsilon|(\su_{ht},\,\bu_{ht})|_\mathcal{A}^2\nonumber\\
& +C\int\limits_0^te^{2\al s}\|\theta_{ht}(s)\|^2|(\bchi_{h}(s),\,\theta_{h}(s))|^2_\mathcal{A}\,ds+\|\bu_{ht}(s)\|^2|(\bchi_{h}(s),\,\theta_{h}(s))|^2_\mathcal{A}+C\int\limits_0^t e^{2\al s} \|g_t(s)\|\, ds.
\end{align} 
Again, differentiate (\ref{3.we1}) and  (\ref{3.we2}) with respect to time, then substitute  $\tu_h=\nu\bchi_{ht}$, $\psi_h=\theta_{ht}$ in the resulting equations and follow the analysis used for  (\ref{3.9ap1}) to arrive at
\begin{align}\label{3.9ap2}
e^{2\al t}\|\bu_{ht}\|^2&+(\nu-2\al \mathcal{N}_\bu C_\bu^1)\int\limits_0^te^{2\al s} |(\su_{hs}(s),\,\bu_{hs}(s))|^2_\mathcal{A} ds\leq \|\bu_{ht}(0)\|^2+C\int\limits_0^t e^{2\al s} \|\brf_t(s)\|\, ds\nonumber\\
&+C\int\limits_0^te^{2\al s}(\|\bu_{ht}(s)\|^2|(\su_{h}(s),\,\bu_{h}(s))|^2_\mathcal{A}(1+\|\bu_h(s)\|^2)+C\|\theta_{hs}(s)\|^2)\,ds.
\end{align} 
Combine (\ref{3.9ap1}) and (\ref{3.9ap2}),  choose $0<\al <\min\{\frac{\kappa}{2(\mathcal{N}_\theta C_\theta^2+\kappa)},\,\frac{\nu}{2(\mathcal{N}_\bu C_\bu^1)}\}$, utilize (\ref{3.uap7vel}) and employ Grownwall's lemma on the resulting equation. Then,  multiply by $e^{-2\al t}$to obtain
 \begin{align}\label{3.uap15}
\|\theta_{ht}(t)\|^2+\|\bu_{ht}(t)\|^2+e^{-2\al t}\int\limits_0^te^{2\al s} (|(\su_{hs}(s),\,\bu_{hs}(s))|^2_\mathcal{A}+ |(\bchi_{hs}(s),\,\theta_{hs}(s))|^2_\mathcal{A}) ds\leq C(t).
\end{align}
\end{proof}
\noindent
	\noindent We now  introduce the   $\bL^2$-projection $\bP_h^\bu:\bL^2(\Omega)\rightarrow \bX_h^l$, $\underline{\bL}^2$--projection $\underline{\bP}_h^\sib:\underline{\bL}^2(\Omega)\rightarrow \Y_h$, $\bL^2$--projection $\bP_h^{\bchi}:\bL^2(\Omega)\rightarrow \bZ_h$ and $L^2$--projection $\underline{\bP}_h^\theta: L^2(\Omega)\rightarrow Z_h$, such that 
    \begin{align*}
     (\brho-\bP_h^\bu\brho,\,\brho_h)&=0\quad \forall \brho\in\bL^2(\Omega),\,\forall\brho_h\in \bX_h^l,\\
     (\su-\underline{\bP}_h^\sib\su\,\su_h)&=0\quad \forall \su\in\underline{\bL}^2(\Omega)\, \forall \su_h\in \Y_h,\\
      (\bchi-\bP_h^{\bchi}\bchi,\,\bchi_h)&=0 \quad \forall \bchi\in\bL^2(\Omega),\,\forall\bchi_h\in \bZ_h,\\
     (\psi-\mathcal{P}_h^\theta\psi,\,\psi_h)&=0\quad \forall \psi\in L^2(\Omega),\, \forall \psi_h\in Z_h.
    \end{align*}
    \begin{lemma}\label{ldgphproj}
        The projections $\bP_h^\bu$ and $\underline{\bP}_h^\sib$ satisfy
		\begin{align}\label{ldgP_h}
			\|\bv-\bP_h^\bu\bv\|+h\|\bv-\bP_h^\bu\bv\|_\varepsilon&\leq Ch^{r+1}|\bv|_{r+1}\quad \forall\,\, \bv\in   \bJ_1\cap\bH^{r+1}(\Omega),\\
            \|\su-\underline{\bP}_h^\sib\su\|&\leq Ch^r|\su|_r\quad \forall \su\in   \Y,\label{P_hsu}\\
            |(\su-\underline{\bP}_h^\sib\su, \bu-\bP_h^\bu\bu)|_\mathcal{A}&\leq Ch^r(|\su|_r+|\bu|_{r+1})\quad \forall \su \in \underline{\bH}^r, \, \forall \bu \in \bJ_1\cap\bH^{r+1}(\Omega),\label{P_hsubu}\\
            \|\su -\underline{\bP}_h^\sib\su\|_{0,e}&\leq Ch_E^{r-1/2}\|\su\|_{r,E}\quad \forall e\in \partial E, \,\,\forall \su\in \underline{\bH}^r.\label{T5suph}
		\end{align}
         \end{lemma}
    \noindent
   The estimate \eqref{ldgP_h} is a direct consequence of the 
orthogonality property of $\boldsymbol{P}_h^{\boldsymbol{u}}$ 
and Lemma~\ref{ldgih}.
     We refer the reader to \cite[Lemma 3.1]{cockburn02} for the justification of (\ref{P_hsu})--(\ref{T5suph}). \\
     Note that $\bP_h^\bchi$ and $\mathcal{P}_h^\theta$ follow the similar properties as those used for  $\underline{\bP}_h^\sib$ and $\bP_h^\bu$.\\
The following theorem presents the main result of this work 
concerning the semidiscrete error analysis.
\begin{theorem}\label{ldgsemierror}
	Under assumption \textbf{(A2)}, with initial data  $\bu_h(0) = \bP_h^\bu \bu_0$ and $\theta_h(0) = \mathcal{P}_h^\theta \theta_0$, it holds for $0<t\leq T$ that 
   \begin{align*}
			&\|(\bu-\bu_h)(t)\|+ h|((\su-\su_h)(t),\,(\bu-\bu_h)(t))|_\mathcal{A}+h\|(p-p_h)(t)\|\leq C(t)h^{r+1},\\
			&\hspace{1cm}\|(\theta-\theta_h)(t)\|+ h|((\bchi-\bchi_h)(t),\,(\theta-\theta_h)(t))|_\mathcal{A}\leq C(t)h^{r+1},
		\end{align*}
	where	C(t) is a positive constant.
	\end{theorem}
\noindent
The proof of Theorem~\ref{ldgsemierror} is carried out in 
Sections~\ref{ldgsection4} and \ref{ldgsection5}.
\section{ LDG Error Estimates for Velocity and Temperature}
\label{ldgsection4}
\se
In this section, we derive optimal error bounds for the stress tensor error $\underline{\e}^\sib=\su-\su_h$, velocity error $\e^\bu=\bu-\bu_h$, {  temperature gradient error} $\e^\bchi=\bchi-\bchi_h$ and  { temperature error} $e^\theta=\theta-\theta_h$,   for $t>0$.
To begin the analysis, we first consider the associated linearized error problem. 
For this purpose, we introduce  $(\vphi_h,\,\bv_h,\,\vtheta_h,\,w_h)\in\Y_h\times\bV_h\times\bZ_h\times Z_h$, which denotes the LDG approximation of a linearized Stokes--type system. 
More precisely,   $(\vphi_h,\,\bv_h,\,\vtheta_h,\,w_h)$ is defined as the solution of  
\begin{align}
	a_1^l(\vphi_h,\,\tu_h)+b_1^l(\bv_h,\,\tu_h)&=0\:\: \forall \tu_h\in\Y_h,\label{4.we1}\\
	a_2^l(\vtheta_h,\,\biota_h)+b_2^l(\theta_h,\,\biota_h)&=0\:\: \forall \biota_h\in\bZ_h,\label{4.we1th}\\
	(\bv_{ht},\,\bphi_h)-\nu(b_1^l(\bphi_h,\,\vphi_h)-c_1^l(\bv_h,\,\bphi_h))&+({ \boldsymbol{\mu}}w_h, \,\bphi_h)
		= (\textbf{f},\,\bphi_h)
		-\mo_1^\bu(\,\bu,\,\bu,\,\bphi_h) \:\:\forall \bphi_h\in \bX_h^l,\label{4.we2}\\
		(w_{ht},\,\psi_h)-\kappa(b_2^l(\psi_h,\,\vtheta_h)-c_2^l(w_h,\,\psi_h))
		&= (g,\,\psi_h)
		-\mo_2^{\theta}(\bu,\,\theta,\,\psi_h) \:\:\forall \psi_h\in Z_h\label{4.we2th}.
	\end{align}
 We split $\underline{\e}^\sib$, $\e^\bchi,\,\e^\bu$ and $e^\theta$ as follows:
	$\underline{\e}^\sib=\su-\vphi_h+\vphi_h-\su_h=\underline{\bxi}^\sib+\underline{\be}^\sib,\:
	\e^\bu=\bu-\bv_h+\bv_h-\bu_h=\bxi^\bu+\be^\bu$, 
	$\e^\bchi=\bchi-\vtheta_h+\vtheta_h-\bchi_h=\bxi^\bchi+\be^\bchi,\:
	e^\theta=\theta-w_h+w_h-\theta_h=\omega^\theta+\eta^\theta$.\\
 Note that $\underline{\bxi}^\sib$, $\bxi^\bu$, $\bxi^\bchi$ and $\omega^\theta$
correspond to the errors arising from the 
approximation of the linearized (Stokes) problem. 
In contrast,   $\underline{\be}^\sib$, $\be^\bu$, $\be^\bchi$ and $\eta^\theta$ account for the errors generated by the 
nonlinear terms in \eqref{1.1}--\eqref{1.2}.  By comparing  equations (\ref{2.we1})-(\ref{2.we2th}) and (\ref{4.we1})-(\ref{4.we2th}), we obtain the following system 
governing $\underline{\bxi}^\sib$, $\bxi^\bu$, $\bxi^\bchi$ and $\omega^\theta$ 
\begin{align}
   a_1^l(\underline{\bxi}^\sib,\,\tu_h)+b_1^l(\bxi^\bu,\,\tu_h)&=0 \qquad \forall \tu_h\in\Y_h,\label{4.ie1}\\
   a_2^l(\bxi^\bchi,\,\biota_h)+b_2^l(\omega^\theta,\,\biota_h)&=0 \qquad \forall \biota_h\in\bZ_h,\label{4.ie1th}\\
		(\bxi^\bu_t,\,\bphi_h)-\nu(b_1^l(\bphi_h,\,\underline{\bxi}^\sib)-c_1^l(\bxi^\bu,\,\bphi_h))&+(\boldsymbol{\mu}\omega^\theta,\,\bphi_h)= -d\,(\bphi_h,\,p) \qquad \forall \bphi_h\in \bX_h^l\label{4.ie2},\\
		(\omega^\theta_t,\,\psi_h)-\kappa(b_2^l(\psi_h,\,\bxi ^\bchi)-c_2^l(\omega^\theta,\,\psi_h))
		&= 0 \qquad \forall \psi_h\in Z_h\label{4.ie2th}. 
\end{align}
%
Prior to establishing the semidiscrete LDG error bounds stated in  in Theorem \ref{ldgsemierror}, we define  the approximations $\boldsymbol{\bS}_h^\bu\bu \in \bX_h^l$, $\underline{\boldsymbol{\bS}}_h^\sib\su \in \Y_h$, $\boldsymbol{\bS}_h^\bchi\bchi \in \bZ_h$ and $\mathcal{S}_h^\theta\theta \in Z_h$ associated with  $(\su,\, \bu,\,\bchi,\,\theta)$ of (\ref{w0})--(\ref{w1th}), which satisfy
\noindent
\begin{align}
	a_1^l(\su-\underline{\boldsymbol{\bS}}_h^\sib\su,\,\tu_h)+b_1^l(\bu-\boldsymbol{\bS}_h^\bu\bu,\,\tu_h)&=0\label{4.xi00001},\\
	a_2^l(\bchi-\boldsymbol{\bS}_h^\bchi\bchi,\,\biota_h)+b_2^l(\theta-\mathcal{S}_h^\theta\theta,\,\biota_h)&=0\label{4.xi00001th},\\
		-\nu(b_1^l(\bphi_h,\,\su_h-\underline{\boldsymbol{\bS}}_h^\sib\su)-c_1^l(\bu-\boldsymbol{\bS}_h^\bu\bu,\,\bphi_h))&+ (\boldsymbol{\mu}(\theta-\mathcal{S}_h^\theta\theta),\,\bphi_h)= -d\,(\bphi_h,\,p),\label{4.xi2}\\
		-\kappa(b_2^l(\psi_h,\,\bchi_h-\boldsymbol{\bS}_h^\bchi\bchi)-c_2^l(\theta-\mathcal{S}_h^\theta\theta,\,\psi_h))
		&=0,\label{4.xi2th}
		\end{align}
  for all    $ (\tu_h,\,\bphi_h,\,\biota_h,\,\psi_h) \in \Y_h\times \bX_h^l \times \bZ_h\times Z_h$.\\ 
More specifically, given $(\su,\, \bu,\,\bchi,\,\theta, p)$, seek $\boldsymbol{\bS}_h^\bu\bu : [0, \infty) \to \bX_h^l$, $\underline{\boldsymbol{\bS}}_h^\sib\su : [0, \infty) \to \Y_h$, $\boldsymbol{\bS}_h^\bchi\bchi : [0, \infty) \to \bZ_h$ and $\mathcal{S}_h^\theta\theta : [0, \infty) \to Z_h$, such that they satisfy (\ref{4.xi00001})--(\ref{4.xi2th}). With these approximations defined, we decompose $\underline{\bxi}^\sib$, $\bxi^\bu$, $\bxi^\bchi$ and $\omega^\theta$ as 
\begin{align*}
\underline{\bxi}^\sib=\su-\underline{\boldsymbol{\bS}}_h^\sib\su+\underline{\boldsymbol{\bS}}_h^\sib\su-\vphi_h:=\underline{\boldsymbol{\bzeta}}^\sib+\underline{\brho}^\sib, \,
\bxi^\bu=\bu-\boldsymbol{\bS}_h^\bu\bu+\boldsymbol{\bS}_h^\bu\bu-\bv_h:=\bzeta^\bu+\brho^\bu,\\
\bxi^\bchi=\bchi-\boldsymbol{\bS}_h^\bchi\bchi+\boldsymbol{\bS}_h^\bchi\bchi-\vtheta_h:=\bzeta^\bchi+\brho^\bchi,\, 
\omega^\theta=\theta-\mathcal{S}_h^\theta\theta+\mathcal{S}_h^\theta\theta-w_h:=\zeta^\theta+\rho^\theta.
\end{align*}
Applying (\ref{4.xi00001})-(\ref{4.xi2th}), we obtain the equations in $\bzeta^\bu$, $\underline{\bzeta}^\sib$, $\bzeta^\bchi$ and $\zeta^\theta$ as 
\begin{align}
	a_1^l(\underline{\bzeta}^\sib,\,\tu_h)+b_1^l(\bzeta^\bu,\,\tu_h)&=0 \quad \forall \tu_h\in\Y_h,\label{4.zeta1}\\
	a_2^l(\bzeta^\bchi,\,\biota_h)+b_2^l(\zeta^\theta,\,\biota_h)&=0 \quad  \forall \biota_h\in\bZ_h,\label{4.zeta1th}\\
		-\nu(b_1^l(\bphi_h,\,\underline{\bzeta}^\sib)-c_1^l(\bzeta^\bu,\,\bphi_h))+ (\boldsymbol{\mu}\zeta^\theta,\,\bphi_h)
		&= -d\,(\bphi_h,\,\,p) \quad \forall \bphi_h\in \bX_h^l\label{4.zeta2},\\
		-\kappa(b_2^l(\psi_h,\,\bzeta^\bchi)-c_2^l(\zeta^\theta,\,\psi_h)) 
		&= 0 \quad \forall \psi_h\in Z_h\label{4.zeta2th}.
		\end{align}
 \begin{lemma}\label{4.sue0ns}
     Under the conditions stated in Theorem \ref{ldgsemierror}, it follows that for all $0 < t \leq T$,
     \begin{align}\label{0zeen}
          \|\bzeta^\bu(t)\|^2+h^2|(\underline{\bzeta}^\sib(t),\,   \bzeta^\bu(t))|^2_\mathcal{A} &\leq Ch^{2r+2}(|\su|^{2}_{r}+|\bu|^2_{r+1}+|\bchi|^{2}_{r}+|\theta|^2_{r+1}+|p|^2_r),\\
           \|\zeta^\theta(t)\|^2+h^2|(\bzeta^\bchi(t),\,   \zeta^\theta(t))|^2_\mathcal{A} &\leq Ch^{2r+2}(|\bchi|^{2r}_{r}+|\theta|^2_{r+1}),\label{0zeenth}\\
          \label{00zeen}
   \|\bzeta_t^\bu(t)\|^2+h^2|(\underline{\bzeta}_t^\sib(t),\,   \bzeta_t^\bu(t))|^2_\mathcal{A} &\leq Ch^{2r+2}(|\su_t|^{2}_{r}+|\bu_t|^2_{r+1}+|\bchi_t|^{2}_{r}+|\theta_t|^2_{r+1}+|p_t|^2_r),\\
    \|\zeta_t^\theta(t)\|^2+h^2|(\bzeta_t^\bchi(t),\,   \zeta_t^\theta(t))|^2_\mathcal{A} &\leq Ch^{2r+2}(|\bchi_t|^{2r}_{r}+|\theta_t|^2_{r+1}).\nonumber        
      \end{align}
  \end{lemma}	
  \begin{proof}
    Substitute $\biota_h=\kappa \bP_h^\bchi\bzeta^\bchi=\kappa(\bzeta^\bchi-(\bchi-\bP_h^\bchi\bchi$))  in (\ref{4.zeta1th}), and $\psi_h=\mathcal{P}_h^\theta\zeta^\theta=\zeta^\theta-(\theta-\mathcal{P}_h^\theta\theta$)  in (\ref{4.zeta2th}). Next, by summing all the above equations, we obtain
  \begin{align}\label{zeen1}
  \kappa & (a_2^l(\bP_h^\bchi\bzeta^\bchi,\,\bP_h^\bchi\bzeta^\bchi)+c_2^l(\mathcal{P}_h^\theta\zeta^\theta,\,\mathcal{P}_h^\theta\zeta^\theta))=-\kappa a_2^l(\bchi-\bP_h^\bchi\bchi,\,\bP_h^\bchi\bzeta^\bchi)\nonumber\\
  &-\kappa b_2^l(\theta-\mathcal{P}_h^\theta\theta,\,\bP_h^\bchi\bzeta^\bchi )-\kappa c_2^l(\theta-\mathcal{P}_h^\theta\theta,\,{P}_h^\theta {\zeta}^\theta ) 
  = A_1+A_2+A_3, { (\text{say})}.
  \end{align}
 By utilizing Lemmas \ref{ldgcont}, \ref{ldgphproj} and  Young's inequality, we find that 
  \begin{align*}
  A_1+A_3\leq  \epsilon |(\bP_h^\bchi\bzeta^\bchi,\, \mathcal{P}_h^\theta\zeta^\theta)|^2_\mathcal{A} +  Ch^{2r}(|\bchi|_{r}^2+|\theta|_{r+1}^2).
  \end{align*}
  A use of  Cauchy-Schwarz, Young's inequalities, (\ref{ldgP_h}) and Lemma \ref{ldgdistrace} leads to
\begin{align*}
   A_2=&
   \sum_{T\in\mathcal{T}_h}\int\limits_T\nabla(\theta-\mathcal{P}_h^\theta\theta)\cdot \bP_h^\bchi\bzeta^\bchi \,d\bx +\sum_{e\in \mathcal{N}_h\backslash \mathcal{N}_D}\int\limits_e(\{\bP_h^\bchi\bzeta^\bchi\}
   -\bsq \bP_h^\bchi\bzeta^\bchi \bqs \textbf{E}_{12}){ \cdot{\bsq \theta-\mathcal{P}_h^\theta\theta\bqs} }\,ds\nonumber\\
    &+\sum_{e\in \mathcal{N}_D}\int\limits_e \bP_h^\bchi\bzeta^\bchi{ \cdot}((\theta-\mathcal{P}_h^\theta\theta)\bn)\,ds   
   \leq  \epsilon|(\bP_h^\bchi\bzeta^\bchi, \mathcal{P}_h^\theta \zeta^\theta )|_\mathcal{A}^2+ Ch^{2r}|\theta|_{r+1}^2.
\end{align*}
  Utilizing $A_1$-$A_3$ in (\ref{zeen1}), we obtain
  \begin{align}\label{zeen2th}
  \kappa|(\bP_h^\bchi\bzeta^\bchi,\, \mathcal{P}_h^\theta\zeta^\theta)|^2_\mathcal{A}\leq Ch^{2r}(|\bchi|^{2}_{r}+|\theta|^2_{r+1}). 
  \end{align}
  Using (\ref{P_hsubu}), \eqref{zeen2th} and triangle inequality, we find that 
  \begin{align}\label{zeen3th}
 |(\bzeta^\bchi(t),\,   \zeta^\theta(t))|^2_\mathcal{A}\leq Ch^{2r}(|\bchi|^{2}_{r}+|\theta|^2_{r+1}). 
  \end{align}
Taking 
$\tu_h=\nu \underline{\bP}_h^\sib\underline{\bzeta}^\sib=\nu(\underline{\bzeta}^\sib-(\su-\underline{\bP}_h^\sib\su$))  in (\ref{4.zeta1}),  $\bphi_h=\bP_h^\bu\bzeta^\bu=\bzeta^\bu-(\bu-\bP_h^\bu\bu)$  in (\ref{4.zeta2})
and applying arguments similar to those presented above. The only modification arises due to the presence of pressure terms. Utilizing  Cauchy--Schwarz and Young's inequalities along with (\ref{ldg2}), (\ref{ldg2.1}), and Lemma~\ref{ldgtrace}, we arrive at
\begin{align}\label{ldgb2}
  |-d(\bP_h^\bu\bzeta^\bu,\,p-r_h(p))|&=\sum_{e\in \mathcal{N}_h \backslash \mathcal{N}_D}\int\limits_e(\{p-r_h(p)\}-\textbf{D}_{12}\cdot \bsq p-r_h(p) \bqs)[ \bP_h^\bu\bzeta^\bu]\,ds +\sum_{e\in \mathcal{N}_D}\int\limits_e (p-r_h(p))  \bP_h^\bu\bzeta^\bu\cdot \bn\,ds\nonumber\\
  &\leq \epsilon |(\underline{\bP}_h^\sib\underline{\bzeta}^\sib,\, \bP_h^\bu\bzeta^\bu)|^2_\mathcal{A}+ Ch^{2r}|p|^{2}_r.
\end{align} 
An application of  
  (\ref{P_hsubu}), (\ref{zeen3th}), \eqref{ldgb2} and triangle inequality leads to 
  \begin{align*}
 |(\underline{\bzeta}^\sib(t),\,   \bzeta^\bu(t))|^2_\mathcal{A}\leq Ch^{2r}(|\su|^{2}_{r}+|\bu|^2_{r+1}+|\bchi|^{2}_{r}+|\theta|^2_{r+1}+|p|^2_r). 
  \end{align*}
  To estimate $\|\bzeta^\bu\|$, we employ the Aubin--Nitsche duality technique. 
Let $(\bw, r)$ denote the solution of the  Stokes problem \cite[pp.~329]{cockburn02},
   \begin{align}
    -\Delta\bw+\nabla r&=\bzeta^\bu \:\:\:\:\text{in}\:\: \Omega,\label{4.de1}\\
    \nabla\cdot\bw&=0 \:\:\:\:\:\:\text{in}\:\: \Omega,\nonumber\\
    \bw&=0 \:\:\:\:\text{on} \:\:\partial\Omega,\nonumber
    \end{align}
    and assume that the solution satisfies the regularity estimate
    \begin{align}\label{4.de02}
        \|\bw\|_2+\|r\|_1\leq C \|\bzeta^\bu\|.
    \end{align}
    	Further, suppose $y$ be the solution of elliptic problem: 
		\begin{align}
		  -\Delta y =\zeta^\theta, \quad y|_{\partial\Omega}=0,\label{ldg9.9th}
		\end{align}
		satisfying the  regularity bound as 
		\begin{align}
		\|y\|_{2,\Omega}  &\leq C\|\zeta^\theta\|.\label{ldg10.0th}
		\end{align}
    Set $\boldsymbol{\upsilon}=-\nabla y$  in (\ref{ldg9.9th}), we find that
    \begin{align}
       \boldsymbol{\upsilon}&=-\nabla y \:\:\text{in} \:\: \Omega,\label{4.de3th}\\
        \nabla\cdot\boldsymbol{\upsilon} &=\zeta^\theta \:\:\text{in} \:\: \Omega.\label{4.de4th}
       \end{align}
       Taking $L^2$--inner product of  (\ref{4.de3th}) with $\bzeta^\bchi$ and (\ref{4.de4th}) with $\zeta^\theta$, we have 
       \begin{align}
           a_2^l( \boldsymbol{\upsilon},\,\bzeta^\bchi)--b_2^l(y,\,\bzeta^\bchi)&=0, \label{4.de5th}\\
          b_2^l(\zeta^\theta,\,\boldsymbol{\upsilon})+c_2^l(\zeta^\theta,\, y)&=\|\zeta^\theta\|^2.
        \label{4.de6th}
       \end{align} 
               Substituting $\biota_h=\kappa\bP_h^\bchi\boldsymbol{\upsilon}$ in (\ref{4.zeta1th}) and $\psi_h=\mathcal{P}_h^\theta y$ in (\ref{4.zeta2th}), we obtain 
       \begin{align}
       a_2^l(\bzeta^\bchi,\,\bP_h^\bchi\boldsymbol{\upsilon})+b_2^l(\zeta^\theta,\,\bP_h^\bchi\boldsymbol{\upsilon})&=0,\label{4.de7th}\\
       -\kappa(b_2^l(\mathcal{P}_h^\theta y,\,\bzeta^\bchi)-c_2^l(\zeta^\theta,\,\mathcal{P}_h^\theta y)) 
		&= 0\label{4.de8th}.
\end{align}  
From (\ref{4.de5th})--(\ref{4.de8th}), and the continuity of $y$, we arrive at   
  \begin{align}\label{4.de9th}
\kappa\|\zeta^\theta\|^2=& \kappa a_2^l(\bzeta^\bchi,\,\boldsymbol{\upsilon}-- \bP_h^\bchi\boldsymbol{\upsilon})+\kappa c_2^l(\zeta^\theta,\,\theta--\mathcal{P}_h^\theta y)-\kappa b_2^l(y-\mathcal{P}_h^\theta y,\,\bzeta^\bchi)+\kappa b_2^l(\zeta^\theta,\,\boldsymbol{\upsilon}-\bP_h^\bchi\boldsymbol{\upsilon}).
\end{align} 
The estimation of the first three terms on the right--hand side of (\ref{4.de9th}) follows along the same lines as that of the terms in (\ref{zeen1}). Subsequently, applying (\ref{ldg10.0th}) results in
\begin{align}\label{4.de13th}
|\kappa a_2^l(\bzeta^\bchi,\,\boldsymbol{\upsilon}- \bP_h^\bchi\boldsymbol{\upsilon})+\kappa c_2^l(\zeta^\theta,\,\theta-\mathcal{P}_h^\theta y)-\kappa b_2^l(y-\mathcal{P}_h^\theta y,\,\bzeta^\bchi)|  
\leq   Ch^r\|\zeta^\theta\||(\bzeta^\bchi,\,\zeta^\theta)|_\mathcal{A}.
\end{align}   
Applying  Cauchy--Schwarz's  inequality, (\ref{P_hsu}), (\ref{T5suph}), (\ref{ldg10.0th}) and  Lemma \ref{Enrscript}, we find that 
\begin{align}\label{4.de1300th}
  b_2^l&(\zeta^\theta,\,\boldsymbol{\upsilon}-\bP_h^\bchi\boldsymbol{\upsilon})=|\sum_{T\in\mathcal{T}_h}\int\limits_T\nabla\zeta^\theta{ \cdot}(\boldsymbol{\upsilon}-\bP_h^\bchi\boldsymbol{\upsilon}) \,d\bx +\sum_{e\in \mathcal{N}_h\backslash \mathcal{N}_D}\int\limits_e(\{\boldsymbol{\upsilon}-\bP_h^\bchi\boldsymbol{\upsilon}\}-\bsq \boldsymbol{\upsilon}-\bP_h^\bchi\boldsymbol{\upsilon} \bqs \textbf{E}_{12}){ \cdot}{\bsq\zeta^\theta\bqs} \,ds\nonumber\\
    &+\sum_{e\in \mathcal{N}_D}\int\limits_e (\boldsymbol{\upsilon}-\bP_h^\bchi\boldsymbol{\upsilon}).(\bP_h^\bu\zeta^\theta\bn)\,ds|
    \leq \|\nabla\zeta^\theta\|\|\boldsymbol{\upsilon}-\bP_h^\bchi\boldsymbol{\upsilon}\| + Ch|\boldsymbol{\upsilon}|_1c_2^l(\zeta^\theta, \zeta^\theta)^{1/2}
    \leq  Ch^r\|\zeta^\theta\||(\bzeta^\bchi,\,\zeta^\theta)|_\mathcal{A}.
\end{align}    
  A use of (\ref{zeen3th}), (\ref{4.de13th}), (\ref{4.de1300th})  and \textbf{(A2)} yields 
\begin{align}\label{4.de160th}
\|\zeta^\theta\|\leq Ch^{2r+2}.
\end{align}  
  Combining (\ref{zeen3th}) and (\ref{4.de160th}), we find the estimate (\ref{0zeenth}).\\
  By setting $\betau=-\nabla\bw$, $\Tilde{r}=-r$ in \eqref{4.de1}, and choosing $\tu_h=\underline{\bP}_h^\sib\betau$, $\bphi_h=\bP_h^\bu\bz$ in (\ref{4.zeta1})--(\ref{4.zeta2}), 
and 
applying the Cauchy--Schwarz and Young inequalities, projection estimates, trace inequalities, and the regularity of the dual Stokes solution \eqref{4.de02} yields  the estimate (\ref{0zeen}).\\
  Repeating the established set of arguments, we derive the estimates similar to (\ref{0zeen}) and (\ref{0zeenth}) for \(\bzeta^\bu_t\) and \(\zeta^\theta_t\), respectively. The { main} distinction lies in utilizing the equations obtained by differentiating (\ref{4.zeta1}) and (\ref{4.zeta2}), (\ref{4.zeta1th}) and (\ref{4.zeta2th}) in time. Additionally, we substitute $\tu_h=\underline{\bP}_h^\sib\betau_t$ and $\bphi_h=\bP_h^\bu\bz_t$, $\biota_h=\bP_h^\bchi\boldsymbol{\upsilon}_t$ and $\psi_h=\mathcal{P}_h^\theta y_t$ within it. Lastly, for the dual problem, we consider the right-hand side as \(\bzeta^\bu_t\) and \(\zeta^\theta_t\). This concludes the proof of Lemma \ref{4.sue0ns}. 
 \end{proof} 
\begin{lemma}\label{4.2l0}
 Let the assumption \textbf {(A2)} holds { true}, and \( 0<\alpha<\min\{\frac{\nu}{2\mathcal{N}_\bu C_\bu^1},\,\frac{\kappa}{2\mathcal{N}_\theta C_\theta^1}\} \). 
 Then, for \( 0 \leq t < T \), there exists a constant \( C(t)> 0\), such that
\[
\int\limits_{0}^{t} e^{2\alpha s} (\| \bxi^\bu(s) \|^2 +\| \omega^\theta(s) \|^2)ds \leq C(t) h^{2r+2}.
\]
\end{lemma}
 
\begin{proof}
Substitutue $\biota_h=\kappa \bP_h^\bchi \bxi^\bchi$ in (\ref{4.ie1th}), $\psi_h=\mathcal{P}_h^\theta\omega^\theta$ in  (\ref{4.ie2th}) and sum the resulting { equations} to obtain 
{ \begin{align}\label{4.2l1th}
\frac{1}{2}\frac{d}{dt}\|\omega^\theta\|^2&+\kappa|(\bP_h^\bchi \bxi^\bchi,\,  \mathcal{P}_h^\theta\omega^\theta)|^2_{\mathcal{A}}= (\omega^\theta_t,\,\theta-P^\theta_h\theta) -\kappa a_2^l(\bchi-\bP_h^\bchi\bchi,\,\bP_h^\bchi \bxi^\bchi)\nonumber\\
&-\kappa b_2^l(\theta- \mathcal{P}_h^\theta\theta,\,\bP_h^\bchi \bxi^\bchi) + \kappa(b_2^l( \mathcal{P}_h^\theta\omega^\theta,\,\bchi-\bP_h^\bchi\bchi)-\kappa c_2^l(\theta- \mathcal{P}_h^\theta\theta,\,  \mathcal{P}_h^\theta\omega^\theta).
\end{align}}
Observe that 
\begin{align*}
(\omega^\theta_t,\,\theta-P^\theta_h\theta)=\frac{1}{2}\frac{d}{dt}\|\theta- \mathcal{P}_h^\theta\theta\|^2.
\end{align*}
Applying arguments analogous to those used for the right--hand side of (\ref{4.de9th}), we arrive at 
\begin{align*}
|-\kappa a_2^l(\bchi-\bP_h^\bchi\bchi,\,\bP_h^\bchi \bxi^\bchi)-&\kappa b_2^l(\theta- \mathcal{P}_h^\theta\theta,\,\bP_h^\bchi \bxi^\bchi)+ \kappa(b_2^l( \mathcal{P}_h^\theta\omega^\theta,\,\bchi-\bP_h^\bchi\bchi)\\&
-\kappa c_2^l(\theta- \mathcal{P}_h^\theta\theta,\,  \mathcal{P}_h^\theta\omega^\theta)|\leq C h^{2r} (|\bchi|^2_r+|\theta|^2_{r+1})+ \epsilon |(\bP_h^\bchi \bxi^\bchi,\,  \mathcal{P}_h^\theta\omega^\theta)|^2_{\mathcal{A}}.
\end{align*}
Combining all the above estimates in (\ref{4.2l1th}) and applying the triangle inequality along with (\ref{P_hsubu}), we { arrive at} 
\begin{align}\label{4.2l2th}
\frac{1}{2}\frac{d}{dt}\|\omega^\theta\|^2+\kappa|(\bxi^\bchi,\, \omega^\theta)|^2_{\mathcal{A}}\leq &\frac{d}{dt}\|\theta- \mathcal{P}_h^\theta\theta\|^2+ C h^{2r} (|\bchi|^2_r+|\theta|^2_{r+1}).
\end{align}
Multiplying (\ref{4.2l2th}) by $e^{2\alpha t}$ and using (\ref{ldgLp}) along with (\ref{ldgP_h}), we obtain
  \begin{align*}
\frac{d}{dt}(e^{2\alpha t}\|\omega^\theta\|^2)+(\kappa-2\alpha \mathcal{N}_\theta C_\theta^1)|(\bxi^\bchi,\, \omega^\theta)|^2_{\mathcal{A}}\leq &\frac{d}{dt}(e^{2\alpha t}\|\theta- \mathcal{P}_h^\theta\theta\|^2)+ C h^{2r}e^{2\alpha t} (|\bchi|^2_r+|\theta|^2_{r+1}).
\end{align*}
With $0<\alpha<\min\{\frac{\nu}{2\mathcal{N}_\bu C_\bu^1},\,\frac{\kappa}{2\mathcal{N}_\theta C_\theta^1}\}$ and integrating from $0$ to $t$ and observing $\|\omega^\theta(0)\|$ is of the order $h^r$, we find that 
 \begin{align}\label{4.2l3th}
\int\limits_0^te^{2\alpha s}|(\bxi^\bchi(s),\, \omega^\theta(s))|^2_{\mathcal{A}}{ ds}\leq C h^{2r}\int\limits_0^t e^{2\alpha s} (|\bchi(s)|^2_r+|\theta(s)|^2_{r+1})ds.
\end{align}
Substituting $\tu_h=\nu \underline{\bP}_h^\sib \underline{\bxi}^\sib$ in (\ref{4.ie1}), $\bphi_h=\bP_h^\bu\bxi^\bu$ in  (\ref{4.ie2}) and following the similar analysis { as} used in (\ref{4.2l3th}), we obtain 
 \begin{align}\label{4.2l3}
\int\limits_0^te^{2\alpha s}|(\underline{\bxi}^\sib(s),\, \bxi^\bu(s))|^2_{\mathcal{A}}{ ds}\leq & C h^{2r}\int\limits_0^t e^{2\alpha s} (|\su(s)|^2_r
+|\bu(s)|^2_{r+1}+|p(s)|^2_{r})ds.
\end{align}
To bound $\|\bxi^\bu\|$, $\|\omega^\theta\|$, we use  Aubin--Nitsche's duality argument: For fixed $h>0$ and $t\in (0,\,T)${ , let}  $\bz(s) \in\bJ_1 $, $q(s)\in L^2(\Omega)/\R,\, y \in H_0^1(\Omega)$ be the unique solution of   the  backward problem 
   \begin{align}
   \bz_s +\nu\Delta\bz-\nabla q+\boldsymbol{\mu}y&=e^{2\alpha s}\bxi^\bu,\, 0\leq s\leq t, \label{4.2.de1}\\
   y_s +\kappa\Delta y&=e^{2\alpha s}\omega^\theta,\, 0\leq s\leq t, \label{4.2.de1th}
    \end{align}
    with { $(\bz(0),\,y(0)) =\textbf{0}$}  
    satisfying the following regularity:
    \begin{align}\label{ldgback stokes bound}
       &\int\limits_0^t e^{-2\alpha s}(\|\Delta\bz(s)\|+\|\bz_s(s)\|^2+\|\nabla q(s)\|^2)\, ds \leq C\int\limits_0^t e^{2\alpha s}\|\bxi^\bu(s)\|^2\,ds,\\
       &\hspace{1cm}\int\limits_0^t e^{-2\alpha s}(\|\Delta y(s)\|^2+\|y_s(s)\|^2)\, ds \leq C\int\limits_0^t e^{2\alpha s}\|\omega^\theta(s)\|^2\,ds.\label{back elliptic bound}
    \end{align}
    Set $\boldsymbol{\upsilon}={ -}\nabla y$ in (\ref{4.2.de1th}), we find that
    \begin{align}
       \boldsymbol{\upsilon}&=-\nabla y \:\:\text{in} \:\: \Omega,\label{4.2de3th}\\
        y_s-\nu\nabla\cdot\boldsymbol{\upsilon}&=e^{2\alpha s}\omega^\theta \:\:\text{in} \:\: \Omega.\label{4.2de4th}
       \end{align}
       Multiply (\ref{4.2de3th}) by $\bxi^\bchi$ and (\ref{4.2de4th}) by $\omega^\theta$, we have 
       \begin{align}
           a_2^l( \boldsymbol{\upsilon},\,\bxi^\bchi)-b_2^l( y,\,\bxi^\bchi)&=0, \label{4.2de5th}\\
         ( y_s,\,  \omega^\theta)- \kappa  b_2^l( \omega^\theta,\,\boldsymbol{\upsilon})-\kappa  c_2^l(\omega^\theta,\,  y)&=e^{2\alpha s}\| \omega^\theta\|^2.
        \label{4.2de6th}
       \end{align}     
Substituting $\biota_h=\bP_h^\bchi\boldsymbol{\upsilon}$ in (\ref{4.ie1th}) and $\psi_h=\mathcal{P}_h^\theta y$ in (\ref{4.ie2th}), we obtain      
 \begin{align}
   a_2^l(\bxi^\bchi,\,\bP_h^\bchi\boldsymbol{\upsilon})+b_2^l( \omega^\theta,\,\bP_h^\bchi\boldsymbol{\upsilon})&=0,\:\: \label{4.2ie1th}\\
		( \omega^\theta_t,\,\mathcal{P}_h^\theta y)-\kappa (b_2^l(\mathcal{P}_h^\theta y,\,\bxi^\bchi)-c_2^l( \omega^\theta,\,\mathcal{P}_h^\theta y))
		&= 0 \label{4.2ie2th}. 
\end{align}   
From (\ref{4.2de5th})-(\ref{4.2ie2th}), we have 
\begin{align} \label{4.2l4th}   
 e^{2\alpha s}\| \omega^\theta\|^2=&  ( \omega^\theta_t,\,\mathcal{P}_h^\theta y)+  ( y_s,\,  \omega^\theta) -a_2^l(\bxi^\bchi,\,\boldsymbol{\upsilon}-\bP_h^\bchi\boldsymbol{\upsilon})+ \kappa (b_2^l( y-\mathcal{P}_h^\theta y,\,\bxi^\bchi)\nonumber\\
 &-b_2^l( \omega^\theta,\,\boldsymbol{\upsilon}-\bP_h^\bchi\boldsymbol{\upsilon})-c_2^l( \omega^\theta,\, y-\mathcal{P}_h^\theta y). 
\end{align}    
 By employing the definition of $\mathcal{P}_h^\theta$, we express the first two terms on the right--hand side of (\ref{4.2l4th}) in the form
		\begin{align}\label{41th1th}
			(  \omega^\theta,  y _s)+(  \omega^\theta_s,\mathcal{P}_h^\theta   y )&=\frac{d}{ds}(  \omega^\theta,  y )-(  \omega^\theta_s,  y -\mathcal{P}_h^\theta   y )=\frac{d}{ds}(  \omega^\theta,  y )-\frac{d}{ds}(  \omega^\theta,  y -\mathcal{P}_h^\theta  y )+{ (\theta-\mathcal{P}_h^\theta\theta,\,y _s)}.
		\end{align}   
Applying  (\ref{41th1th}), the techniques as utilized on the right hand side of (\ref{4.de9th}) in   (\ref{4.2l4th}), we find that 
\begin{align} \label{4.2l5th} 
 e^{2\alpha s}\| \omega^\theta\|^2\leq &  \frac{d}{ds}(  \omega^\theta, \mathcal{P}_h^\theta y)+ \epsilon e^{-2\alpha s}(\|\boldsymbol{\upsilon}\|^2_1+\| y\|^2_2+\| y_s\|^2)+C\epsilon h^2 e^{2\alpha s}|(\bxi^\bchi,\, \omega^\theta)|^2_{\mathcal{A}}\quad+C\epsilon^2h^{2r+2}e^{2\alpha s}|\theta|^2_{r+1}.
\end{align} 
 On integrating (\ref{4.2l5th}) with respect to { time} from 0 to t and applying (\ref{back elliptic bound}), we obtain 
 \begin{align*}
 \int\limits_0^te^{2\alpha s}\| \omega^\theta(s)\|^2\,ds\leq   (  \omega^\theta(t), \mathcal{P}_h^\theta y(t))&-(  \omega^\theta(0), \mathcal{P}_h^\theta y(0))
 +C\epsilon h^2 \int\limits_0^t e^{2\alpha s}|(\bxi^\bchi(s),\, \omega^\theta(s))|^2_{\mathcal{A}}\,ds\nonumber\\
 &+C\epsilon \int\limits_0^t e^{2\alpha s}\| \omega^\theta(s)\|^2\,ds+C\epsilon^2h^{2r+2}\int\limits_0^t e^{2\alpha s}(|\theta(s)|^2_{r+1}\,ds).
\end{align*}   
  Choosing $\epsilon$ appropriately, we then apply (\ref{4.2l3th}) and the assumption \textbf{(A2)} to  arrive at
  \begin{align}\label{4.2shu1}
  \int\limits_{0}^{t} e^{2\alpha s} \| \omega^\theta(s) \|^2 ds \leq C h^{2r+2}.
  \end{align}    
   Moreover, by choosing $\betau=\nabla\bw$ and $\Tilde{r}=-r$ in (\ref{4.2.de1}), applying arguments analogous to those used in (\ref{4.2shu1}) and employing (\ref{ldgback stokes bound}), we arrive at 
 \begin{align*}
 \int\limits_0^te^{2\alpha s}\|\bxi^\bu(s)\|^2\,ds\leq   ( \bxi^\bu(t), \bP_h^\bu\bz(t))&-( \bxi^\bu(0), \bP_h^\bu\bz(0))
 +C\epsilon h^2 \int\limits_0^t e^{2\alpha s}|(\underline{\bxi}^\sib(s),\,\bxi^\bu(s))|^2_{\mathcal{A}}\,ds\nonumber\\
 &+C\epsilon \int\limits_0^t e^{2\alpha s}\|\bxi^\bu(s)\|^2\,ds+C\epsilon h^{2r+2}\int\limits_0^t e^{2\alpha s}(|\bu(s)|^2_{r+1}+|p(s)|^2_{r})\,ds.
\end{align*}   
  Choosing $\epsilon$ appropriately, we then apply (\ref{4.2l3}) and \textbf{(A2)} to complete the rest of the proof.  
\end{proof}	

\begin{lemma}\label{4.3l0}
Under the assumptions of Lemma \ref{4.sue0ns},  for $0<t\leq T$
\[
\| \bxi^\bu(t) \| + h |(\underline{\bxi}^\sib(t),\, \bxi^\bu(t) )|_{\mathcal{A}} \leq C(t) h^{r+1},
\qquad
\| \omega^\theta(t) \| + h |(\bxi^\bchi(t),\, \omega^\theta(t) )|_{\mathcal{A}} \leq C(t) h^{r+1}.
\]
\end{lemma}
\begin{proof}
From (\ref{4.ie1th}), (\ref{4.ie2th}), (\ref{4.zeta1th}) and (\ref{4.zeta2th}), we have 
\begin{align}
 a_2^l(\brho^\bchi,\,\biota_h)+b_2^l(\rho^\theta,\,\biota_h)&=0 \qquad  \forall \biota_h\in\bZ_h,\label{4.3ie1th}\\
		(\rho^\theta_t,\,\psi_h)-\kappa(b_2^l(\psi_h,\,\brho^\bchi)-c_2^l(\rho^\theta,\,\psi_h))
		&=-({ \zeta^\theta_t},\,\psi_h)\qquad  \forall \psi_h\in Z_h\label{4.3sie2th}.
\end{align} 
Substitute $\biota_h=\kappa\brho^\bchi$  in (\ref{4.3ie1th}),  $\psi_h=\rho^\theta$
 in (\ref{4.3sie2th}) and sum all the resulting equations  to obtain
 \begin{align}\label{ldgrho100}
 \frac{1}{2}\frac{d}{dt}\|\rho^\theta\|^2+\kappa |(\brho^\bchi,\,\rho^\theta)|
 ^2_{\mathcal{A}}\leq -(\bzeta^\bu_t,\,\rho^\theta).
 \end{align}
Multiply (\ref{ldgrho100}) by $e^{2\alpha t}$, integrate the resulting inequality with respect to time, and apply the Cauchy--Schwarz  and Young  inequalities to obtain
\begin{align}\label{4.3l1th}
e^{2\alpha t} \|\rho^\theta(t)\|^2+\kappa\int\limits_0^t &e^{2\alpha s} |(\brho^\bchi(s),\,\rho^\theta(s))|
 ^2_{\mathcal{A}}\, ds\leq C \int\limits_0^t e^{2\alpha s} (\|\rho^\theta(s)\|^2+ \|\zeta_s^\theta(s)\|^2)\,ds\nonumber\\
 &\leq C \int\limits_0^t e^{2\alpha s} (\|\zeta^\theta(s)\|^2+ \|\omega^\theta(s)\|^2+ \|\zeta_s^\theta(s)\|^2)\,ds.
\end{align}
Apply  Lemmas \ref{4.sue0ns}, \ref{4.2l0} in (\ref{4.3l1th}) to arrive at 
\begin{align*}
\|\rho^\theta(t)\|^2+\kappa e^{--2\alpha t}\int\limits_0^t e^{2\alpha s} |(\brho^\bchi(s),\,\rho^\theta(s))|
 ^2_{\mathcal{A}}\, ds\leq Ch^{2r+2}.
\end{align*}
By applying Lemma \ref{ldgdistrace}, we find that
\begin{align}\label{4.3rho1}
\|\rho^\theta(t)\|^2+h^2 |(\brho^\bchi(s),\,\rho^\theta(s))|
 ^2_{\mathcal{A}}\leq Ch^{2r+2}.
\end{align}
A use of (\ref{4.3l1th}), Lemma \ref{4.sue0ns} and \textbf{(A2)} leads to second estimate of this lemma.\\
Furthermore,
from (\ref{4.ie1}), (\ref{4.ie2}), (\ref{4.zeta1}) and (\ref{4.zeta2}), we have 
\begin{align}
 a_1^l(\underline{\brho}^\sib,\,\tu_h)+b_1^l(\brho^\bu,\,\tu_h)&=0\qquad \forall \tu_h\in\Y_h,\label{4.3ie1}\\
		(\brho^\bu_t,\,\bphi_h)-\nu(b_1^l(\bphi_h,\,\underline{\brho}^\sib)-c_1^l(\brho^\bu,\,\bphi_h))
		&=-(\bzeta^\bu_t,\,\bphi_h)\qquad \forall \bphi_h\in \bX_h^l\label{4.3sie2}.
\end{align} 
Substitute $\tu_h=\nu\underline{\brho}^\sib$  in (\ref{4.3ie1}) and $\bphi_h= \brho^\bu$
 in (\ref{4.3sie2})  apply the similar techniques used for (\ref{4.3rho1}) to arrive at
\begin{align*}
\|\brho^\bu(t)\|^2+h^2 |(\underline{\brho}^\sib(s),\,\brho^\bu(s))|
 ^2_{\mathcal{A}}\leq Ch^{2r+2}.
\end{align*}
A use of  Lemma \ref{4.sue0ns} and \textbf{(A2)} leads to the required result.
\end{proof}

\noindent 
{ Now,} we provide the estimates for  $\underline{\be}^\sib$, $\be^\bu$, $\be^\bchi$ and $\eta^\theta$. 
\begin{lemma}\label{4.snl0}
Suppose the assumptions of Theorem \ref{ldgsemierror} hold true. Then, for $t > 0$
\begin{align*}
\|\be^\bu(t)\|+h|(\underline{\be}^\sib(t),\,\be^\bu(t))|^2_\mathcal{A}\leq C(t)h^{r+1},
\quad
\|\eta^\theta(t)\|+h|({\be}^\bchi(t),\,\eta^\theta(t))|^2_\mathcal{A}\leq C(t)h^{r+1}.
\end{align*}
\end{lemma}
\begin{proof}
From (\ref{3.we1th}), (\ref{3.we2th}), (\ref{4.we1th}) and (\ref{4.we2th}), we obtain  
\begin{align}
	&a_2^l(\be^\bchi,\biota_h)+b_2^l(\eta^\theta,\biota_h)=0\:\:\:\:\:\:\: \forall \biota_h\in\bZ_h,\label{4.snl1th}\\
	(\eta^\theta_t,\psi_h)-\kappa(b_2^l(\psi_h,\be^\bchi)&-c_2^l(\eta^\theta,\psi_h))=-(\mo_2^{\theta}(\,\bu,\,\theta,\,\psi_h)-\mo_2^{\theta_h}(\bu_h,\,\theta_h,\,\psi_h)) \:\:\forall \psi_h\in Z_h\label{4.snl2th}.
	\end{align}
	Substitute  $\biota_h=\kappa \be^\bchi$ in (\ref{4.snl1th}), $\psi_h= \eta^\theta$ in (\ref{4.snl2th}) and sum all the resulting equations to obtain
	\begin{align}\label{4.snl4th}
	\frac{1}{2}\frac{d}{dt}\|\eta^\theta\|^2+&\kappa|(\be^\bchi,\,\eta^\theta)|^2_\mathcal{A} \leq|-(\mo_2^{\theta}(\,\bu,\,\theta,\,\eta^\theta)-\mo_2^{\theta_h}(\bu_h,\,\theta_h,\,\eta^\theta))|.
	\end{align}
As $\theta$ is continuous, we have
		\[\mo_2^{\theta}( \bu,\theta,\eta^\theta)=\mo_2^{\theta_h}(\bu,\theta,\eta^\theta).\]
    The terms on the right--hand side of (\ref{4.snl4th}) can be reformulated as 
\begin{align}\label{4non1}
			\mo_2^{\theta_h}(\bu_h,\theta_h,\eta^\theta)-\mo_2^{\theta_h}(\bu,\theta,\eta^\theta) =&-\mo_2^{\theta_h}(\bu_h,\eta^\theta,\eta^\theta)+\mo_2^{\theta_h}(\bxi^\bu,\omega^\theta,\eta^\theta)-\mo_2^{\theta_h}(\bxi^\bu,\theta,\eta^\theta)\nonumber\\&+\mo_2^{\theta_h}(\be^\bu,\omega^\theta,\eta^\theta)-\mo_2^{\theta_h}(\be^\bu,\theta,\eta^\theta)
			-\mo_2^{\theta_h}(\bu,\omega^\theta,\eta^\theta).
		\end{align}
		It is worth noting that the first term on the right--hand side of \eqref{4non1} is non--positive by virtue of (\ref{ldg5.5}) and can therefore be discarded.
 	\begin{align}		
			\mo_2^{\theta_h}(\bu_h,\theta_h,\eta^\theta)-\mo_2^{\theta_h}(\bu,\theta,\eta^\theta)\le &  -\mo_2^{\theta_h}(\bu,\omega^\theta,\eta^\theta)+\mo_2^{\theta_h}(\bxi^\bu,\omega^\theta,\eta^\theta)\nonumber\\
			&-\mo_2^{\theta_h}(\bxi^\bu,\theta,\eta^\theta)+\mo_2^{\theta_h}(\be^\bu,\omega^\theta,\eta^\theta)-\mo_2^{\theta_h}(\be^\bu,\theta,\eta^\theta). \label{ldgc0th}
		\end{align}
		By applying (\ref{ldgtri1}), Lemmas~\ref{ldgl4} and \ref{Enrscript}, together with Young's inequality, we obtain
		\begin{align}
			|\mo_2^{\theta_h}(\be^\bu,\theta,\eta^\theta)|&\leq C \|\be^\bu\|_{L^2(\Omega)}|\theta|_{\bW^{1,4}(\Omega)}\|\eta^\theta\|_\varepsilon \leq \epsilon|(\be^\bchi,\,\eta^\theta|^2_\mathcal{A}+C\|\theta\|_2^2\|\be^\bu\|^2.\label{c4th}
		\end{align}
 By employing the Cauchy--Schwarz and Young inequalities, together with (\ref{ldgLp}) and Lemmas~\ref{ldgtrace}, \ref{ldgl4}, \ref{ldgdistrace}, and \ref{Enrscript}, we arrive at
		\begin{align}\label{c5th}
			|\mo_1^{\theta_h}(\bxi^\bu,\omega^\theta,\eta^\theta)|=& \bigg|\sum_{T\in\mathcal{T}_h} \int\limits_T (\bxi^\bu\cdot \nabla\omega^\theta)\cdot \eta^\theta\,d\bx
			+\sum_{T\in\mathcal{T}_h}\int\limits_{\partial T_-}|\{\bxi^\bu\}\cdot \bn_T|(\omega^{\theta int}-\omega^{\theta ext})\cdot \eta^{\theta int}\,ds \nonumber\\
			&+\frac{1}{2}\sum_{T\in\mathcal{T}_h} \int\limits_T(\nabla\cdot \bxi^\bu)\omega^\theta\cdot \eta^\theta\,d\bx -\frac{1}{2}\sum_{e\in\mathcal{N}_h}\int\limits_e[\bxi^\bu]\cdot\{\omega^\theta\cdot\eta^\theta\}\,ds\bigg|\nonumber\\
			&\leq \epsilon|(\be^\bchi,\,\eta^\theta|^2_\mathcal{A}  +C (\|\omega^\theta\|^2_\varepsilon+\|\bxi^\bu\|_\varepsilon^2)
		\end{align}
		and
        \noindent
		\begin{align}
			|\mo_2^{\theta_h}(\be^\bu,\omega^\theta,\eta^\theta)|=& \bigg|\sum_{T\in\mathcal{T}_h}\int\limits_T (\be^\bu\cdot \nabla\omega^\theta)\cdot \eta^\theta\,d\bx+\sum_{T\in\mathcal{T}_h}\int\limits_{\partial T_-}|\{\be^\bu\}\cdot \bn_T|(\omega^{\theta int}-\omega^{\theta ext})\cdot \eta^{\theta int}\,ds \nonumber\\
			&+\frac{1}{2}\sum_{T\in\mathcal{T}_h} \int\limits_T(\nabla\cdot \be^\bu)\omega^\theta\cdot \eta^\theta\,d\bx -\frac{1}{2}\sum_{e\in\mathcal{N}_h}\int\limits_e[\be^\bu]\cdot \bn\{\omega^\theta\cdot\eta^\theta\}\,ds\bigg|\nonumber\\
			 &\leq \epsilon|(\be^\bchi,\,\eta^\theta|^2_\mathcal{A}+C h^{-2}\|\omega^\theta\|_\varepsilon^2\|\be^\bu\|^2.\label{c3th}
\end{align}
		Utilizing (\ref{ldgintegration}), the first nonlinear term on the right-hand side of (\ref{ldgc0th}) takes the form
		\begin{align}
			\mo_2^{\theta_h}(\bu;\omega^\theta,\eta^\theta)=& -\sum_{T\in\mathcal{T}_h}\int\limits_T(\bu\cdot \nabla \eta^\theta)\cdot \omega^\theta\,d\bx-\frac{1}{2}\sum_{T\in\mathcal{T}_h}\int\limits_T(\nabla\cdot \bu)\omega^\theta\cdot \eta^\theta\,d\bx \nonumber\\
			&+\int\limits_{\mathcal{N}_+}|\bu\cdot \bn|\omega^\theta\cdot \eta^\theta \,ds +\frac{1}{2}\sum_{e\in\mathcal{N}_h}\int\limits_e[\bu]\cdot \bn_e\{\omega^\theta\cdot\eta^\theta\}\,ds \nonumber\\
			& -\sum_{T\in \mathcal{T}_h}\int\limits_{\partial T_-}|\bu\cdot \bn_T|\omega^{\theta ext}\cdot (\eta^{\theta int}-\eta^{\theta ext})\,ds=I_1+I_2+I_3+I_4+I_5, (\text{say}).\label{c1th}
		\end{align}
		As $\bu$ is continuous and $\mathcal{N}_+$ is subset of $\partial \Omega$, we have $I_3 = I_4 = 0$. 
Combining Lemmas~\ref{ldgsobolev2.1} and \ref{Enrscript} with the Cauchy--Schwarz and Young inequalities and (\ref{ldgLp}) yields
		\begin{align}
			|I_1|+|I_2|&\leq C \|\bu\|_{L^\infty(\Omega)}\|\nabla \eta^\theta\|_{L^2(\Omega)}\|\omega^\theta\|+ C\|\nabla \bu\|_{L^4(\Omega)}\|\omega^\theta\|_{L^2(\Omega)}\|\eta^\theta\|_{L^4(\Omega)} \nonumber\\
			&\leq C\|\bu\|_2\|\omega^\theta\|\|\eta^\theta\|_\varepsilon\leq C\|\bu\|_2^2\|\omega^\theta\|^2+ \epsilon|(\be^\bchi,\,\eta^\theta)|^2_\mathcal{A}. \label{c11th}
		\end{align}
		\noindent 
        Utilizing Lemmas \ref{ldgsobolev2.1}, \ref{ldgtrace} and Young's inequality, we { find that} 
		\begin{align}
			|I_5| &\leq C\|\bu\|_{L^\infty(\Omega)}\sum_{T\in \mathcal{T}_h}\|\omega^\theta\|_{L^2(\partial E)}|e|^{1/2-1/2}\|[\eta^\theta]\|_{L^2(\partial E)}
			\nonumber\\
			 &\leq C\|\bu\|_2^2(\|\omega^\theta\|^2+h^2\|\omega^\theta\|_\varepsilon^2)+\epsilon|(\be^\bchi,\,\eta^\theta|^2_\mathcal{A}.\label{c13th}
		\end{align}
				Next, to handle the third term on the right--hand side of (\ref{ldgc0th}), we employ the identity 
\(
(\nabla\cdot \bv)\rho = \nabla\cdot (\bv \rho) - \bv\cdot \nabla \rho,
\)
thereby
\begin{align*}
\int\limits_T(\nabla \cdot \bv)\rho \upsilon\,{ d\bx} 
			=&-\int\limits_T(\bu\cdot \nabla \upsilon)\rho\,{ d\bx}+\int\limits_{\partial T}(\bv^{int} \cdot \bn_T)\rho^{int} \upsilon^{int}\,{ ds} -\int\limits_T(\bv\cdot \nabla \rho) \upsilon\,{ d\bx}.
		\end{align*} 
			By means of the above relation, we find that
\begin{align}
	\mo_2^{\theta_h}&(\bxi^\bu, \theta, \eta^\theta)= \frac{1}{2}\sum_{T\in\mathcal{T}_h}\int\limits_T (\bxi^\bu\cdot \nabla\theta) \eta^\theta\,dx- \frac{1}{2}\sum_{T\in\mathcal{T}_h}\int\limits_T (\bxi^\bu\cdot \nabla \eta^\theta) \theta\,dx\nonumber\\
			&- \sum_{T\in\mathcal{T}_h}\int\limits_{\partial T_-} |\{\bxi^\bu\}\cdot \bn_T|(\theta^{int}-\theta^{ext}) \eta^{\theta int}\,ds-\frac{1}{2}\sum_{e\in\mathcal{N}_h}\int\limits_e[\bxi^\bu]\cdot \bn_e\{\theta \eta^\theta\}\,ds\nonumber\\
			&+ \frac{1}{2}\sum_{T\in \mathcal{T}_h}\int\limits_{\partial E} (\bxi^{\bu int}\cdot \bn_T)\theta^{int}\eta^{\theta int} \,ds
			= I_6+I_7+I_8+I_9+I_{10},\quad  \text{(say)}.\label{c2th}
		\end{align}
		Since $\theta$ is continuous, it follows that $I_8 = 0$.  Applying (\ref{ldgLpth}), we bound $I_6$ and $I_7$ as follows:
		\begin{align}
			|I_{6}| & \leq \|\bxi^\bu\|_{L^2(\Omega)}\|\nabla\theta\|_{L^4(\Omega)}\|\eta^\theta\|_{L^4(\Omega)}\leq \epsilon|(\be^\bchi,\,\eta^\theta)|^2_\mathcal{A}+C\|\theta\|_{2,\Omega}^2\|\bxi^\bu\|_{L^2(\Omega)}^2, \label{c21th} \\
			|I_{7}|&\leq  \|\bxi^\bu\|_{L^2(\Omega)}\|\rho^\theta\|_\varepsilon\|\theta\|_{L^\infty(\Omega)} \leq \epsilon|(\be^\bchi,\,\eta^\theta)|^2_\mathcal{A}+C\|\theta\|_{2,\Omega}^2\|\bxi^\bu\|_{L^2(\Omega)}^2. \label{c22th}
		\end{align}
	 We next express the sum in $I_{10}$ over edges rather than elements and analyze its contribution on an interior edge $e$. 
Let $e$ be shared by $T_r$ and $T_s$ with outward normals $\bn_r$ and $\bn_s$. 
Hence,
		\[\int\limits_e (\bxi^\bu|_{T_r}\cdot \bn_r)\theta|_{T_r} \eta^\theta|_{T_r}\,dr+\int\limits_e(\bxi^\bu|_{T_s}\cdot\bn_s)\theta|_{T_s}\eta^\theta|_{T_s}\,dr=\int\limits_e[(\bxi^\bu\cdot \bn_e)\theta \eta^\theta]\,dr.\]
		By applying  Lemmas \ref{ldgsobolev2.1}, \ref{ldgtrace}, we write
		\begin{align}
			I_{9}+I_{10}=&\frac{1}{2}\sum_{e\in\mathcal{N}_h}\int\limits_e\{\bxi^\bu\}\cdot \bn_e[\theta \eta^\theta]\,ds\, \leq C\, \|\,\theta\,\|_{L^\infty(\Omega)}\,\sum_{e\in \mathcal{N}_h}\,\frac{1}{|e|^{1/2}}\,\|[\eta^\theta]\|_{L^2(e)}\nonumber\\
			& \times |e|^{1/2}\,\|\bxi^\bu\|_{L^2(e)}\,\leq C\,\|\theta\|_{2,\Omega}\, \sum_{T\in \mathcal{T}_h}\,\|\eta^\theta\|_\varepsilon\,(\|\bxi^\bu\|_{L^2(E)}\,+\,h_E\|\nabla\bxi^\bu\|_{L^2(E)})\nonumber\\
			 \leq & \epsilon|(\be^\bchi,\,\eta^\theta)|^2_\mathcal{A}\,+C\, \|\theta\|_{2,\Omega}^2(\|\bxi^\bu\|_{L^2(\Omega)}^2\,+\,h^2\|\bxi^\bu\|_\varepsilon^2).\label{c23th}
		\end{align}
		\noindent
	Inserting (\ref{c4th})--(\ref{c23th}) into (\ref{ldgc0th}), and hence into (\ref{4.snl4}), followed by multiplication of the resulting inequality by $e^{2\alpha t}$, yields 
		\begin{align}
			\frac{d}{dt}&(e^{2\al t}\|\eta^\theta\|+(\kappa - 2\al \mathcal{N}_\theta C_\theta^1)e^{2\al t}|(\be^\bchi,\,\eta^\theta)|^2_\mathcal{A} \leq C(\|\bu\|_{2,\Omega}^2+h^{-2}\|\omega^\theta\|_\varepsilon^2)e^{2\al t}\|\eta^\theta\|_{L^2(\Omega)}^2\nonumber\\&+ C\|\theta\|_2^2\|\be^\bu\|^2
			+Ce^{2\al t}\|\bu\|_{2,\Omega}^2(\|\omega^\theta\|_{L^2(\Omega)}^2+h^2\|\omega^\theta\|_\varepsilon^2)+Ce^{2\al t}\|\brho^\bu\|_\varepsilon^2.\label{ldgc6th}
			\end{align}	
Furthermore, from (\ref{3.we1}), (\ref{3.we2}) and (\ref{4.we1}), (\ref{4.we2}), we obtain  
\begin{align}
	&a_1^l(\underline{\be}^\sib,\tu_h)+b_1^l(\be^\bu,\tu_h)=0\:\:\:\:\:\:\: \forall \tu_h\in\Y_h,\label{4.snl1}\\
	(\be^\bu_t,_h)-\nu(b_1^l(\bphi_h,\underline{\be}^\sib)&-c_1^l(\be^\bu,\bphi_h))=-(\mo_1^\bu(\,\bu,\,\bu,\,\bphi_h)-\mo_1^{\bu_h}(\,\bu_h,\,\bu_h,\,\bphi_h)) \:\:\forall \bphi_h\in \bX_h^l\label{4.snl2}.
	\end{align}
Furthermore,	substitute  $\tu_h=\nu \underline{\be}^\sib$ in (\ref{4.snl1}), $\bphi_h= \be^\bu$ in (\ref{4.snl2}) and sum all the resulting equations to obtain
\begin{align}\label{4.snl4}
	\frac{1}{2}\frac{d}{dt}\|\be^\bu\|^2+\nu|(\underline{\be}^\sib,\,\be^\bu)|^2_\mathcal{A}\leq &|-(\mo_1^\bu(\,\bu,\,\bu,\,\be^\bu)-\mo_1^{\bu_h}(\,\bu_h,\,\bu_h,\,\be^\bu))|-(\boldsymbol{\mu}\eta^\theta,\,\be^\bu).
	\end{align}	
	By  following the similar steps as those applied to obtaining (\ref{ldgc6th}), we find that 
\begin{align}\label{c6thvel}
\frac{1}{2}\frac{d}{dt}e^{2\alpha t}(\|\be^\bu\|^2+(\nu-2\alpha(\mathcal{N}_\bu C_\bu^1))|(\underline{\be}^\sib,\,\be^\bu)|^2_\mathcal{A}\leq& e^{2\al t}C(\|\bu\|^2_2+h^{-2}\|\bxi^\bu\|^2)e^{2\alpha t}(\|\be^\bu\|^2)+C\|\eta^\theta\|\nonumber\\
&+Ce^{2\alpha t}\|\bu\|^2_2(\|\bxi^\bu\|^2+h^2\|\bxi^\bu\|^2_\varepsilon)+Ce^{2\alpha t}\|\bxi^\bu\|^4.
\end{align}	
From (\ref{ldgc6th}) and (\ref{c6thvel}), integrating with respect to time from $0$, $t$ and using $\eta^\theta(0)=0$, $\be^\bu(0)=0$, an application of Gr\"onwall's inequality along with (\ref{4.2l3}) yields the required estimate.
\end{proof}
\noindent
\textbf{Proof of velocity and temperature error estimate in Theorem \ref{ldgsemierror}.} Utilizing $\e^\bu = \bxi^\bu + \be^\bu$, $e^\theta = \omega^\theta + { \eta^\theta}$, the triangle inequality, and Lemmas  \ref{4.3l0} and \ref{4.snl0}, we arrive the desired result. \hspace{9.5cm}
$\square$\\
{\bf{ Uniform in time error estimates:}}
Assume that the data satisfy the following smallness conditions:
		{\fontsize{10}{12}\selectfont \begin{align}
	N_1=\sup_{\bz_h,\bw_h\in\bV_h}\frac{\mo_1^{\bw_h} (\bw_h,\bv_h,\bw_h)}{\|\bw_h\|_\varepsilon^2\|\bv_h\|_\varepsilon}\, \text{and}\,
	 \frac{N_1 \mathcal{N}_\bu (C
_\bu^1)^{5/2}}{\nu^2}\|\brf\|_{ L^\infty(0,\infty;\textbf{L}^2(\Omega))}<1,\label{ldguniqueness}\\
	N_2=\sup_{\bw_h\in\bV_h,\,\psi_h,\, y_h \in Z_h}\frac{\mo_2^{l_h}(\bw_h,\psi_h,y_h)}{\|\bw_h\|_\varepsilon\|\psi_h\|_\varepsilon\|y_h\|_\varepsilon}\,\text{and}\,  \frac{N_2 \mathcal{N}_\theta (C
_\theta^1)^{5/2}}{\kappa^2}\|g\|_{ L^\infty(0,\infty;L^2(\Omega))}<1.	\nonumber	
	\end{align}}
	Then, the error estimates stated in Theorem~\ref{ldgsemierror} hold uniformly in time. More precisely,
\begin{align*}
		\|\,(\bu-\bu_h)(t)\,\|+\|\,(\theta-\theta_h)(t)\,\| \,+h(\,|((\su-\su_h)(t),(\bu-\bu_h)(t))\,|_\mathcal{A}
		+|((\bchi-\bchi_h)(t),(\theta-\theta_h)(t))\,|_\mathcal{A})\,\leq  C\,h^{r+1},
	\end{align*}
where $C > 0$ is a constant independent of $h$ and the time $t$.
\\
We reformulate the nonlinear terms as follows:
\begin{align}\label{univel}
		\mo_1^{\bu_h}(\bu_h,\,\bu_h,\,\be^\bu)-\mo_1^{\bu_h}(\,\bu,\,\bu,\,&\be^\bu)=-\mo_1^{\bu_h}(\,{ \bu_h},\,\be^\bu,\,\be^\bu)-\mo_1^{\bu_h}(\,\be^\bu,\,\bu_h,\,\be^\bu)\nonumber\\&+\mo_1^{\bu_h}(\,\bxi^\bu,\,\bxi^\bu,\,\be^\bu)- \mo_1^{\bu_h}(\,\bu,\,\bxi^\bu,\,\be^\bu)-\mo_1^{\bu_h}(\,\bxi^\bu,\,\bu,\,\be^\bu).
		\end{align}
By virtue of (\ref{ldg5.5}), the first term on the right--hand side of (\ref{univel}) is non--positive and therefore discarded. 
Employing the arguments from the proof of Lemma~\ref{4.snl0}, we estimate the third, fourth, and fifth terms on the right--hand side of (\ref{univel}) as \begin{align*}
		|\mo_1^{\bu_h}(\,\bxi^\bu,\,\bxi^\bu,\,\be^\bu)&+ \mo_1^{\bu_h}(\,\bu,\,\bxi^\bu,\,\be^\bu)+\mo_1^{\bu_h}(\,\bxi^\bu,\,\bu,\,\be^\bu)|\leq C(\|\bxi^\bu\|+h\|\bxi^\bu\|_\varepsilon)(\|\bu\|_2+h^{-1}\|\bxi^\bu\|_\varepsilon)\|\be^\bu\|_\varepsilon.
		\end{align*}
		Using (\ref{ldguniqueness}) together with Lemma~\ref{Enrscript}, we find that
\begin{align*}
|\mo_1^{\bu_h}(\,\be^\bu,\,\bu_h,\,\be^\bu)| \leq N_1\|\be^\bu\|^2_\varepsilon\|\bu_h\|_\varepsilon\leq N_1C
_1^{3/2}|(\underline{\be}^\sib,\,\be^\bu)|^2_\mathcal{A}|(\su_h,\,\bu_h)|_\mathcal{A}.
\end{align*}		
Following a modification of the proof of Lemma~\ref{4.snl0}, we first start (\ref{4.snl4}) in which the last term is handled using (\ref{ldgc6th}), leading to
\begin{align*}
\frac{d}{dt}\|\be^\bu\|^2+2(\nu-{N_1}C_1^{3/2})|(\underline{\be}^\sib,\,\be^\bu)|^2_\mathcal{A}\leq Ch^2|(\underline{\be}^\sib,\,\be^\bu)|^2_\mathcal{A}.
\end{align*}	
Multiply the above equation by
$e^{2\alpha t}$ and then, integrate with respect to time from $0$ to $t$. Dividing the resulting equation by $e^{2\alpha t}$, we arrive at
\begin{align*}
\|\be^\bu\|^2+2e^{-2\alpha t}\int\limits_0^te^{2\alpha s}&(\nu-N_1(C
_\bu^1)^{3/2}|(\su_h,\,\bu_h)|_\mathcal{A})|(\underline{\be}^\sib,\,\be^\bu)|^2_\mathcal{A} ds\\ &\leq  e^{-2\alpha t}(\|\be^\bu(0)\|^2)
+2\alpha e^{-2\alpha t}\int\limits_0^te^{2\alpha s}\|\be^\bu\|^2ds +Ch^2e^{-2\alpha t}\int\limits_0^te^{2\alpha s}|(\underline{\be}^\sib,\,\be^\bu)|^2_\mathcal{A}\,ds.
\end{align*}
Taking the limit as $t \to \infty$, and employing (\ref{3.uap9}) in conjunction with L'H\"opital's rule, it follows that
\begin{align*}
\frac{1}{\al}\left(\nu-\frac{NC
_1^{5/2}\mathcal{N}  \|\brf\|_{L^\infty(\bL^2(\Omega))}}{\nu}\right)\limsup_{t\rightarrow \infty}|(\underline{\be}^\sib,\,\be^\bu)|^2_\mathcal{A}\leq Ch^{2r+2}.
\end{align*}
Using (\ref{ldgLp}), (\ref{ldguniqueness}), and Lemma~\ref{Enrscript}, we conclude that
\begin{align}
\limsup_{t\rightarrow \infty}|(\underline{\be}^\sib,\,\be^\bu)|^{2}_\mathcal{A}\leq Ch^{2r+2}
 \quad\text{and}\quad \limsup_{t\rightarrow \infty}\|\be^\bu(t)\|\leq Ch^{r+1}.\label{univel2}
\end{align}
A combination of Lemma \ref{4.sue0ns} and (\ref{univel2} ) yields 
\begin{align*}
\limsup_{t\rightarrow \infty}\|\bu(t)-\bu_h(t)\|\leq Ch^{r+1}.  
\end{align*}
By an analogous argument, we obtain the following estimates for $\theta$:
		\[\limsup_{t\rightarrow \infty}\|\theta(t)-\theta_h(t)\|\leq  Ch^{r+1}. \hspace{6cm} \quad \square\] 
We next derive the LDG semidiscrete error estimates for the pressure.
\section{   Error Estimates for Pressure  }\label{ldgsection5}
\se
This section is devoted to proving a lemma that plays a fundamental role in the derivation of the semidiscrete pressure error estimates.
\begin{lemma}\label{ldgetL2L2err}
 Suppose the assumption {\bf{(A2)}} is satisfied and $0<\alpha<\min\{\frac{\nu}{2\mathcal{N}_\bu C_\bu^1},\,\frac{\kappa}{2\mathcal{N}_\theta C_\theta^1}\}$. Then for $0<t\leq T$, 
		\begin{equation*}
			\int\limits_0^te^{2\al s}(\|\e_s^\bu(s)\|^2+\|e_s^\theta(s)\|)\,ds\,\leq\,  C(t)  h^{2r},
		\end{equation*}
		where   $C(t)$ is a positive constant. 
	\end{lemma}
	\begin{proof}	
	From (\ref{2.we1}), (\ref{2.we2}), (\ref{3.we1}), (\ref{3.we2}), (\ref{4.xi00001}), and (\ref{4.xi2}) with 
	$\underline{\boldsymbol{\Upsilon}}^\sib=\underline{\boldsymbol{\bS}}_h^\sib\su-\su_h$, $\boldsymbol{\Upsilon}^\bu=\boldsymbol{\bS}_h^\bu\bu-\bu_h$ and $\Upsilon^\theta=\mathcal{S}_h^\theta\theta-\theta_h$, we have 
	\begin{align}
	\label{4.5ps1}
	a_1^l(\underline{\boldsymbol{\Upsilon}}^\sib_t,\tu_h)+b_1^l(\boldsymbol{\Upsilon}^\bu_t,\tu_h)&=0,\\
	\label{4.5ps2}
	\hspace{-0.4cm}(\boldsymbol{\Upsilon}^\bu_t,\bphi_h)-\nu(b_1^l(\bphi_h,\underline{\boldsymbol{\Upsilon}}^\sib)-c_1^l(\boldsymbol{\Upsilon}^\bu,\bphi_h))+\mo_1^\bu(\,\bu,\,\bu,\,\bphi_h)
		&-\mo_1^{\bu_h}(\,\bu_h,\,\bu_h,\,\bphi_h)+ (\boldsymbol{\mu}\Upsilon^\theta,\,\bphi_h)= -(\bzeta^\bu_t,\,\bphi_h),
	\end{align}
for all $(\tu_h,\,\bphi_h)\in (\Y_h,\,\bX_h^l)$.	
	Substitute $\tu_h=\nu \underline{\boldsymbol{\Upsilon}}^\sib$ in (\ref{4.5ps1}), $\bu_h=\boldsymbol{\Upsilon}^\bu_t$ in (\ref{4.5ps2}) and sum all the resulting equations to obtain
	\begin{align}\label{4.5ps3}
	\|\boldsymbol{\Upsilon}^\bu_t\|+\frac{\nu}{2}\frac{d}{dt}(a_1^l(\underline{\boldsymbol{\Upsilon}}^\sib,\, \underline{\boldsymbol{\Upsilon}}^\sib)&+ c_1^l(\boldsymbol{\Upsilon}^\bu,\,\boldsymbol{\Upsilon}^\bu)) +\mo_1^\bu(\,\bu,\,\bu,\,\boldsymbol{\Upsilon}^\bu_t)
		-\mo_1^{\bu_h}(\,\bu_h,\,\bu_h,\,\boldsymbol{\Upsilon}^\bu_t)+ (\boldsymbol{\mu}\Upsilon^\theta,\,\boldsymbol{\Upsilon}_t^\bu)= -(\bzeta^\bu_t,\,\boldsymbol{\Upsilon}^\bu_t).
	\end{align}
	Applying  Cauchy-Schwarz  and Young's inequalities, we obtain
	\begin{align}\label{4.5ps4vel}
	|(\boldsymbol{\mu}\Upsilon^\theta,\,\boldsymbol{\Upsilon}_t^\bu)|\leq \epsilon\|\boldsymbol{\Upsilon}_t^\bu\|^2+ C\|\Upsilon^\theta\|^2.
	\end{align}
{ The trilinear term can be expressed in the following equivalent form:}
		\begin{align}
		\mo_1^\bu(\,\bu,\,\bu,\,\boldsymbol{\Upsilon}^\bu_t)
		-\mo_1^{\bu_h}(\,\bu_h,\,\bu_h,\,\boldsymbol{\Upsilon}^\bu_t)=-\mo_1^{\bu_h}(\,\e^\bu,\,\e^\bu,\,\boldsymbol{\Upsilon}^\bu_t)+\mo_1^{\bu_h}(\,\e^\bu,\,\bu,\,\boldsymbol{\Upsilon}^\bu_t)+\mo_1^{\bu_h}(\,\bu,\,\e^\bu,\,\boldsymbol{\Upsilon}^\bu_t).
		\end{align}	
				We apply (\ref{ldgtri20}), Lemma \ref{ldgdistrace}, Theorem \ref{ldgsemierror} and Young's inequality to estimate $\mo_1^{\bu_h}(\,\de^\bu, \e^\bu,\boldsymbol{\Upsilon}^\bu_t)$ as 
		\begin{align}
	\hspace{-0.4cm}		|\mo_1^{\bu_h}(\,\de^\bu,\,\e^\bu,\,\boldsymbol{\Upsilon}^\bu_t)|&\leq\, C\,\|\de^\bu\|_\varepsilon\|\e^\bu\|_\varepsilon\,\|\boldsymbol{\Upsilon}^\bu_t\|_\varepsilon\,\leq Ch^{-1}\|\e^\bu\|_\varepsilon\,\|\boldsymbol{\Upsilon}^\bu_t\|\leq C \|\e^\bu\|_\varepsilon\,\|\boldsymbol{\Upsilon}^\bu_t\|\leq\, \epsilon\|\boldsymbol{\Upsilon}^\bu_t\|^2+C\|\e^\bu\|_\varepsilon^2.\label{4.5ps4}
		\end{align}
Employing the continuity of $\bu$, along with the H\"older, Young, and Sobolev inequalities, Lemma~\ref{ldgdistrace}, and assumption {\bf (A2)}, we conclude that
		\begin{align}
			|\mo_1^{\bu_h}(\,\de^\bu,\bu,\boldsymbol{\Upsilon}^\bu_t)|\,\leq \epsilon\|\boldsymbol{\Upsilon}^\bu_t\|^2\,+\,C\|\de^\bu\|_\varepsilon^2.\label{4.5ps5}
		\end{align}
	By an analogous argument, we estimate
		\begin{align}
			|\mo_1^{\bu_h}(\,\bu,\e^\bu,\boldsymbol{\Upsilon}^\bu_t)|\leq \epsilon\|\boldsymbol{\Upsilon}^\bu_t\|^2+C\|\e^\bu\|_\varepsilon^2.\label{4.5ps6}
		\end{align}
Substitute (\ref{4.5ps4vel})-(\ref{4.5ps6})	in (\ref{4.5ps3}) and multiply the resulting inequality by $e^{2\alpha t}$. Then, integrate from 0 to t with respect to time and applying Lemma \ref{ldgcont}  to arrive at
\begin{align}\label{ldgb4}
\int\limits_0^t e^{2\alpha s}\|\boldsymbol{\Upsilon}^\bu_t(s)\|^2\,ds+\frac{\nu}{2}e^{2\alpha s}|(\underline{\boldsymbol{\Upsilon}}^\sib,&\,\boldsymbol{\Upsilon}^\bu)|^2_{\mathcal{A}}\leq \int\limits_0^t e^{2\alpha s} (|(\underline{\boldsymbol{\Upsilon}}^\sib(s),\,\boldsymbol{\Upsilon}^\bu(s))|^2_{\mathcal{A}}+\|\Upsilon^\theta(s)\|^2)ds\nonumber\\
&+\int\limits_0^t e^{2\alpha s} \|\bzeta^\bu_s(s)\|^2\,ds
+  \int\limits_0^t e^{2\alpha s} \|\e^\bu(s)\|^2_\varepsilon\,ds+|(\underline{\boldsymbol{\Upsilon}}^\sib(0),\,\boldsymbol{\Upsilon}^\bu(0))|^2_{\mathcal{A}}.
\end{align}
Furthermore, by applying (\ref{0zeen}), (\ref{0zeenth}) and \textbf{(A2)} in \eqref{ldgb4}, we have
\begin{align*}
\int\limits_0^t e^{2\alpha s} |(\underline{\boldsymbol{\Upsilon}}^\sib(s),\,\boldsymbol{\Upsilon}^\bu(s))|^2_{\mathcal{A}}ds\leq C e^{2\alpha t} h^{2r} + C \int\limits_0^t  e^{2\alpha s} |(\underline{\e}^\sib(s),\,\e^\bu(s))|^2_\mathcal{A}\,ds.
\end{align*}
By utilizing (\ref{00zeen}) and \textbf{(A2)}, we obtain
\begin{align*}
\int\limits_0^t e^{2\alpha s}\|\boldsymbol{\Upsilon}^\bu_t(s)\|^2\,ds\leq C e^{2\alpha t} h^{2r}.
\end{align*}
By an analogous argument, we bound
		\begin{align*}
\int\limits_0^t e^{2\alpha s}\|\Upsilon^\theta_t(s)\|^2\,ds\leq C e^{2\alpha t} h^{2r}.
\end{align*}
Moreover, combining the triangle inequality with (\ref{00zeen}) and assumption {\bf (A2)} leads to the required estimate.	
\end{proof}
\begin{lemma}\label{5.et0} Under assumption {\bf (A2)}, and for 
\(
0 < \alpha < \min\left\{\frac{\nu}{2\mathcal{N}_\bu C_\bu^1},\,\frac{\kappa}{2\mathcal{N}_\theta C_\theta^1}\right\}
\)
and $t>0$
\begin{align*}
\|\e_t^\bu(t)\|\leq Ce^{Ct}h^r. 
\end{align*}
\end{lemma}
\begin{proof}
From (\ref{2.we1}),  (\ref{2.we2}), (\ref{3.we1}) and (\ref{3.we2}), we find that
	\begin{align}
	&\hspace{2.5cm}a_1^l(\underline{\e}^\sib,\,\tu_h)+b_1^l(\e^\bu,\,\tu_h)=0 \label{5.et1}\quad \forall \tu_h \in \Y_h,\\
	&(\e^\bu_t,\,\bphi_h)-\nu(b_1^l(\bphi_h,\,\underline{\e}^\sib)-c_1^l(\e^\bu,\,\bphi_h))
		+d(\bphi_h,\,p)+(\mu e^\theta,\,\bphi_h)\nonumber\\
		&\hspace{3cm}= \mo_1^\bu(\,\bu,\,\bu,\,\bphi_h)-\mo_1^{\bu_h}(\,\bu_h,\,\bu_h,\,\bphi_h) \quad \forall \bphi_h,\in \bX_h^l. \label{5.et2}
	\end{align}
Differentiating (\ref{5.et1}) and (\ref{5.et2}) with respect to time, we arrive at
\begin{align}
	&\hspace{3cm}a_1^l(\underline{\e}^\sib_t,\,\tu_h)+b_1^l(\e^\bu_t,\,\tu_h)=0 \label{5.et3}\quad \forall \tu_h \in \Y_h,\\
	&(\e^\bu_{tt},\,\bphi_h)-\nu(b_1^l(\bphi_h,\,\underline{\e}^\sib_t)-c_1^l(\e^\bu_t,\,\bphi_h))
		+d(\bphi_h,\,p_t)+(\boldsymbol{\mu} e^\theta_t,\,\bphi_h)
		= \mo_1^\bu(\,\bu_t,\,\bu,\,\bphi_h)\nonumber\\
		&\hspace{3cm}-\mo_1^{\bu_h}(\,\bu_{ht},\,\bu_h,\,\bphi_h)
		+\mo_1^\bu(\,\bu,\,\bu_t,\,\bphi_h)-\mo_1^{\bu_h}(\,\bu_h,\,\bu_{ht},\,\bphi_h) \quad \forall \bphi_h \in \bX_h^l.\label{5.et30}
	\end{align}	
	Substitute $\tu_h=\nu\, \underline{\bP}_h^\sib \underline{\e}^\sib_t=\nu\,(\underline{\e}^\sib_t- (\su_t-\underline{\bP}_h^\sib\su_t))$
in (\ref{5.et3})  and $\bphi_h= \bP_h^\bu\e^\bu_t=\e^\bu_t- (\bu_t-\bP_h^\bu\bu_t)$ in (\ref{5.et30}). Then, sum all the resulting equations to obtain
\begin{align}\label{5.et4}
\frac{d}{dt}&\|\e^\bu_t\|^2 +\nu\, (a_1^l(\underline{\bP}_h^\sib \underline{\e}^\sib_t,\, \underline{\bP}_h^\sib \underline{\e}^\sib_t)+c_1^l(\bP_h^\bu\e^\bu_t,\,\bP_h^\bu\e^\bu_t))=\frac{d}{dt}\|\bu-\bP_h^\bu\bu_t\|^2-\nu\, a_1^l(\su_t-\underline{\bP}_h^\sib\su_t,\,\underline{\bP}_h^\sib \underline{\e}^\sib_t)\nonumber\\
&- \nu\, c_1^l(\bu_t-\bP_h^\bu\bu_t,\,\bP_h^\bu\e^\bu_t)-\nu\, b_1^l(\bu_t-\bP_h^\bu\bu_t,\,\underline{\bP}_h^\sib \underline{\e}^\sib_t)
+\nu\, b_1^l( \bP_h^\bu\e^\bu_t,\,\su_t-\underline{\bP}_h^\sib\su_t)
-d(\bP_h^\bu\e^\bu_t,\,p_t-r_h(p_t))\nonumber\\
&-(\mu e^\theta_t,\,\bP_h^\bu\e^\bu_t)+\mo_1^\bu(\,\bu_t,\,\bu,\,\bP_h^\bu\e^\bu_t)-\mo_1^{\bu_h}(\,\bu_{ht},\,\bu_h,\,\bP_h^\bu\e^\bu_t)+\mo_1^\bu(\,\bu,\,\bu_t,\,\bP_h^\bu\e^\bu_t)-\mo_1^{\bu_h}(\,\bu_h,\,\bu_{ht},\,\bP_h^\bu\e^\bu_t)\nonumber\\
=& F_1+F_2+F_3+F_4+F_5+F_6+F_7+F_8+F_9+F_{10}.\end{align}
We proceed to estimate the bilinear terms on the right--hand side of (\ref{5.et4}). 
Using the Cauchy--Schwarz and Young inequalities, (\ref{ldgP_h}), (\ref{P_hsubu}), Lemma~\ref{Enrscript}, and arguments similar to those in (\ref{4.2l1th}), we deduce
 \begin{align}
\label{prb1}
    \hspace{-0.2cm}|F_2+F_3+...+F_7| \leq  \epsilon |( \underline{\bP}_h^\sib \underline{\e}^\sib_t, \bp_h^\bu\e^\bu_t)|^2_\mathcal{A}+ Ch^2(|\su_t|^2_1+|\bu_t|^2_2+|p_t|^2_1)+C\|e^\theta_t\|.
\end{align}
\noindent
Using techniques similar to those in Lemma~\ref{4.snl0}, along with Lemmas~\ref{ldgtrace} and \ref{ldgdistrace} and the Young and Cauchy--Schwarz inequalities, we obtain 	 \begin{align}
			&|\mo_1^{\bu_h}(\,\bu_{ht},\bu_h,\boldsymbol{\bP}_h^\bu\de^\bu_t)+(\mo_1^{\bu_h}(\,\bu_h,\bu_{ht},\boldsymbol{\bP}_h^\bu\de^\bu_t)-\mo_1^{\bu_h}(\,\bu_t,\bu,\boldsymbol{\bP}_h^\bu\de^\bu_t)-\mo_1^{\bu_h}(\,\bu,\bu_t,\boldsymbol{\bP}_h^\bu\de^\bu_t)|\nonumber\\
			&=|-\mo_1^{\bu_h}(\,\bu_{ht},\de^\bu,\boldsymbol{\bP}_h^\bu\de^\bu_t)-\mo_1^{\bu_h}(\,\de^\bu_t,\bu,\boldsymbol{\bP}_h^\bu\bu\de^\bu_t)
			-\mo_1^{\bu_h}(\,\de^\bu,\bu_{ht},\boldsymbol{\bP}_h^\bu\de^\bu_t)-\mo_1^{\bu_h}(\,\bu,\de^\bu_t,\boldsymbol{\bP}_h^\bu\de^\bu_t)|\nonumber\\
			&\leq  \epsilon|( \underline{\bP}_h^\sib \underline{\e}^\sib_t, \bp_h^\bu\e^\bu_t)|^2_\mathcal{A}+C(|( \underline{\e}^\sib,\, \e^\bu)|^2_\mathcal{A}|( \su_{ht},\, \bu_{ht})|^2_\mathcal{A}+(\|\de^\bu_t\|+ h^{2r}\|\bu_t\|_{r+1}^2)\|\bu\|_{2}^2).\label{4.6nolinear}
		\end{align}
Applying (\ref{prb1})--(\ref{4.6nolinear}) in (\ref{5.et4}) together with {\bf (A2)}, and multiplying by $e^{2\alpha t}$ followed by integration over $[0,t]$, we arrive at
		 \begin{align*}
			&e^{2\al t}\|\de^\bu_t(t)\|+\nu \int\limits_0^t e^{2\al s}|( \underline{\bP}_h^\sib \underline{\e}^\sib_t, \bp_h^\bu\e^\bu_t)|^2_\mathcal{A}\,ds\leq C\int\limits_0^t e^{2\al s}\|\de^\bu_s(s)\|\,ds+e^{2\al t}\|(\bu_t-\boldsymbol{\bP}_h^\bu\bu_t)(t)\|+C\int\limits_0^t e^{2\al s}\|e^\theta_s(s)\|\,ds \\ &+C\int\limits_0^t e^{2\al s}|( \underline{\e}^\sib(s),\, \e^\bu(s))|^2_\mathcal{A}|( \su_{hs}(s),\, \bu_{hs}(s))|^2_\mathcal{A}\,ds+Ch^{2r}\int\limits_0^t e^{2\al s}\big(\|\su_t\|_{r}+\|\bu_s(s)\|_{r+1}^2+\|p_s(s)\|_{r,\Omega}^2\big).
		\end{align*}
		Using  Lemmas \ref{3.uap0}, \ref{ldgphproj} and \ref{ldgetL2L2err}, along with Theorem \ref{ldgsemierror}  and {\bf{(A2)}}, we  arrive at the required result.		\end{proof}
\noindent
\noindent
	{\bf{{Proof of pressure error estimate in Theorem 3.1.}}}
	An application of (\ref{2.we2}), and (\ref{3.w2}) leads to
	\begin{align}\label{5.sp0}
-d\,(\bphi_h,&\,p_h-r_h (p))= (\bu_{ht}-\bu_t,\,\bphi_h)-\nu(b_1^l(\bphi_h,\,\su_h-\su)+c_1^l(\bu_h-\bu,\,\bphi_h))\nonumber\\
		&+\mo_1^{\bu_h}(\bu_h,\,\bu_h,\,\bphi_h)-\mo_1^\bu(\,\bu,\,\bu,\,\bphi_h)-d\,(\bphi_h,\,p-r_h (p))+(\boldsymbol{\mu}(\theta_h-\theta),\,\bphi_h)\quad \forall \bphi_h \in \bX_h. 
\end{align}
Utilizing the inf--sup condition stated in Lemma \ref{ldginf-sup}, we find that  there is a \( \bphi_h \in \bX_h \), such that
\begin{align}\label{5.sp1}
d(\bphi_h, p_h-r_h (p))=-\|p_h-r_h (p)\|^2, \quad \|\bphi_h\|_\varepsilon\leq \frac{1}{\beta_0}\|p_h-r_h (p)\|.
\end{align}
Applying (\ref{5.sp1}) in 	(\ref{5.sp0}), we observe that 
\begin{align*}
\|p_h-&r_h (p)\|^2=(\bu_{ht}-\bu_t,\,\bphi_h)-\nu(b_1^l(\bphi_h,\,\su_h-\su)+c_1^l(\bu_h-\bu,\,\bphi_h))\nonumber\\
		&+\mo_1^{\bu_h}(\bu_h,\,\bu_h,\,\bphi_h)-\mo_1^\bu(\,\bu,\,\bu,\,\bphi_h)-d\,(\bphi_h,\,p-r_h (p))+(\boldsymbol{\mu}(\theta_h-\theta),\,\bphi_h).
\end{align*}
 Following the arguments presented in Lemma~\ref{4.snl0}, we arrive at
\begin{align}\label{5.sp5}
\|p_h-r_h (p)\|^2&\leq C\big(\|\e^\bu_t\|^2+\|\e^\bu\|^2+\|e^\theta\|^2+|(  \underline{\e}^\sib,\, \e^\bu)|^2_\mathcal{A} +h^2(|\bu_t|^2+|\bu|^2+|p|^2_1)\big).
\end{align}
An application of the triangle inequality, (\ref{5.sp5}), Lemmas~\ref{5.et0}, and Theorem~\ref{ldgsemierror} leads to
\begin{align*}
\|p-p_h\|^2\leq C(t)h^r. 
\end{align*} 
This completes the proof. \hspace{10cm} 
 $\square$
 
 \section{Fully Discrete  Formulation}
\label{ldgsection6}
\se
In this section, we apply the backward Euler method to discretize the semidiscrete LDG Boussinesq system of equations in time, as given in (\ref{3.w1})--(\ref{3.w3}). Consider a uniform division of the time interval \([0, T]\) given by \(\{t_m\}_{m=0}^M\), with time step \(k>0\) and \(t_m = m k\) for all \(0 \le m \le M\), where \(T = M k\) is the final time. For a sequence \(\{\phi^m\}_{m \geq 0}\),  we define \(\bar\partial_t \phi^m\) as
\begin{align*}
\bar\partial_t\phi^m =\frac{1}{k}(\phi^m-\phi^{m-1}),\quad \,\,m\geq { 1},
	 \end{align*}
     where $\phi^m=\phi(t_m)$.\\
We now define the  backward  Euler scheme  for (\ref{3.w1})--(\ref{3.w3}) as \\
	Given $\bU^0$ and  $\Theta^0$, seek $(\Su^m,\,\bU^m,\,\boldsymbol{\Psi}^m,\,\Theta^m,\,\bP^m)_{m\ge 1}\in$ $\Y_h\times\bX_h\times\bZ_h\times Z_h\times M_h$, such that
\begin{align}
a_1^l(\Su^m,\,\tu_h)&+b_1^l(\bU^m,\,\tu_h)=0\quad\quad \forall \tu_h \in \Y_h,\label{6.w1}\\
a_2^l(\boldsymbol{\Psi}^m,\,\biota_h)&+b_2^l(\Theta^m,\,\biota_h)=0\quad\quad \forall \biota_h \in \bZ_h,\label{6.w1th}\\
	(\bar\partial_t\bU^m,\,\bphi_h)--\nu(b_1^l(\bphi_h,\,\Su^m)&-c_1^l(\bU^m,\,\bphi_h))
		+\mo_1^{\bU^m}(\bU^m,\,\bU^m,\,\bphi_h)+d\,(\bphi_h,\,\bP^m)\nonumber\\
		&+(\boldsymbol{\mu} \Theta_h^m,\,\bphi_h)= (\textbf{f}^m,\,\bphi_h)\,\forall \bphi_h\in\bX_h, \label{6.w2}\\
		(\bar\partial_t\Theta^m,\,\psi_h)-\nu(b_2^l(\psi_h,\,\boldsymbol{\Psi}^m)-c_2^l(\Theta^m,\,\psi_h))&+\mo_2^{\Theta^m}(\bU^m,\,\Theta^m,\,\psi_h)= (g^m,\psi_h) \quad\forall \psi_h\in Z_h, \label{6.w2th}\\
		 d\,(\bU^m,\,q_h)
         &= 0 \quad \forall q_h\in M_h, \label{6.w3}
\end{align} 
where $\bU^0=\bu_h(0)=\bP_h^\bu\bu_0$ and $\Theta^0=\theta_h(0)=\mathcal{P}_h^\theta\theta_0.$\\
Equivalently (\ref{6.w1})--(\ref{6.w3}) can be rewritten as follows: \\
 Find $(\Su^m,\,\bU^m,\,\boldsymbol{\Psi}^m,\,\Theta^m,\,\bP^m)_{m\ge 1}\in \Y_h\times\bX^l_h\times\bZ_h\times Z_h\times M_h$, such that
{ \begin{align}
a_1^l(\Su^m,\,&\tu_h)+b_1^l(\bU^m,\,\tu_h)=0,\label{6.we1}\\
a_2^l(\boldsymbol{\Psi}^m,\,&\biota_h)+b_2^l(\Theta^m,\,\biota_h)=0, \label{6.we1th}\\
	(\bar\partial_t\bU^m,\,\bphi_h)-\nu(b_1^l(\bphi_h,\,\Su^m)-&c_1^l(\bU^m,\,\bphi_h))
		+\mo_1^{\bU^m}(\bU^m,\,\bU^m,\,\bphi_h)+d\,(\bphi_h,\,\bP^m)\nonumber\\&\hspace{3cm}+(\boldsymbol{\mu} \Theta_h^m,\,\bphi_h)=(\textbf{f}^m,\,\bphi_h), \label{6.we2}\\
(\bar\partial_t\Theta^m,\,\psi_h)-\nu(b_2^l(\psi_h,\,\boldsymbol{\Psi}^m)-&c_2^l(\Theta^m,\,\psi_h))
		+\mo_2^{\Theta^m}(\bU^m,\,\Theta^m,\,\psi_h)= (g^m,\psi_h),\label{6.we2th}
\end{align}
 for all $(\tu_h,\,\bphi_h,\,\biota_h,\psi_h,\,q_h)\in \Y_h\times\bX^l_h\times\bZ_h\times Z\times M_h$.}
\subsection{\emph{A Priori} Bounds}
In Lemma \ref{6.apU0}, we present {\it a priori} regularity estimates for the solutions of (\ref{6.we1})--(\ref{6.we2th}).\begin{lemma}\label{6.apU0}	
Under the assumption of  {\bf{(A2)}} and
\( 
0 < \alpha < \min \{\frac{\nu}{2\mathcal{N}_\bu C_\bu^1},\,\frac{\kappa}{2\mathcal{N}_\theta C_\theta^1}\}
\)
with  \(\bU^0 = \bP_h^\bu \bu_0,\,  \Theta^0 = \mathcal{P}_h^\theta \theta_0\), we have 
\begin{align}
\|\bU^m\|^2+e^{-2\al t_M}k\sum_{m=1}^Me^{-2\al t_m}|(\Su^m,\,\bU^m)|^2_\mathcal{A}\leq C_T, \quad m=0,1,...,M,\label{6.apU010}\\
\|\Theta^m\|^2+e^{-2\al t_M}k\sum_{m=1}^Me^{-2\al t_m}|(\boldsymbol{\Psi}^m,\,\Theta^m)|^2_\mathcal{A}\leq C_T, \quad m=0,1,...,M. \label{6.apU011}
\end{align}	
where $C_T=C(\nu,\,\kappa,\,\alpha,\,\mathcal{N}_\bu,\, \mathcal{N}_\theta, \,C_\theta^1, C_\bu^1,\,\,T) > 0$,
\end{lemma}	
\begin{proof}
Choose $\psi_h=\Theta^m$ in (\ref{6.we2th})  and  $\biota_h=\kappa\boldsymbol{\Psi}^m$ in (\ref{6.we1th}). Then, sum all the resulting equations, use the fact
\begin{align}
(\bar\partial_t\Theta^m,\,\Theta^m)=\frac{1}{k}(\Theta^m-\Theta^{m-1},\,\Theta^m)&\ge \frac{1}{2k}(\|\Theta^m\|^2-\|\Theta^{m-1}\|^2)=\frac{1}{2}\bar\partial_t\|\Theta^m\|^2,\label{6.apU01}
\end{align}
(\ref{ldg5.5}) and the Cauchy--Schwarz inequality to arrive at
 \begin{align}\label{6.apU1th}
\dbt \|\Theta^m\|^2 +2 \kappa |(\boldsymbol{\Psi}^m,\,\Theta^m)|^2_\mathcal{A}\leq 2 \|g^m\|\|\Theta^m\|.\end{align}
Observe that,
\begin{align}\label{6.apU2th}
\sum_{m=1}^nke^{2\al t_m}\dbt\|\Theta^m\|^2=&-\sum_{m=1}^{n-1}e^{2\al t_m}(e^{2\al k}-1)\|\Theta^m\|^2-e^{2\al k}\|\Theta^0\|^2\nonumber\\
&+e^{2\al t_n}\|\Theta^m\|^2.
\end{align}
A multiplication of (\ref{6.apU1th}) by $k e^{2\al t_m}$,  { a} addition over $m=1$ to $M$ and a  utilization of (\ref{6.apU2th}) yields 
\begin{align}\label{6.apU3th}
e^{2\al t_M}\|\Theta^M\|^2&+\bigg(\kappa-\mathcal{N}_\theta C_\theta^1\frac{(e^{2\al k}-1)}{k}\bigg)k\sum_{m=1}^M  e^{2\al t_m}|(\boldsymbol{\Psi}^m,\Theta^m)|^2_\mathcal{A}\leq e^{2\al k}\|\Theta^0\|^2+\frac{2\mathcal{N}_\theta^2}{\kappa}\sum_{m=1}^M e^{2\al t_m}\|\brf\|^2.
\end{align}
Choose $0<\al \le \frac{\kappa}{2\mathcal{N}_\theta C_\theta^1}$, such that
\begin{align*}
1+\frac{\kappa k}{\mathcal{N}_\theta C_\theta^1}\ge e^{2\al k}.
\end{align*}
Multiply (\ref{6.apU3th}) by $e^{-2\al t_M}$ and utilize \textbf{(A2)} to establish (\ref{6.apU011}).\\
Furthermore, substitute $\bphi_h=\bU^m$ in (\ref{6.we2})  and  $\tu_h=\nu\Su^m$ in (\ref{6.we1}). Then, sum all the resulting equations and utilize  the techniques those applied for (\ref{6.apU011})  to arrive at
\begin{align}\label{6.apU3}
e^{2\al t_M}\|\bU^M\|^2&+\bigg(\nu-\mathcal{N}_\bu C_\bu^1\frac{(e^{2\al k}-1)}{k}\bigg)k\sum_{m=1}^M  e^{2\al t_m}|(\Su^m,\bU^m)|^2_\mathcal{A}\nonumber\\ &\leq e^{2\al k}\|\bU^0\|^2+\frac{2\mathcal{N}_\bu^2}{\nu}\sum_{m=1}^M e^{2\al t_m}\|\brf^m\|^2+C\sum_{m=1}^M  e^{2\al t_m}\|\Theta^m\|.
\end{align}
Pick $0<\al \le  \frac{\nu}{2\mathcal{N}_\bu C_\bu^1}$, such that
\begin{align*}
1+\frac{\nu k}{\mathcal{N}_\bu C_\bu^1}\ge e^{2\al k}.
\end{align*}
Multiply (\ref{6.apU3}) by $e^{-2\al t_M}$ and utilize (\ref{6.apU011}), and \textbf{(A2)} to establish (\ref{6.apU010}).
\end{proof}	
\noindent
\begin{lemma}\label{5.2uinv}
Let  the assumption $\textbf{(A2)}$ holds true.  Then, for every $ 0<t\leq T$, there exists a constant $C(t)>0$, such that 
\begin{align*}
e^{-2\alpha t}\int\limits_0^t(\|\bu_{hss}(s)\|^2_{-1,\,h}\, ds+ \|\theta_{hss}(s)\|^2_{-1,\,h}\,ds) \leq C,
\end{align*}
\begin{align}
\|\theta_{htt}\|_{-1,h}=\sup \bigg\{ \frac{<\theta_{htt},\,\psi_h>}{\|\psi_h\|_\varepsilon},\,\psi_h \in Z_h,\, \psi_h \neq 0 \bigg\},\label{3.uap3}
\end{align}
\end{lemma}  
  \begin{proof}
  Rewrite (\ref{3.uap11th}) as follows:
\begin{align}\label{3.uap19th}
	(\theta_{htt},\,\psi_h)=\kappa(b_2^l(\psi_h,\,\bchi_{ht})-c_2^l(\theta_{ht},\,\psi_h))-\mo_2^{\theta_h}(\bu_{ht},\,\theta_h,\,\psi_h)
	-\mo_2^{\theta_h}(\bu_h,\,\theta_{ht},\,\psi_h)+(g_t,\,\psi_h)\, \forall\psi_h \in Z_h. 
	\end{align}
Utlizing the Cauchy--Schwarz inequality and Lemma \ref{ldgdistrace}, we arrive at
\begin{align}\label{3.uap019th}
   | b_2^l(\psi_h,\,\bchi_{ht})|
    &\leq C \|\psi_{h}\|_\varepsilon |(\bchi_{ht},\theta_{ht})|_\mathcal{A}.
\end{align}
Utilizing   the Cauchy--Schwarz inequality, we have 
{ \begin{align}\label{3.uap119th}
|c_2^l(\theta_{ht},\psi_h)|&=\left|\sum_{e\in\mathcal{N}_h}\int\limits_e\frac{c_{11}}{|e|}[\theta_{ht}]\cdot[\psi_h]\,ds\right| 
    \leq |(\bchi_{ht},\theta_{ht})|_\mathcal{A} \|\psi_h\|_\varepsilon.
\end{align}}
 Applying (\ref{ldgtri2}), Lemma \ref{Enrscript}, (\ref{3.uap019th}),  (\ref{3.uap119th}) in (\ref{3.uap19th}), we find that
\begin{align*}
(\theta_{htt},\,\psi_h)\leq &C(|(\bchi_{ht},\,\theta_{ht})|_\mathcal{A}+|(\bchi_{ht},\,\theta_{ht})|_\mathcal{A}|(\su_{h},\,\bu_{h})|_\mathcal{A}+|(\su_{ht},\,\bu_{ht})|_\mathcal{A}|(\bchi_{h},\,\theta_{h})|_\mathcal{A} +\|g_t\|^2)\|\psi_h\|_\varepsilon.
\end{align*}
Additionally, from Theorem \ref{ldgsemierror}
, we have $|(\bchi_{h},\,\theta_{h})|_\mathcal{A}\leq C$. Therefore,
\begin{align}\label{3.uap21th}
\|\theta_{htt}\|_{-1,h}^2 \leq C(|(\bchi_{ht},\,\theta_{ht})|^2_\mathcal{A}+\|g_t\|^2).
\end{align}
In a similar manner, we can bound
\begin{align}\label{3.uap21}
\|\bu_{htt}\|_{-1,h}^2 \leq C(|(\su_{ht},\,\bu_{ht})|_\mathcal{A}+\|\brf_t\|^2).
\end{align} 
Combine (\ref{3.uap21th}) and (\ref{3.uap21}) and  multiply the resulting inequality by $e^{2\alpha t}$. Integrate from $0$ to $t$ and utilize  (\ref{3.uap15}),  (\ref{3.uap3}). Then,  again multiply by  $e^{-2\alpha t}$ to arrive at 
  the desired result.
 \end{proof}
\subsection{Fully Discrete Error Estimates}
  We now derive fully discrete error estimates for backward  Euler scheme. For a fixed \(m \in \mathbb{N} \) with \( 1 \leq m \leq M \), suppose  \( \underline{\bE}_\sib^m = \Su^m - \su_h(t_m) \), \( \bE_\bu^m = \bU^m - \bu_h(t_m) \), \( \bE_\bchi^m = \boldsymbol{\Psi}^m - \bchi_h(t_m) \) and \( E_\theta^m = \Theta^m - \theta_h(t_m)\). Evaluate (\ref{3.we1})--(\ref{3.we2th}) at \( t = t_m \), then subtract these equations from (\ref{6.we1})--(\ref{6.we2th}) , respectively, to obtain
\begin{align}
a_1^l(\underline{\bE}_\sib^m,\tu_h)+b_1^l(\bE_\bu^m,\tu_h)&=0 \quad \forall \tu_h\in \Y_h,\label{6.apU4}\\
a_2^l(\bE_\bchi^m,\biota_h)+b_2^l(E_\theta^m,\biota_h)& = 0 \quad \forall \biota_h\in \bZ_h,\label{6.apU4th}\\
	(\dbt\bE_\bu^m,\bphi_h)-\nu(b_1^l(\bphi_h,\underline{\bE}_\sib^m)-c_1^l(\boldsymbol{\varphi}_m^\bu,\bphi_h))
		&=(\bu_{ht}^m,\bphi_h)- (\dbt\bu_{h}^m,\bphi_h)+\Lambda_h^\bu(\bphi_h) \quad \forall \bphi_h\in   \bX_h^l, \label{6.apU5}\\
		 (\dbt E_\theta^m,\psi_h)-\kappa(b_2^l(\psi_h,\bE_\bchi^m)-c_2^l(T_m^\theta,\psi_h))
		&=(\theta_{ht}^m,\psi_h)- (\dbt\theta_{h}^m,\psi_h)+\Lambda_h^\theta(\psi_h) \quad \forall \psi_h\in   Z_h, \label{6.apU5th}
\end{align}
where 	$\Lambda_h^\bu(\bphi_h)=\Lambda_{h1}^\bu(\bphi_h)+\Lambda_{h2}^\bu(\bphi_h)$ with 
\begin{equation}\label{6.apU6}
 \left\{
\begin{array}{l}
 \Lam^\bu_{h1}(\bphi_h) =  \mo_1^{\bu_h^m}(\bu_h^m,\,\bu_h^m,\,\bphi_h) -\mo_1^{\bu_h^m}(\bU^m,\bU^m,\,\bphi_h),\\
  \Lam^\bu_{h2}(\bphi_h)=  \mo_1^{\bu_h^m}(\bU^m,\,\bU^m,\,\bphi_h)-\mo_1^{\bU^m}(\bU^m,\,\bU^m,\,\bphi_h)
  \end{array}\right.
  \end{equation}
 and  $\Lambda_h^\theta(\psi_h)=\Lambda_{h1}^\theta(\psi_h)+\Lambda_{h2}^\theta(\psi_h)$ with
  \begin{equation*}
 \left\{
\begin{array}{l}
 \Lambda_{h1}^\theta(\psi_h) =  \mo_2^{\theta_h^m}(\bu_h^m,\,\theta_h^m,\,\psi_h) -\mo_2^{\theta_h^m}(\bU^m,\Theta^m,\,\psi_h),\\
  \Lambda_{h2}^\theta(\psi_h)=  \mo_2^{\theta_h^m}(\bU^m,\Theta^m,\,\psi_h)-\mo_2^{\Theta^m}(\bU^m,\Theta^m,\,\psi_h).
  \end{array}\right.
  \end{equation*}

  \begin{lemma}\label{6.beu0}
Under the assumptions stated in Lemma \ref{6.apU0} together with those of Theorem \ref{ldgsemierror}, it follows that
\begin{align*}
\|\bE_\bu^m\|+\|E_\theta^m\|&+\bigg(\nu e^{-2\al t_M}k\sum_{m=1}^M e^{2\al t_m}|(\underline{\bE}_\sib^m, \bE_\bu^m)|^2_\mathcal{A}\bigg)^{1/2} + \bigg(\kappa e^{-2\al t_M}k\sum_{m=1}^M e^{2\al t_m}|(\bE_\bchi^m, E_\theta^m)|^2_\mathcal{A}\bigg)^{1/2}\leq C_Tk,
\end{align*}
where the constant $C_T=C(\nu,\,\kappa,\,\alpha,\,\mathcal{N}_\bu,\, \mathcal{N}_\theta, \,C_\theta^1, C_\bu^1,\,\,T) > 0$ is independent of $h$ and $k$.
\end{lemma}  
	\begin{proof}
Choose $\biota_h=\kappa \bE_\bchi^m$ in (\ref{6.apU4th}), $\psi_h=E_\theta^m$ in  (\ref{6.apU5}) and sum all the resulting equations to obtain 
\begin{align}
	(\dbt E_\theta^m, E_\theta^m)+\kappa(a_2^l(\bE_\bchi^m,\bE_\bchi^m)+c_2^l(E_\theta^m,E_\theta^m))
		=&(\theta_{ht}^m,\, E_\theta^m)- (\dbt\theta_{h}^m,\, E_\theta^m)+\Lambda_h^{\theta}(E_\theta^m).\label{6.apU05th}
\end{align}
 Then, use (\ref{6.apU01}) in (\ref{6.apU05th}) to arrive at
\begin{align}\label{6.beu1th}
\hspace{-0.5cm}\dbt \| E_\theta^m\|^2+2\kappa |(\bE_\bchi^m,  E_\theta^m)|^2_\mathcal{A} =  2(\theta_{ht}^m, E_\theta^m)- 2(\dbt\theta_{h}^m, E_\theta^m)
 +2\Lambda_h^\theta( E_\theta^m).
\end{align}
Using (\ref{6.apU6}), the nonlinear terms in \(\Lambda_{h1}^\theta(E_\theta^m)\) can be expressed as 
\begin{align}\label{6.beu2}
\Lambda_{h1}^{\theta}( E_\theta^m)= &-\mo_2(\bU^m,\, E_\theta^m,\, E_\theta^m)-\mo_2(\bE_\bu^m,\,\theta_h^m-\theta^m,\, E_\theta^m)-
\mo_2(\bE_\bu^m,\,\theta^m,\, E_\theta^m).
\end{align}
Owing to the positivity property (\ref{ldg5.5}), the first term on the right--hand side of (\ref{6.beu2}) can be omitted.\\
The second term on the right--hand side of (\ref{6.beu2}) can be bounded similarly to \(\mo_2^{\theta_h}(\be^\bu, \omega^\theta, \eta^\theta)\) in Lemma \ref{4.snl0}. Using the Cauchy–Schwarz and Young inequalities together with Theorem \ref{ldgsemierror} gives
\begin{align}\label{6.beu3th}
|--\mo_2(\bE_\bu^m,\theta_h^m--\theta^m, E_\theta^m)|\leq C\|\bE_\bu^m\|^2+{\epsilon}\|E_\theta^m\|^2_\varepsilon.
\end{align}
Apply (\ref{ldgtri1}) and  Cauchy--Schwarz's inequality to estimate the third term in (\ref{6.beu2}) as 
\begin{align}\label{6.beu4th}
|-\mo_2(\bE_\bu^m,\bu^m, E_\theta^m)|\leq C\|\bE_\bu^m\|^2+{\epsilon}\|E_\theta^m\|^2_\varepsilon+ C\|E_\theta^m\|^2+{\epsilon}\|\bE_\bu^m\|^2_\varepsilon.
\end{align}
Utilizing (\ref{ldgtri200}), Young's inequality and Lemmas \ref{ldgl4}, \ref{6.apU0}, we find that
\begin{align}\label{6.beu5th}
&|\mo_2^{\theta_h^m}(\bU^m,\Theta^m,\, E_\theta^m)-\mo_2^{\Theta^m}(\bU^m,\Theta^m,\, E_\theta^m)|
\leq \|\bE_\bu^m\|_{L^4(\Omega)} \|\Theta^m\|_\varepsilon \|E_\theta^m\|_{L^4(\Omega)}\nonumber\\
&\leq \|\bE_\bu^m\|^{1/2}\|\bE_\bu^m\|_\varepsilon^{1/2} \|\Theta^m\|_\varepsilon\|E_\theta^m\|^{1/2}\|E_\theta^m\|^{1/2}_\varepsilon\leq C\|\bE_\bu^m\|\| \|\Theta^m\|_\varepsilon\|E_\theta^m\|_\varepsilon+C\|\bE_\bu^m\|_\varepsilon \|\Theta^m\|_\varepsilon\|E_\theta^m\|.
\end{align}
By applying Young's inequality, we obtain 
\begin{align}\label{6.beu6th}
|\Lam^\theta_{h2`}( E_\theta^m)|\leq C\|\Theta^m\|_\varepsilon^2\|\bE_\bu^m\|^2+{\epsilon}\|E_\theta^m\|^2_\varepsilon+ C\|\Theta^m\|_\varepsilon^2\|E_\theta^m\|^2+{\epsilon}\|\bE_\bu^m\|^2_\varepsilon.
\end{align}
Applying the Taylor expansion to the first term in (\ref{6.beu1th}), we deduce that
\begin{align}\label{6.beu7th}
(\theta_{ht}^m, E_\theta^m)- (\dbt\theta_{h}^m, E_\theta^m)=\frac{1}{k}\int\limits_{t_{m-1}}^{t_m}(t-t_{m-1})(\theta_{htt}(t),  E_\theta^m)dt.
\end{align}
Employing the Cauchy--Schwarz inequality, Young's inequality, along with Lemma \ref{Enrscript}, it follows that  
		\begin{align}\label{6.beu8th}
			2(\theta_{ht}^m, E_\theta^m)- 2(\dbt\theta_{h}^m, E_\theta^m)& \leq C\,k^{1/2}\,\left(\int\limits_{t_{m-1}}^{t_m}\,\|\theta_{htt}(t)\|^2_{-1,h}\,dt\right)^{1/2}\| E_\theta^m\|_\varepsilon\nonumber\\
			&\leq {\epsilon}|(\bE_\bchi^m, E_\theta^m)|^2_\mathcal{A}\,+\,Ck\int\limits_{t_{m-1}}^{t_m}\,\|\theta_{htt}(t)\|^2_{-1,h}\,dt.
		\end{align}
Substitute (\ref{6.beu3th})-(\ref{6.beu8th}) in  (\ref{6.beu1th}), multiply both sides of the resulting inequality by $k e^{2\alpha mk}$, and sum over $m = 1$ to $n$ $(\leq M)$. Next, apply the  techniques those used for  (\ref{6.apU2th})
to find   that
\begin{align}\label{6.beu10th}
e^{2\al t_n}\| E_\theta^n\|^2&+ \kappa k  \sum_{m=1}^{n}e^{2\alpha mk} |(\bE_\bchi^m, E_\theta^m)|^2_\mathcal{A} \leq  \sum_{m=1}^{n-1}ke^{2\al mk}(e^{2\al k}-1)\| E_\theta^m\|^2\nonumber\\
&+Ck\sum_{m=1}^{n}e^{2\alpha mk}(1+\|\Theta^m\|^2_\varepsilon)(\| E_\theta^m\|^2+\| \bE_\bu^m\|^2)+ Ck^2\sum_{m=1}^{n}e^{2\alpha mk}\int\limits_{t_{m-1}}^{t_m}\,\|\theta_{htt}(t)\|^2_{-1,h}\,dt.
\end{align}
Use Lemma \ref{5.2uinv} to find that
\begin{align}\label{6.beu11th}
Ck^2\sum_{m=1}^{n}e^{2\alpha mk}\int\limits_{t_{m-1}}^{t_m}\,\|\theta_{htt}(t)\|^2_{-1,h}\,dt&= Ck^2\sum_{m=1}^{n}\int\limits_{t_{m-1}}^{t_m}e^{2\alpha (t_m-t)}e^{2\alpha t}\,\|\theta_{htt}(t)\|^2_{-1,h}\,dt \nonumber\\
&\leq Ck^2e^{2\alpha k}\sum_{m=1}^{n}\int\limits_{t_0}^{t_m}e^{2\alpha t}\,\|\theta_{htt}(t)\|^2_{-1,h}\,dt  \leq Ck^2e^{2\alpha (n+1)k}.
\end{align}
A use of (\ref{6.beu11th}) in (\ref{6.beu10th}) and the fact $(e^{2\al k}-1)\leq Ck$ yields
\begin{align}
\label{6.beu13th}
e^{2\al t_n}\| E_\theta^n\|^2+ \nu k  \sum_{m=1}^{n}e^{2\alpha mk} |(\bE_\bchi^m,\,  E_\theta^m)|_\mathcal{A}^2 \leq  C \sum_{m=1}^{n}e^{2\al t_m}\|\Theta^m\|_\varepsilon^2\| E_\theta^m\|^2
&+ Ck^2e^{2\alpha (n+1)k}\nonumber\\
+ \epsilon \sum_{m=1}^{n}e^{2\al t_m}|(\underline{\bE}_\sib^m,\, \bE_\bu^m)|_\mathcal{A}^2.
\end{align}
Additionally, substitute $\tu_h=\nu \underline{\bE}_\sib^m$ in (\ref{6.apU4}), $\bphi_h=\bE_\bu^m$ in  (\ref{6.apU5}) and apply the techniques involved for (\ref{6.beu13th}) to arrive at 
\begin{align}
\label{6.beu13}
e^{2\al t_n}\|\bE_\bu^n\|^2+ \nu k  \sum_{m=1}^{n}e^{2\alpha mk} |(\underline{\bE}_\sib^m,\, \bE_\bu^m)|_\mathcal{A}^2 \leq &  C \sum_{m=1}^{n}e^{2\al t_m}(\|\bE_\bu^m\|^2+\|E_\theta^m\|^2 )+ Ck^2e^{2\alpha (m+1)k}.
\end{align}
Now, the desired result follows from  (\ref{6.beu13th}), (\ref{6.beu13})  and the discrete Grönwall's lemma.
\end{proof}
\begin{theorem}\label{bevel}
Under the conditions stated in Theorem \ref{ldgsemierror} together with those of Lemma \ref{6.beu0}, it follows that
	\begin{align*}
	\|\bu(t_m)-\bU^m\|+ \|\theta(t_m)-\Theta^m\|\leq C_T(h^{r+1}\,+\,k),
	\end{align*}
where the constant
$C_T=C(\nu,\,\kappa,\,\alpha,\,\mathcal{N}_\bu,\, \mathcal{N}_\theta, \,C_\theta^1, C_\bu^1,\,\,T) > 0.$	\end{theorem}
	\begin{proof}
The assertion follows immediately from Theorem 3.1 and Lemma \ref{6.beu0}.
	\end{proof}
\begin{lemma}\label{6.dbtu0}
Under the hypotheses of Lemma \ref{6.apU0} together with those of Theorem \ref{ldgsemierror}, it follows that 
	\begin{align*}
	e^{-2\al t_n}\sum_{m=1}^{n}e^{2\al t_n}(\|\dbt\bE_\bu^m\|^2+\|\dbt E_\theta^m\|^2)+\nu|( \underline{\bE}_\sib^m,\,\bE_\bu^m)|^2_\mathcal{A}+\kappa|( \bE_\bchi^n,\,E_\theta^n)|^2_\mathcal{A} \leq C_Tk,
	\end{align*}
	where the constant $C_T=C(\nu,\,\kappa,\,\alpha,\,\mathcal{N}_\bu,\, \mathcal{N}_\theta, \,C_\theta^1, C_\bu^1,\,\,T) > 0$, independent of $h$ and $k$.
	\end{lemma}
\begin{proof}
By repeating the arguments employed in (\ref{6.apU4th}), we obtain
\begin{align}\label{6.dbtu1th}
a_2^l(\bE_\bchi^{m-1},\,\biota_h)+b_2^l(E_\theta^{m-1},\,\biota_h)=0 \quad \forall \biota_h\in \bZ_h.\
\end{align}
A use of (\ref{6.apU4th}) and (\ref{6.dbtu1th}) leads to
\begin{align}\label{6.dbtu2th}
a_2^l(\dbt \bE_\bchi^m,\,\biota_h)+b_2^l(\dbt E_\theta^m,\,\biota_h)=0 \quad \forall \biota_h\in \bZ_h.\
\end{align}
Substitute  $\biota_h=\kappa \bE_\bchi^m$ in (\ref{6.dbtu2th}), and  $\psi_h=\dbt E_\theta^m$ in  (\ref{6.apU5th}). Then, sum  all the resulting equations to obtain 
\begin{align}\label{6.dbtu3th}
\|\dbt E_\theta^m\|^2+\kappa(a_2^l(\bE_\bchi^m,\dbt\bE_\bchi^m )&+c_2^l({ E_\theta^m},\,\dbt E_\theta^m)) =2(\theta_{ht}^m,\,\dbt E_\theta^m)- 2(\dbt\theta_{h}^m,\,\dbt E_\theta^m)+2\Lambda_h^\theta(\dbt E_\theta^m).
\end{align}
We reformulate the nonlinear part of (\ref{6.dbtu3th}) as 
\begin{align} \label{6.dbtu4th}
\Lambda_h^\theta(\dbt E_\theta^m)=&\Lam^\theta_{h1}(\dbt E_\theta^m)+\Lam^\theta_{h2}(\dbt E_\theta^m)=-\mo_2^{\theta_h^m}( \bE_\bu^m,\, E_\theta^m,\,\dbt E_\theta^m)\nonumber\\
&-\mo_2^{\theta_h^m}(\bu_h^m,\, E_\theta^m,\,\dbt E_\theta^m)
-\mo_2^{\theta_h^m}(\bE_\bu^m,\,\bu_h^m,\,\dbt E_\theta^m)+\Lam^\theta_{h2}(\dbt E_\theta^m).
\end{align}
The first term on the right hand side of (\ref{6.dbtu4th}) can be written as 
\begin{align*}
-\mo_2^{\theta_h^m}( \bE_\bu^m,\, E_\theta^m,\,\dbt E_\theta^m)=-\frac{1}{k}\mo_2^{\bu_h^m}(\bE_\bu^m,\, E_\theta^m,\, E_\theta^m)+ \frac{1}{k}\mo_2^{\theta_h^m} (\bE_\bu^m,\, E_\theta^m,\, E_\theta^{m-1}).
\end{align*}
Applying (\ref{ldg5.5}) and  (\ref{ldgtri2}), we find that 
\begin{align}\label{5.4pre0th}
-\mo_2^{\theta_h^m}( \bE_\bu^m,&\, E_\theta^m,\,\dbt E_\theta^m)\leq \frac{1}{k}\mo_2^{\theta_h^m}( \bE_\bu^m,\, E_\theta^m,\, E_\theta^{m-1})  \leq \frac{C}{k}|(\bE_\bchi^m,\, E_\theta^m)|_\mathcal{A}|(\bE_\bchi^{m},\, E_\theta^{m})|_\mathcal{A} |(\bE_\bchi^{m-1},\, E_\theta^{m-1})|_\mathcal{A}. 
\end{align}
An utilization of Lemmas \ref{ldgsobolev2.1}, \ref{ldgTrilinear} and Theorem \ref{ldgsemierror} leads to
\begin{align}\label{5.4pre1th}
|-\mo_2^{\theta_h^m}(\bu_h^m,\, E_\theta^m,\,\dbt E_\theta^m)-\mo_2^{\theta_h^m}(\bE_\bu^m,\,\theta_h^m,\,\dbt E_\theta^m)|
\leq \epsilon \|\dbt E_\theta^m\|^2 + C  |(\bE_\bchi^m,\, E_\theta^m)|^2_\mathcal{A}. 
\end{align}
As a consequence of  $ [\theta(t_m)]=0$, the difference of the nonlinear term 
$\mo_2^{\theta_h^m}(\bU^m,\theta^m,\,\bar\partial_t E^m_\theta)-\mo_2^{\Theta^m}(\bU^m,\theta^m,\,\bar\partial_tE^m_\theta)=0$.
	 	 We subtract this zero term with $\theta(t_m)$ from $\Lam_h^{\theta 2}(\bar\partial_t E^m_\theta)$ and use Lemma \ref{ldgTrilinear} to arrive at
		\begin{align*}
			&|\Lam_h^{\theta 2}(\bar\partial_t\bE^m_\bu)|
		\,	\leq \,C \|\bu_h(t_m)-\bU_h^{m}\|_{L^4(\Omega)}\,\,\|\Theta_h^m-\theta(t_m)\|_\varepsilon\,\,
			\|\bar\partial_t E^m_\theta\|_{L^4(\Omega)}.
			\end{align*}
			Using the triangle inequality in conjunction with (\ref{ldgLp}) and Lemmas \ref{ldgsobolev2.1} and \ref{ldgdistrace}, we deduce
			\begin{align*}
			|\Lam_h^{\theta 2}&(\bar\partial_t E^m_\theta)| \leq  \,C \|\bu_h(t_m)-\bU_h^{m}\|_{L^4(\Omega)}\,\,\|\theta_h(t_m)-\theta(t_m)\|_\varepsilon\,\,\|\bar\partial_t E^m_\theta\|_{L^4(\Omega)}\\
			&+C \|\bu_h(t_m)\,-\bU_h^{m}\|_{L^4(\Omega)}\,\,\| E^m_\theta\|_\varepsilon\,\,\|\bar\partial_t E^m_\theta\|_{L^4(\Omega)}\\
			 \leq & C\|\bu_h(t_m)\,-\bU_h^{m}\|_\varepsilon\,\,\|\theta_h(t_m)-\theta(t_m)\|_\varepsilon\,\, \left(\frac{1}{\displaystyle\min_{T\in \mathcal{T}_h}h_T^{1/2}}\|\bar\partial_t E^m_\theta\|_{L^2(\Omega)}\right)\\
			&+C\|\bu_h(t_m)\,-\bU_h^{m}\|^{1/2}\|\bu_h(t_m)-\bU_h^{m}\|_\varepsilon^{1/2}\| E^m_\theta\|_\varepsilon \|\bar\partial_t E^m_\theta\|_{L^2(\Omega)}^{1/2}\|\bar\partial_t E^m_\theta\|_\varepsilon^{1/2}.
		\end{align*}
		An application of  Lemmas \ref{3.uap0}, \ref{6.apU0} and Theorem \ref{ldgsemierror}  yields
\begin{align*}
			|\Lam_h^{\theta 2}(\bar\partial_t E^m_\theta)| &\leq  C\|\bu_h(t_m)-\bU_h^{m}\|_\varepsilon\|\bar\partial_t E^m_\theta\|  +\frac{C}{k^{1/2}} \|\bu_h(t_m)-\bU_h^{m}\|_\varepsilon^{1/2}\| E^m_\theta\|_\varepsilon(\| E^m_\theta\|_\varepsilon^{1/2}+\|  E_\theta^{m-1}\|_\varepsilon^{1/2}).
		\end{align*}		
Applying  triangle  and Young's inequalities, we find that
		\begin{align}\label{5.4pre2th}
			| \Lambda_{h2}^\bu(\bar\partial_t \bE^m_\bu)|
			\leq  \epsilon \|\bar\partial_t \bE^m_\bu\|_{L^2(\Omega)}^2+C\|  E_\theta^{m}\|_\varepsilon^2+C\|  \bE_\bu^{m}\|_\varepsilon^2+\frac{C}{k}(\| E_\theta^{m-1}\|_\varepsilon^2+\| E^m_\theta\|_\varepsilon^2).
	 	\end{align}
%
By utilizing (\ref{6.beu7th}), Cauchy--Schwarz's inequality and  Lemma \ref{Enrscript}, we arrive at
\begin{align}\label{6.dbtu8th}
			(\theta_{ht}^m,\,\dbt E_\theta^m)- (\dbt\theta_{h}^m,\,\dbt E_\theta^m)&\leq C\,k^{1/2}\,\bigg(C\int\limits_{t_{n-1}}^{t_n} \|\theta_{hss}(s)\|^2_{-1,h}\,ds \bigg)^{1/2}|(\dbt \bE_\bchi^m,\dbt E_\theta^m)|_\mathcal{A}\nonumber\\
      &\leq  C\int\limits_{t_{n-1}}^{t_n} \|\theta_{hss}(s)\|^2_{-1,h}\,ds + \frac{C}{k}  |( \bE_\bchi^m-\bE_\bchi^{m-1},\, E_\theta^m- E_\theta^{m-1})|^2_\mathcal{A}.
		\end{align}
		Applying symmetricity of ${ a_2^l}(\cdot,\cdot)$, we obtain
		\begin{align}\label{5.4pre3th}
			a_2^l(\bE_\bchi^m,\,\dbt \bE_\bchi^m )=& \frac{1}{2}\bigg(\frac{1}{k}a_2^l( \bE_\bchi^m,\, \bE_\bchi^m)-\frac{1}{k}a_2^l( \bE_\bchi^{m-1},\, \bE_\bchi^{m-1})+k a_2^l(\dbt \bE_\bchi^m,\,\dbt \bE_\bchi^m)\bigg).
		\end{align}
		Further,
		\begin{align}\label{5.4pre4th}
			\sum_{m=1}^n  e^{2\al t_m}&\bigg(a_2^l( \bE_\bchi^m,\, \bE_\bchi^m)+c_2^l(E^m_\theta, E^m_\theta)-(a_2^l( \bE_\bchi^{m-1},\, \bE_\bchi^{m-1})+c_2^l(E^{m-1}_\theta, E^{m-1}_\theta))\bigg)\nonumber\\
			=&e^{2\al t_n}(a_2^l( \bE_\bchi^n,\, \bE_\bchi^n)+c_2^l(E^n_\theta, E^n_\theta))-\sum_{m=1}^{n-1} e^{2\al t_m}(e^{2\al k}-1)(a_2^l( \bE_\bchi^m,\, \bE_\bchi^m)++c_2^l(E^m_\theta, E^m_\theta)).
		\end{align}		
Combinine (\ref{5.4pre0th})--(\ref{5.4pre4th}), multiply   (\ref{6.dbtu3th}) by $ke^{2\alpha t_m}$, sum over $m=1$ to $n(\leq M)$ to arrive at 
\begin{align}\label{5.4pre5th}
k\sum_{m=1}^n &\|\dbt E_\theta^m\|^2 + \nu e^{2\alpha t_n}|(\bE_\bchi^n,\, E_\theta^n)|_\mathcal{A}\leq k\sum_{m=1}^{n-1} e^{2\al t_m}(a_2^l( \bE_\bchi^m,\, \bE_\bchi^m)+c_2^l(E^m_\theta, E^m_\theta))+Ck \sum_{m=1}^{n} e^{2\al t_m} |(\bE_\bchi^m,\, E^m_\theta)|^2_\mathcal{A}\nonumber\\
&
+\frac{C}{k}k \sum_{m=1}^{n-1} e^{2\al t_m} |(\bE_\bchi^m,\, E^m_\theta)|^2_\mathcal{A}|(\bE_\bchi^{m-1},\, E^{m-1}_\theta)|^2_\mathcal{A}
+ \frac{C}{k}k \sum_{m=1}^{n-1} e^{2\al t_m} (|(\bE_\bchi^m,\, E^m_\theta)|^2_\mathcal{A}+|(\bE_\bchi^{m-1},\, E^{m-1}_\theta)|^2_\mathcal{A})\nonumber\\
&+Ck \sum_{m=1}^{n} e^{2\al t_m}\int\limits_{t_{n-1}}^{t_n} \|\theta_{hss}(s)\|^2_{-1,h}\,ds+Ck \sum_{m=1}^{n} e^{2\al t_m} |(\underline{\bE}_\sib^m,\, \bE^m_\bu)|^2_\mathcal{A}.
\end{align}
Applying Lemmas \ref{ldgcont}, \ref{3.uap0}, \ref{5.2uinv} and \ref{6.beu0} in (\ref{5.4pre5th}), we obtain
\begin{align}\label{5.4pre5th0}
k\sum_{m=1}^n \|\dbt E_\theta^m\|^2 &+ \nu e^{2\alpha t_n}|(\bE_\bchi^n,\, E_\theta^n)|_\mathcal{A}\leq \frac{C}{k}k \sum_{m=1}^{n-1} e^{2\al t_m} (|(\bE_\bchi^m,\, \bE^m_\bu)|^2_\mathcal{A}+|(\bE_\bchi^{m-1},\, \bE^{m-1}_\bu)|^2_\mathcal{A})+ Cke^{2 \alpha t_{n+1}}.
\end{align}
Furthermore, following the steps utilized in (\ref{6.apU4}), we obtain
\begin{align}\label{6.dbtu1}
a_1^l(E_\su^{m-1},\,\tu_h)+b_1^l(\bE_\bu^{m-1},\,\tu_h)=0 \quad \forall \tu_h\in \Y_h.\
\end{align}
A use of (\ref{6.apU4}) and (\ref{6.dbtu1}) leads to
\begin{align}\label{6.dbtu2}
a_1^l(\dbt \underline{\bE}_\sib^m,\,\tu_h)+b_1^l(\dbt\bE_\bu^m,\,\tu_h)=0 \quad \forall \tu_h\in \Y_h.\
\end{align}
Substitute  $\tu_h=\nu \underline{\bE}_\sib^m$ in (\ref{6.dbtu2}), and  $\bphi_h=\dbt\bE_\bu^m$ in  (\ref{6.apU5}). Then, sum  all the resulting equations to obtain 
\begin{align}\label{6.dbtu3}
\|\dbt\bE_\bu^m\|^2+\nu(a_1^l(\underline{\bE}_\sib^m,\dbt\bE_\su^m )&+c_1^l(\bT_m^\bu,\,\dbt\bE_\bu^m)) =  2(\bu_{ht}^m,\,\dbt\bE_\bu^m)- 2(\dbt\bu_{h}^m,\,\dbt\bE_\bu^m)\nonumber\\
&+2\Lambda_h(\dbt\bE_\bu^m)+(\boldsymbol{\mu}E_\theta^m,\,\dbt\bE_\bu^m).
\end{align}
Apply those steps utilized for (\ref{5.4pre5th0}) to arrive at
\begin{align}\label{5.4pre5th0vel}
k\sum_{m=1}^n \|\dbt\bE_\bu^m\|^2 +& \nu e^{2\alpha t_n}|(\underline{\bE}_\sib^n,\,\bE_\bu^n)|_\mathcal{A}\leq \frac{C}{k}k \sum_{m=1}^{n-1} e^{2\al t_m} (|(\underline{\bE}_\sib^m,\, \bE^m_\bu)|^2_\mathcal{A}+|(\underline{\bE}_\sib^{m-1},\, \bE^{m-1}_\bu)|^2_\mathcal{A})+ Cke^{2 \alpha t_{m+1}}.
\end{align}
Finally, an application of (\ref{5.4pre5th0}), (\ref{5.4pre5th0vel}),     the discrete Gro\"nwall's inequality and Lemma \ref{6.beu0} leads to the desired estimate. This completes the proof.
\end{proof}	
\begin{lemma}\label{5.5neget}
Assuming the hypotheses of Lemma \ref{6.beu0} hold. Then, 
		\begin{align*}
			\|\bar\partial_t\bE_\bu^m\|_{-1,h}\leq  C_T k^{1/2},
		\end{align*}
		 where $C_T=C(\nu,\,\kappa,\,\alpha,\,\mathcal{N}_\bu,\, \mathcal{N}_\theta, \,C_\theta^1, C_\bu^1,\,\,T) > 0$ is a constant. 
	\end{lemma}

\begin{proof}
Next, we express the nonlinear terms of (\ref{6.dbtu3}) as
\begin{align*} 
\Lambda_h^\bu(\bphi_h)=\Lambda_{h1}^\bu(\bphi_h)+\Lambda_{h2}^\bu(\bphi_h)=&-\mo_2^{\bu_h^m}(\bE_\bu^m,\,\bE_\bu^m,\,\bphi_h)-\mo_2^{\bu_h^m}(\bu_h^m,\,\bE_\bu^m,\,\bphi_h)\nonumber\\
&-\mo^{\bu_h^m}(\bE_\bu^m,\,\bu_h^m,\,\bphi_h)+\Lambda_{h2}^\bu(\bphi_h).
\end{align*}
Applying triangle inequality, Lemma \ref{ldgTrilinear} and Theorem \ref{ldgsemierror}, we have 
\begin{align}\label{5.5preerr1}
\Lambda_{h1}^\bu(\bphi_h)+\Lambda_{h2}^\bu(\bphi_h)\leq C|(\underline{\bE}_\sib^m,\,\bE_\bu^m)|^2_\mathcal{A}\|\bphi_h\|_\varepsilon+ C|(\underline{\bE}_\sib^m,\,\bE_\bu^m)|_\mathcal{A}\|\bphi_h\|_\varepsilon.
\end{align}
Using Lemma \ref{ldgcont}, it follows that 
\begin{align}\label{6.dbtu7}
|\nu(a_1^l(\underline{\bE}_\sib^m,\,\dbt \underline{\bE}_\sib^m )+c_1^l(\bE_\bu^m,\,\dbt\bE_\bu^m))|\leq C |(\dbt \underline{\bE}_\sib^m,\,\dbt\bE_\bu^m)|_\mathcal{A} |( \underline{\bE}_\sib^m,\, \bE_\bu^m)|_\mathcal{A}.
\end{align}
Employing (\ref{6.beu7th}) in combination with the Cauchy--Schwarz inequality and Lemma \ref{Enrscript}, we find that
\begin{align}\label{5.5.dbtu8}
			(\bu_{ht}^m,\,\bphi_h)- (\dbt\bu_{h}^m,\,\bphi_h)
			&
            \leq C\,k^{1/2}\,\bigg(C\int\limits_{t_{m-1}}^{t_m} \|\bu_{hss}(s)\|^2_{-1,h}\,ds \bigg)^{1/2}\|\bphi_h\|_\varepsilon \leq  C k^{1/2}\|\bphi_h\|_\varepsilon.
		\end{align}
		Combining all the bounds (\ref{5.5preerr1})--(\ref{5.5.dbtu8}) in (\ref{6.apU5}) and applying the definition of $\|\cdot\|_{-1,h}$, Lemma \ref{6.dbtu0}, we establish our desired result. This completes the rest of the proof. 
\end{proof}
\begin{lemma}\label{6.bePressure0}
Let the assumptions of   Lemma \ref{6.apU0} and Theorem \ref{ldgsemierror} hold true. Then, 
\begin{align*}
\|P^m-p_h^m\|\leq C_Tk^{1/2} \qquad \text{for} \qquad m = 1, 2, \ldots, M,
\end{align*}
where  the constant $C_T=C(\nu,\,\kappa,\,\alpha,\,\mathcal{N}_\bu,\, \mathcal{N}_\theta, \,C_\theta^1, C_\bu^1,\,\,T) > 0$ independent of $h$ and $k$.
\end{lemma}	
\begin{proof}
Substracting (\ref{3.w2}) from (\ref{6.w2}), we arrive at
\begin{align*}
d\,(\bphi_h,\,\bP^m -p_h^m)=&(\bu_{ht}^m,\bphi_h)-(\dbt\bu_h^m,\bphi_h)-\nu(b_1^l(\bphi_h,\underline{\bE}_\sib^m)-c_1^l(\bE_\bu^m,\bphi_h))	
\nonumber\\
&-(\bar\partial_t\bE_\bu^m,\bphi_h)
		+\Lambda_h^\bu(\bphi_h)-(\boldsymbol{\mu}E_\theta^m,\,\bphi_h)\qquad \forall \bphi_h\in\bX_h.
\end{align*}
Repeating the analysis of Lemma \ref{6.dbtu0} and employing Lemma \ref{ldginf-sup}, we obtain
{ \begin{align*}
\hspace{-0.5cm}\|P^m-p_h^m\| \leq C(\|\dbt\bE_\bu^m\|&+\kappa|(\dbt \underline{\bE}_\sib^m,\,\dbt\bE_\bu^m)|_\mathcal{A} + |( \underline{\bE}_\sib^m,\,\bE_\bu^m)|_\mathcal{A}+\|E_\theta^m\|+k).
\end{align*}
}
The result now follows immediately from Lemmas \ref{6.beu0}--\ref{5.5neget}, which completes the proof.
\end{proof}
\noindent
An immediate consequence of Lemma \ref{6.bePressure0} and Theorem \ref{ldgsemierror} is the following pressure error estimate.
\begin{theorem}\label{FullyP}
Suppose that the assumptions of Lemma \ref{6.apU0} and Theorem \ref{ldgsemierror} are satisfied. Then, for $m=1,2,\ldots,M$, 
 \begin{align*}
	\|p(t_m)-P^m\| \leq  C_T(h^r+k^{1/2}) 
	\end{align*}
where  the constant $C_T=C(\nu,\,\kappa,\,\alpha,\,\mathcal{N}_\bu,\, \mathcal{N}_\theta, \,C_\theta^1, C_\bu^1,\,\,T) > 0$ independent of $h$ and $k$.
	\end{theorem}
\noindent    
\section{Numerical Experiments}\label{ldgsection7}
This section presents numerical simulations for the Boussinesq system to illustrate and verify the theoretical results established in Theorems \ref{bevel} and \ref{FullyP}. 
The spatial discretization is carried out using the finite element method, while the temporal discretization is performed by backward Euler method.
The computational domain is chosen as $\Omega = [0,1]\times[0,1]$, which corresponds to a two--dimensional square region. 
For the approximation of the stress tensor $\su$, velocity $\bu$, and pressure $p$, we employ the mixed finite element spaces $P_0\text{-}P_1\text{-}P_0$ and $P_1\text{-}P_2\text{-}P_1$. For the temperature gradient $\bchi$, temperature  \( \theta \), we adopt finite element spaces \(P_0\text{-} P_1 \) and \(P_1\text{-}P_2\), depending on the desired accuracy and refinement level.  
The time step $k$ is chosen as $\mathbb{O}(h^{r+1})$, where $h$ denotes the spatial discretization size, the parameter r is equal to 1 or 2, $\kappa=1$. 
The constants are chosen as
$\textbf{C}_{12} = \textbf{0}$, $\textbf{D}_{12} = \textbf{0}$, and $\textbf{E}_{12} = \textbf{0}$ \,  (\cite{cockburn02}). Note that $\underline{\e}_\sib^M=\su(t_M)-\Su^M,\,\e_\bu^M=\bu(t_M)-\bU^M$, ${\e}_{\bchi}^M=\bchi(t_M)-\boldsymbol{\Psi}^M,\, e_\theta^M=\theta(t_M)-\Theta^M$, $e_p^M=p(t_M)-P^M$.
\begin{example}\label{ex1ldg}
To verify the numerical scheme, $\mathbf{f}$ and $g$ are chosen so that the exact solution $(\bu,p,\theta)=((u_1,u_2),p,\theta)$ is obtained as
\begin{align*}
		&u_1 = 2e^t  x^2 (x-1)^2 y (y-1) (2y-1), u_2= -2 e^tx(x-1)(2x-1)y^2(y-1)^2, p=2(x-y)e^t,\\
&\theta = e^t 2 x^2 (x-1)^2 y (y-1) (2y-1)-2 e^tx(x-1)(2x-1)y^2(y-1)^2.
		\end{align*}
	\end{example}
\noindent
Tables \ref{ldgpolyP1Vnu1}-\ref{ldgpolyP2Tnu1} present the computed fully discrete errors and corresponding convergence rates for various mixed finite element combinations, stress tensor--velocity--pressure and temperature gradient--temperature, \(P_0\text{-}P_1\text{-}P_0\),  \(P_0\text{-}P_1\) and  \(P_1\text{-}P_2\text{-}P_1\),  \(P_1\text{-}P_2\), all evaluated with the viscosity parameter \(\nu = 1\).
Similarly, Tables \ref{ldgpolyP1Vnu01}--\ref{ldgpolyP2Tnu01} report results for the same element pairs, now considering a lower viscosity value of \(\nu =0.1\), respectively. For a more diffusive scenario, Tables \ref{ldgpolyP1Vnu001}--\ref{ldgpolyP2Tnu001} display the numerical outcomes for the above mixed finite element  spaces under the setting \(\nu = 0.01\). Tables \ref{ldgpolyP1Vnu0001}--\ref{ldgpolyP2Tnu0001} continue this trend, exploring the performance of these finite element methods at a very small viscosity level of \(\nu = 0.001\).
In all cases, the penalty parameters were set to \(c_{11} = c_{22} = 50\). It is important to emphasize that the data presented across Tables \ref{ldgpolyP1Vnu1}--\ref{ldgpolyP2Tnu0001} confirm the theoretical convergence predictions established in Theorems \ref{bevel} and \ref{FullyP} for the final time  $T = 0.25$. 
  {\small
		\renewcommand{\arraystretch}{1.3}
		\begin{table}[ht!]
			\centering
			\caption{
			Numerical errors and corresponding convergence rates for the stress tensor, velocity, and pressure with $P_0$--$P_1$--$P_0$ element for $\nu=1$ in Example \ref{ex1ldg}.}	
		\label{ldgpolyP1Vnu1}
			\vspace{.4cm}
			\begin{tabular}{|c|c| c| c| c| c| c|}
				
				\hline
				$h$ & $\|\e_\bu^M\|$ & order & $|(\underline{\e}_\sib^M,\,\e_\bu^M)|_\mathcal{A}$ &order &$\|e_p^M\|$ & order\\ 
				\hline
\hline 
 1/4 & 0.002728 &  & 0.065534 &  & 0.105901 & \\
 1/8 & 0.001872 & 0.5430 & 0.045605 & 0.5230 & 0.078561 & 0.4308\\
 1/16 & 0.000814 & 1.2009 & 0.029176 & 0.6444 & 0.060763 & 0.3706\\
 1/32 & 0.000260 & 1.6464 & 0.016583 & 0.8151 & 0.037406 & 0.6999\\
 1/64 & 0.000071 & 1.8712 & 0.008743 & 0.9234 & 0.020308 & 0.8812\\ 
	         		\hline
	           
			\end{tabular}
			\vspace{.1cm}
		\end{table}
	}

{\small
		\renewcommand{\arraystretch}{1.3}
		\begin{table}[ht!]
			\centering
			\caption{
Numerical errors and corresponding convergence rates for the temperature gradient and temperature with the $P_0$--$P_1$ element for $\nu=1$ in Example \ref{ex1ldg}.}	
		\label{ldgpolyP1Tnu1}
			\vspace{.4cm}
			\begin{tabular}{|c|c| c| c| c| c| c|}
				
				\hline
				$h$ & $\|e_\theta^M\|$ & order & $|(\e_{\bchi}^M,\, e_\theta^M)|_\mathcal{A}$ &order\\ 
\hline 
\hline
 1/4 & 0.000379 & & 0.023196 & \\
 1/8 & 0.000092 & 2.0338 & 0.011742 & 0.9828\\
 1/16 & 0.000023 & 1.9743 & 0.005911 & 0.9903\\
 1/32 & 0.000006 & 1.9851 & 0.002962 & 0.9967\\
 1/64 & 0.000001 & 1.9969 & 0.001482 & 0.9987\\
	         		\hline
	           
			\end{tabular}
			\vspace{.1cm}
		\end{table}
	}

  {\small
		\renewcommand{\arraystretch}{1.3}
		\begin{table}[ht!]
			\centering
			\caption{
			Numerical errors and corresponding convergence rates for the stress tensor, velocity, and pressure with $P_1$--$P_2$--$P_1$ element for $\nu=1$ in Example \ref{ex1ldg}.}		
		\label{ldgpolyP2Vnu1}
			\vspace{.4cm}
			\begin{tabular}{|c|c| c| c| c| c| c|}
				
				\hline
				$h$ & $\|\e_\bu^M\|$ & order & $|(\underline{\e}_\sib^M,\,\e_\bu^M)|_\mathcal{A}$ &order &$\|e_p^M\|$ & order\\ 
\hline 
\hline
 1/4 & 0.0005514 & & 0.0268531 & & 0.0274516 & \\ 
 1/8 & 0.0000586 & 3.2332 & 0.008108 & 1.7275 & 0.007855 & 1.8052\\ 
 1/16 & 0.0000056 & 3.3715 & 0.002167 & 1.9034 & 0.002052 & 1.9363\\ 
 1/32 & 0.0000006 & 3.2046 & 0.000555 & 1.9652 & 0.000518 & 1.9846\\ 
	         		\hline
	           
			\end{tabular}
			\vspace{.1cm}
		\end{table}
	}

{\small
		\renewcommand{\arraystretch}{1.3}
		\begin{table}[ht!]
			\centering
			\caption{
Numerical errors and corresponding convergence rates for the temperature gradient and temperature with $P_1$--$P_2$ element for $\nu=1$ in Example \ref{ex1ldg}.}	
		\label{ldgpolyP2Tnu1}
			\vspace{.4cm}
			\begin{tabular}{|c|c| c| c| c| c| c|}
				
				\hline
				$h$ & $\|e_\theta^M\|$ & order & $|(\e_{\bchi}^M,\, e_\theta^M)|_\mathcal{A}$ &order\\ 
\hline 
\hline
 1/4 & 0.0000690 &  & 0.012456 & \\ 
 1/8 & 0.0000082 & 3.0566 & 0.003388 & 1.8781\\ 
 1/16 & 0.0000008 & 3.2125 & 0.000866 & 1.9679\\ 
 1/32 & 0.0000001& 3.1316 & 0.000217 & 1.9911\\ 
	         		\hline
	           
			\end{tabular}
			\vspace{.1cm}
		\end{table}
	}

	  {\small
		\renewcommand{\arraystretch}{1.3}
		\begin{table}[ht!]
			\centering
			\caption{
Numerical errors and corresponding convergence rates for the stress tensor, velocity, and pressure with $P_0$--$P_1$--$P_0$ element for $\nu=0.1$ in Example \ref{ex1ldg}.}	
		\label{ldgpolyP1Vnu01}
			\vspace{.4cm}
			\begin{tabular}{|c|c| c| c| c| c| c|}
				
				\hline
				$h$ & $\|\e_\bu^M\|$ & order & $|(\underline{\e}_\sib^M,\,\e_\bu^M)|_\mathcal{A}$ &order &$\|e_p^M\|$ & order\\ 
\hline 
\hline
 1/4 & 0.005722 &  & 0.414158 &  & 0.034825 & \\ 
 1/8 & 0.001794 & 1.6729 & 0.200193 & 1.0487 & 0.022436 & 0.6342\\ 
 1/16 & 0.000644 & 1.4765 & 0.098622 & 1.0214 & 0.013179 & 0.7675\\ 
 1/32 & 0.000195 & 1.7220 & 0.0490 & 1.0082 & 0.007193 & 0.8735\\ 
 1/64 & 0.000052 & 1.8966 & 0.0244 & 1.0050 & 0.003746 & 0.9411\\ 
	         		\hline
	           
			\end{tabular}
			\vspace{.1cm}
		\end{table}
	}

{\small
		\renewcommand{\arraystretch}{1.3}
		\begin{table}[ht!]
			\centering
			\caption{
Numerical errors and corresponding convergence rates for the temperature gradient and temperature with $P_0$-$P_1$ element for $\nu=0.1$ in Example \ref{ex1ldg}.}	
		\label{ldgpolyP1Tnu01}
			\vspace{.4cm}
			\begin{tabular}{|c|c| c| c| c| c| c|}
				
				\hline
				$h$ & $\|e_\theta^M\|$ & order & $|(\e_{\bchi}^M,\, e_\theta^M)|_\mathcal{A}$ &order\\ 
\hline 
\hline
 1/4 & 0.000379 & & 0.023196 & \\ 
 1/8 & 0.000092 & 2.0339 & 0.011742 & 0.9821\\ 
 1/16 & 0.000023 & 1.9746 & 0.005911 & 0.9903\\ 
 1/32 & 0.000005 & 1.9854 & 0.00296207 & 0.9966\\ 
 1/64 & 0.000001 & 1.9971 & 0.001482 & 0.9987\\ 
	         		\hline
	           
			\end{tabular}
			\vspace{.1cm}
		\end{table}
	}

		  {\small
		\renewcommand{\arraystretch}{1.3}
		\begin{table}[ht!]
			\centering
			\caption{
Numerical errors and corresponding convergence rates for the stress tensor, velocity, and pressure with $P_1$--$P_2$--$P_1$ element for $\nu=0.1$ in Example \ref{ex1ldg}.
}	
		\label{ldgpolyP2Vnu01}
			\vspace{.4cm}
			\begin{tabular}{|c|c| c| c| c| c| c|}
				
				\hline
				$h$ & $\|\e_\bu^M\|$ & order & $|(\underline{\e}_\sib^M,\,\e_\bu^M)|_\mathcal{A}$ &order &$\|e_p^M\|$ & order\\ 
\hline 
\hline
 1/4 & 0.0004256 &  & 0.027546 &  & 0.0028337 & \\ 
 1/8 & 0.0000501 & 3.0853 & 0.0081601 & 1.7551 & 0.0007879 & 1.8466\\ 
 1/16 & 0.0000052 & 3.2559 & 0.002171 & 1.9098 & 0.0002052 & 1.9403\\ 
 1/32 & 0.0000006 & 3.1197 & 0.0005554& 1.9669 & 0.0000518 & 1.9849\\ 	
	         		\hline
	           \end{tabular}
			\vspace{.1cm}
		\end{table}
	}

{\small
		\renewcommand{\arraystretch}{1.3}
		\begin{table}[ht!]
			\centering
			\caption{
Numerical errors and corresponding convergence rates for the temperature gradient and temperature with $P_1$--$P_2$ element for $\nu=0.1$ in Example \ref{ex1ldg}.}	
		\label{ldgpolyP2Tnu01}
			\vspace{.4cm}
			\begin{tabular}{|c|c| c| c| c| c| c|}
				
				\hline
				$h$ & $\|e_\theta^M\|$ & order & $|(\e_{\bchi}^M,\, e_\theta^M)|_\mathcal{A}$ &order\\ 
\hline 
\hline
 1/4 & 0.0000690 & & 0.012456 & \\ 
 1/8 & 0.0000082 & 3.0566 & 0.0033884 & 1.8781\\ 
 1/16 & 0.0000008 & 3.2125 & 0.0008661 & 1.9679\\ 
 1/32 & 0.0000001 & 3.1316 & 0.0002178 & 1.9911\\  
	         		\hline
	           
			\end{tabular}
			\vspace{.1cm}
		\end{table}
	}

	  {\small
		\renewcommand{\arraystretch}{1.3}
		\begin{table}[ht!]
			\centering
			\caption{
Numerical errors and corresponding convergence rates for the stress tensor, velocity, and pressure with $P_0$--$P_1$--$P_0$ element for $\nu=0.01$ in Example \ref{ex1ldg}.}	
		\label{ldgpolyP1Vnu001}
			\vspace{.4cm}
			\begin{tabular}{|c|c| c| c| c| c| c|}
				
				\hline
				$h$ & $\|\e_\bu^M\|$ & order & $|(\underline{\e}_\sib^M,\,\e_\bu^M)|_\mathcal{A}$ &order &$\|e_p^M\|$ & order\\ 
\hline
\hline 
 1/4 & 0.049534 &  & 3.89963 &  & 0.040452 & \\ 
 1/8 & 0.011249 & 2.1385 & 1.94561 & 1.0031 & 0.021209 & 0.9315\\ 
 1/16 & 0.002498 & 2.1708 & 0.944718 & 1.0422 & 0.011711 & 0.8568\\ 
 1/32 & 0.000581 & 2.1046 & 0.463555 & 1.0271 & 0.006169 & 0.9246\\ 
 1/64 & 0.000139 & 2.0561 & 0.229352 & 1.0151 & 0.003167 & 0.9617\\ 
	         		\hline
	           
			\end{tabular}
			\vspace{.1cm}
		\end{table}
	}

{\small
		\renewcommand{\arraystretch}{1.3}
		\begin{table}[ht!]
			\centering
			\caption{
Numerical errors and corresponding convergence rates for the temperature gradient and temperature with $P_0$--$P_1$ element for $\nu=0.01$ in Example \ref{ex1ldg}.}	
		\label{ldgpolyP1Tnu001}
			\vspace{.4cm}
			\begin{tabular}{|c|c| c| c| c| c| c|}
				
				\hline
				$h$ & $\|e_\theta^M\|$ & order & $|(\e_{\bchi}^M,\, e_\theta^M)|_\mathcal{A}$ &order\\ 
\hline 
\hline
 1/4 & 0.000381 &  & 0.023196 &  \\ 
 1/8 & 0.000092 & 2.0327 & 0.011742 & 0.9821\\ 
 1/16 & 0.000023 & 1.9742 & 0.005911 & 0.9903\\ 
 1/32 & 0.000005 & 1.9852 & 0.002962 & 0.9966\\ 
 1/64 & 0.000001 & 1.9969 & 0.001482 & 0.9987\\ 
	         		\hline
	           
			\end{tabular}
			\vspace{.1cm}
		\end{table}
	}

	  {\small
		\renewcommand{\arraystretch}{1.3}
		\begin{table}[ht!]
			\centering
			\caption{
Numerical errors and corresponding convergence rates for the stress tensor, velocity, and pressure with $P_01$--$P_2$--$P_1$ element for $\nu=0.01$ in Example \ref{ex1ldg}.}	
		\label{ldgpolyP2Vnu001}
			\vspace{.4cm}
			\begin{tabular}{|c|c| c| c| c| c| c|}
				
				\hline
				$h$ & $\|\e_\bu^M\|$ & order & $|(\underline{\e}_\sib^M,\,\e_\bu^M)|_\mathcal{A}$ &order &$\|e_p^M\|$ & order\\ 
\hline 
\hline
 1/4 & 0.000310178 &  & 0.0292073 &  & 0.0003141 & \\ 
 1/8 & 0.0000449 & 2.7878 & 0.0082658 & 1.8211 & 0.0000798 & 1.9757\\ 
 1/16 & 0.0000053 & 3.0749 & 0.0021781 & 1.9241 & 0.0000205 & 1.9579\\ 
 1/32 & 0.0000006 & 3.0588 & 0.0005559 & 1.9701 & 0.0000051 & 1.9865\\
	         		\hline
	           
			\end{tabular}
			\vspace{.1cm}
		\end{table}
	}

{\small
		\renewcommand{\arraystretch}{1.3}
		\begin{table}[ht!]
			\centering
			\caption{
Numerical errors and corresponding convergence rates for the temperature gradient and temperature with $P_1$-$P_2$ element for $\nu=0.01$ in Example \ref{ex1ldg}.}	
		\label{ldgpolyP2Tnu001}
			\vspace{.4cm}
			\begin{tabular}{|c|c| c| c| c| c| c|}
				
				\hline
				$h$ & $\|e_\theta^M\|$ & order & $|(\e_{\bchi}^M,\, e_\theta^M)|_\mathcal{A}$ &order\\ 
\hline 
\hline
 1/4 & 0.0000690 & & 0.012456 & \\ 
 1/8 & 0.0000082 & 3.0566 & 0.0033884 & 1.8781\\ 
 1/16 & 0.0000008 & 3.2125 & 0.0008661 & 1.9679\\ 
 1/32 & 0.0000001 & 3.1316 & 0.0002178 & 1.9911\\ 
	         		\hline
	           
			\end{tabular}
			\vspace{.1cm}
		\end{table}
	}

	  {\small
		\renewcommand{\arraystretch}{1.3}
		\begin{table}[ht!]
			\centering
			\caption{
Stress tensor, velocity and  pressure  errors and order of convergence with $P_0$--$P_1$--$P_0$ element, $\nu=0.001$   
for Example \ref{ex1ldg}.}	
		\label{ldgpolyP1Vnu0001}
			\vspace{.4cm}
			\begin{tabular}{|c|c| c| c| c| c| c|}
				
				\hline
				$h$ & $\|\e_\bu^M\|$ & order & $|(\underline{\e}_\sib^M,\,\e_\bu^M)|_\mathcal{A}$ &order &$\|e_p^M\|$ & order\\ 
\hline 
\hline
 1/4 & 0.241519 &  & 22.3938 & & 0.147269 & \\ 
 1/8 & 0.095449 & 1.3393 & 17.648 & 0.3435 & 0.031182 & 2.2396\\ 
 1/16 & 0.024058 & 1.9882 & 9.29434 & 0.92508 & 0.012499 & 1.3188\\ 
 1/32 & 0.005611 & 2.1001 & 4.6155 & 1.0098 & 0.006261 & 0.9975\\ 
 1/64 & 0.001340 & 2.0658 & 2.28945 & 1.0115 & 0.003174 & 0.9798\\ 
 
	         		\hline
	           
			\end{tabular}
			\vspace{.1cm}
		\end{table}
	}

{\small
		\renewcommand{\arraystretch}{1.3}
		\begin{table}[ht!]
			\centering
			\caption{
Temperature gradient and temperature  errors and order of convergence with $P_0$--$P_1$ element, $\nu=0.001$   
for Example \ref{ex1ldg}.}	
		\label{ldgpolyP1Tnu0001}
			\vspace{.4cm}
			\begin{tabular}{|c|c| c| c| c| c| c|}
				
				\hline
				$h$ & $\|e_\theta^M\|$ & order & $|(\e_{\bchi}^M,\, e_\theta^M)|_\mathcal{A}$ &order\\ 
\hline 
\hline
 1/4 & 0.000396 &  & 0.023201 & \\ 
 1/8 & 0.000108 & 1.8708 & 0.011746 & 0.982\\ 
 1/16 & 0.000028 & 1.9324 & 0.005911 & 0.9906\\ 
 1/32 & 0.000007 & 1.9742 & 0.002962 & 0.9968\\ 
 1/64 & 0.000001 & 1.9898 & 0.001482 & 0.9988\\ 
	         		\hline
	           			\end{tabular}
			\vspace{.1cm}
		\end{table}
	}

	  {\small
		\renewcommand{\arraystretch}{1.3}
		\begin{table}[ht!]
			\centering
			\caption{
Stress tensor, velocity and  pressure  errors and order of convergence with $P_1$--$P_2$--$P_1$ element, $\nu=0.001$    
for Example \ref{ex1ldg}.}	
		\label{ldgpolyP2Vnu0001}
			\vspace{.4cm}
			\begin{tabular}{|c|c| c| c| c| c| c|}
				
				\hline
				$h$ & $\|\e_\bu^M\|$ & order & $|(\underline{\e}_\sib^M,\,\e_\bu^M)|_\mathcal{A}$ &order &$\|e_p^M\|$ & order\\ 
\hline 
\hline
 1/4 & 0.0002961 & & 0.030582 &  & 0.0000371 & \\ 
 1/8 & 0.0000413 & 2.8415 & 0.008423 & 1.8602 & 0.0000084 & 2.138\\ 
 1/16 & 0.0000052 & 2.9723 & 0.002184 & 1.9469 & 0.0000021 & 2.0275\\ 
 1/32 & 0.0000006 & 3.0273 & 0.0005562 & 1.9736 & 0.0000005 & 1.9965\\ 
 
	         		\hline
	           
			\end{tabular}
			\vspace{.1cm}
		\end{table}
	}

{\small
		\renewcommand{\arraystretch}{1.3}
		\begin{table}[ht!]
			\centering
			\caption{
Temperature gradient and temperature  errors and order of convergence with $P_1$--$P_2$ element, $\nu=0.001$ 
for Example \ref{ex1ldg}.}	
		\label{ldgpolyP2Tnu0001}
			\vspace{.4cm}
			\begin{tabular}{|c|c| c| c| c| c| c|}
				
				\hline
				$h$ & $\|e_\theta^M\|$ & order & $|(\e_{\bchi}^M,\, e_\theta^M)|_\mathcal{A}$ &order\\ 
\hline 
\hline
 1/4 & 0.0000690 &  & 0.012456 & \\ 
 1/8 & 0.0000082 & 3.0566 & 0.003388 & 1.8781\\ 
 1/16 & 0.0000008 & 3.2125 & 0.000866 & 1.9679\\ 
 1/32 & 0.0000001 & 3.1316 & 0.000217 & 1.9911\\
	         		\hline
	           			\end{tabular}
			\vspace{.1cm}
		\end{table}
	}

\begin{figure}[hbtp]
\centering
\subfloat[]{{\includegraphics[width=4.5cm]{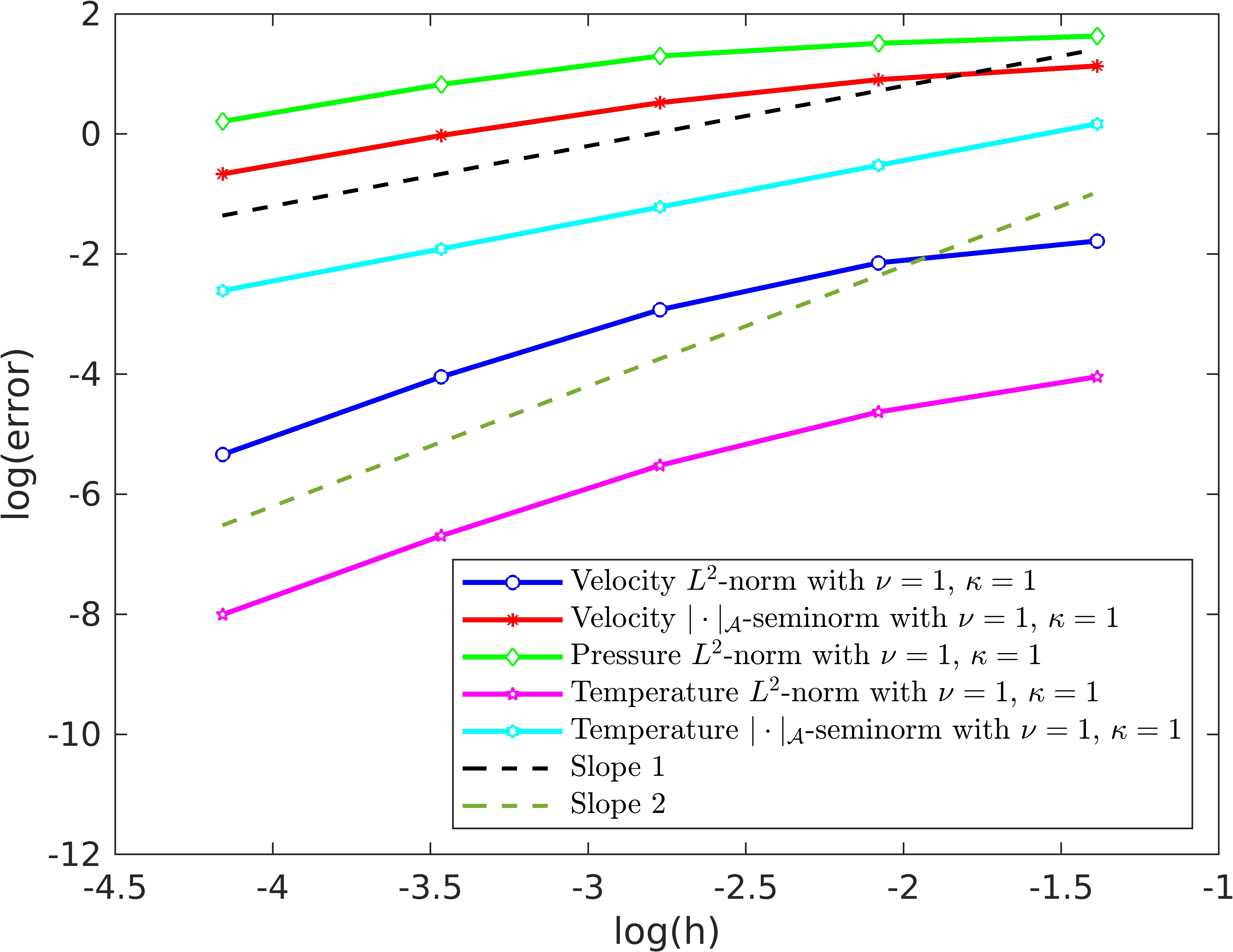}}}
\hspace{0.2cm}
\subfloat[]{{\includegraphics[width=4.5cm]{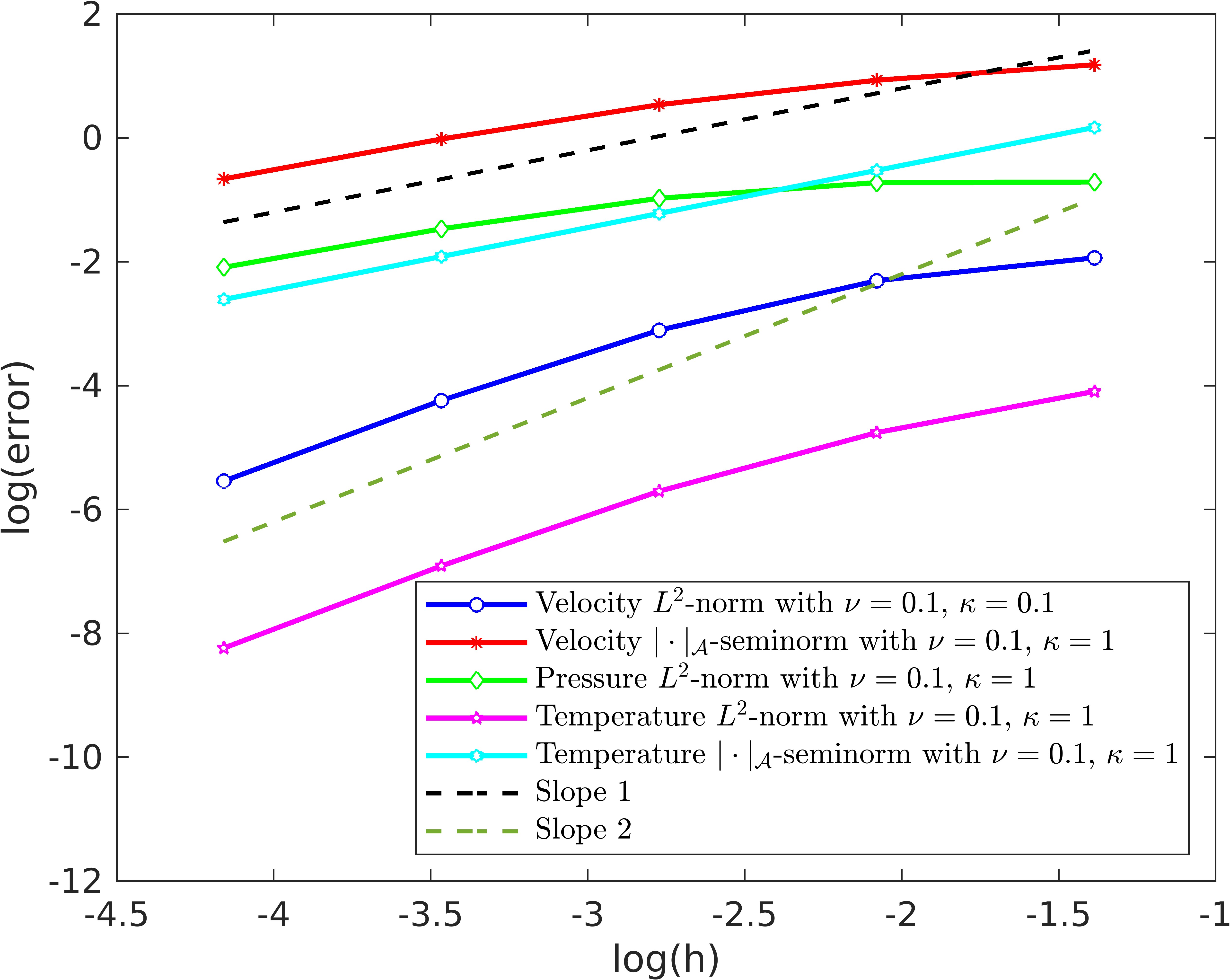}}}
\hspace{0.2cm}
\subfloat[]{{\includegraphics[width=5cm]{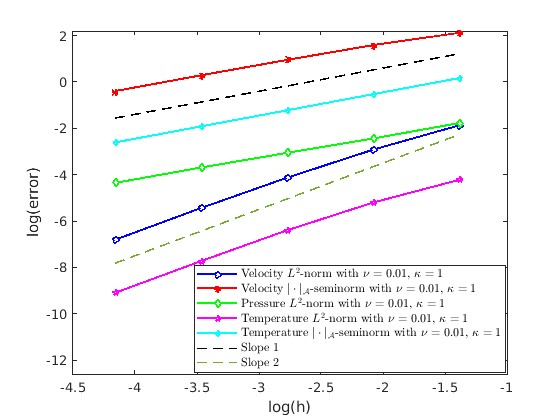}}}
\caption{Convergence graphs (log--errors) for the stress tensor, velocity, pressure, temperature gradient, and temperature with  $P_0$--$P_1$--$P_0$--$P_0$--$P_1$ element for Example \ref{ldgex2}.} 
\label{ldgfig1} 
\centering
\hspace{0.8cm}
\subfloat[]{{\includegraphics[width=4.5cm]{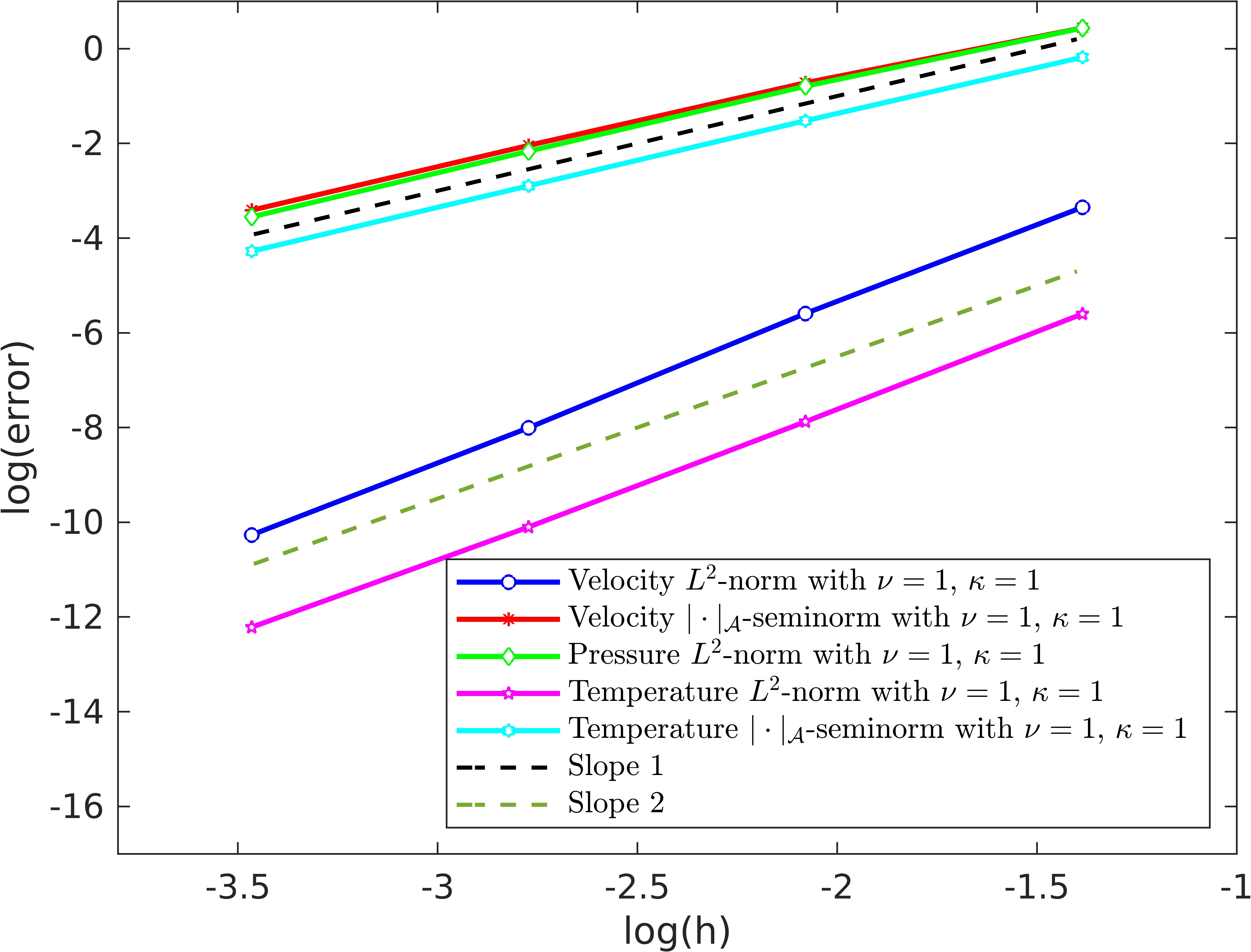}}}
\hspace{0.2cm}
\subfloat[]{{\includegraphics[width=5cm]{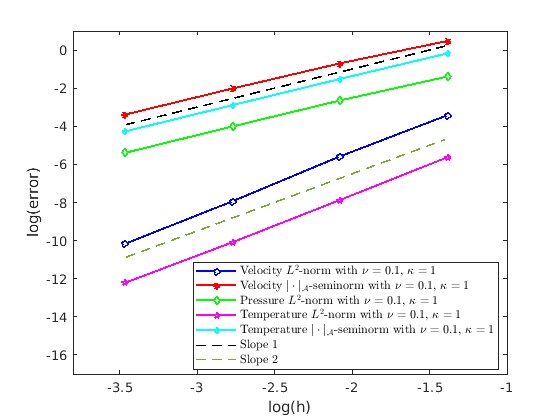}}}
\hspace{0.2cm}
\subfloat[]{{\includegraphics[width=5cm]{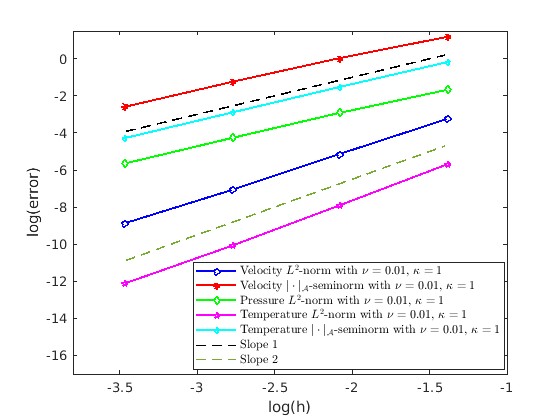}}}
\caption{Convergence graphs (log--errors) for the stress tensor, velocity, pressure, temperature gradient, and temperature with the $P_1$--$P_2$--$P_1$--$P_1$--$P_2$ element for Example \ref{ldgex2}.} 
\label{ldgfig2} 
\end{figure}

\begin{example}\label{ldgex2}
The exact solution is defined as
\begin{align*}
u_1(x, y, t) &= \cos(t) \sin^2(\pi x) \sin(2\pi y),\, 
u_2(x, y, t) = -\cos( t) \sin(2\pi x) \sin^2(\pi y),\\
p(x,y,t)&= \cos(t)\sin(2\pi x)\sin(2\pi y),\\
\theta(x,\,y,\,t)&= \cos(t) \sin^2(\pi x) \sin(2\pi y)-\cos( t) \sin(2\pi x) \sin^2(\pi y).
\end{align*}

\end{example}	
\noindent
Figure \ref{ldgfig1} showcases the numerical error behavior of the LDG method for the stress tensor, velocity, pressure, temperature gradient, and temperature. The analysis is conducted for varying viscosity values \( \nu = \{1, 0.1, 0.01\} \), using the mixed element configurations \( P_0-P_1-P_0 \) and \( P_0-P_1 \).
Figure \ref{ldgfig2}
  presents the LDG errors using the \( P_1-P_2-P_1 \) and \( P_1\text{-}P_2 \) finite element pairs. The simulations are performed up to a final time of \( T = 0.25 \), and the observed error behavior supports the theoretical estimates established in Theorems \ref{bevel} and \ref{FullyP}.

\begin{figure}[hbtp]
\centering
\subfloat[ Velocity vector field]
{{\includegraphics[width=3.5cm]{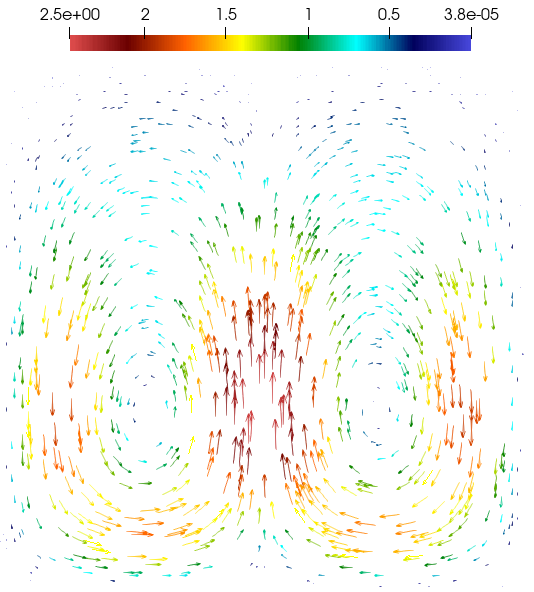}}}
\hspace{0.2 cm}
\subfloat[First component of velocity]
{{\includegraphics[width=3.5cm]{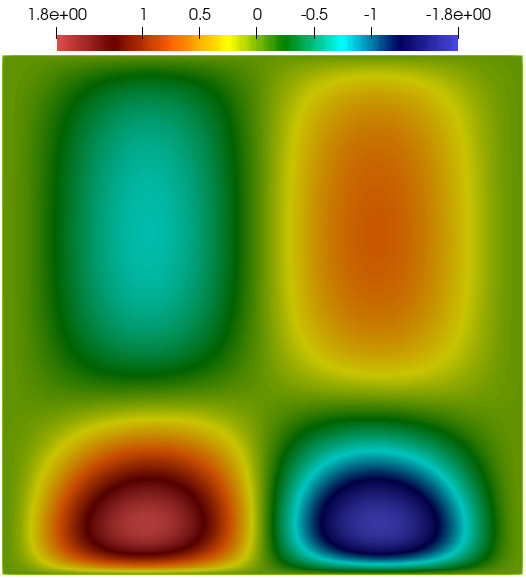}}}
\hspace{0.2 cm}
\subfloat[Second component of velocity]
{{\includegraphics[width=3.5cm]{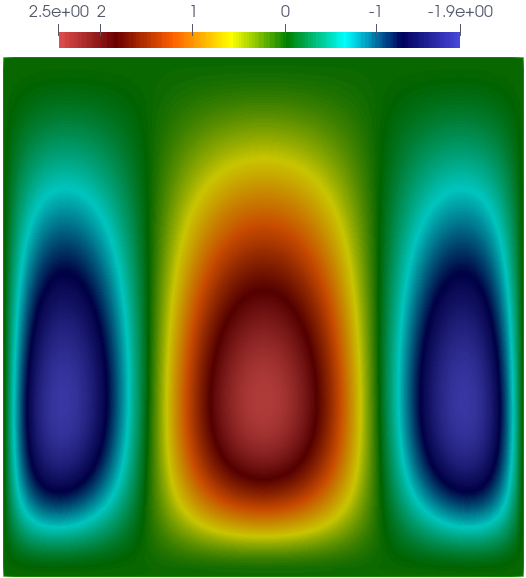}}}
\hspace{0.2cm}
\subfloat[Pressure]
{{\includegraphics[width=3.5cm]{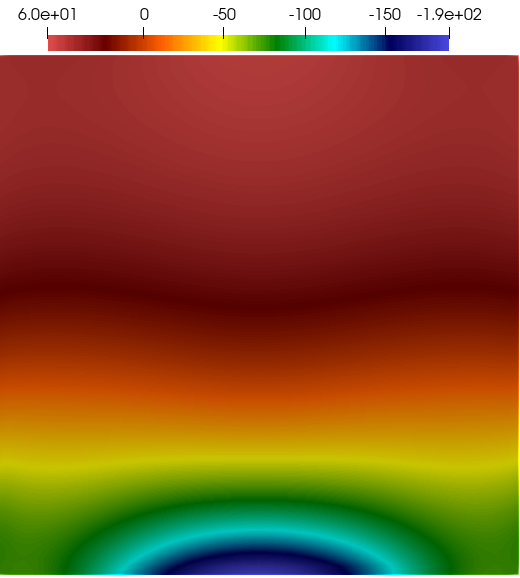}}}\\
\subfloat[Velocity vector field]
{{\includegraphics[width=3.5cm]{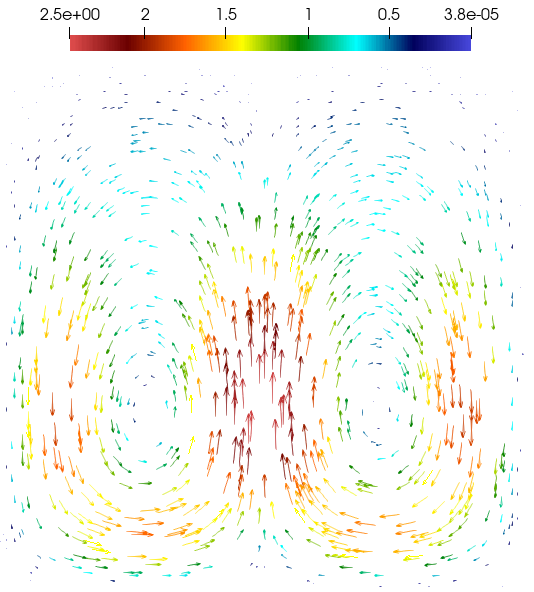}}}
\hspace{0.2cm}
\subfloat[First component of velocity]
{{\includegraphics[width=3.5cm]{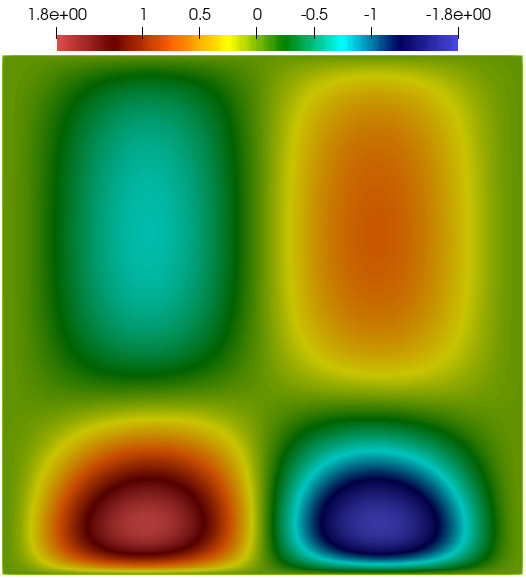}}}
\hspace{0.2cm}
\subfloat[Second component of velocity]
{{\includegraphics[width=3.5cm]{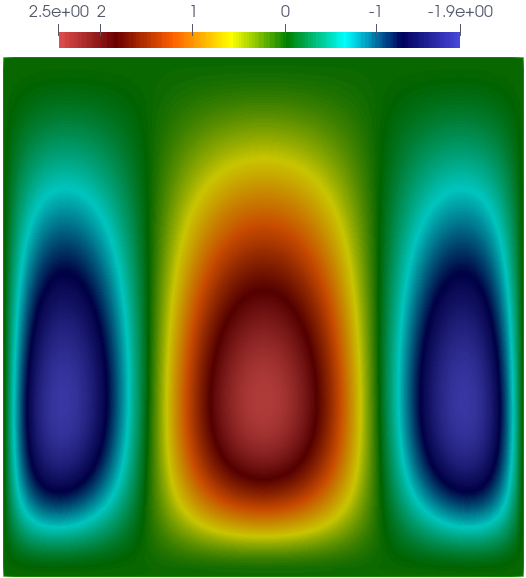}}}
\hspace{0.2cm}
\subfloat[Pressure]
{{\includegraphics[width=3.5cm]{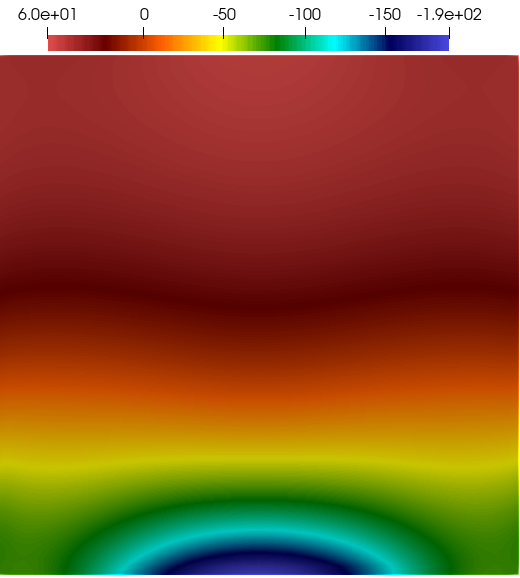}}}
\caption{{ LDG  method for stationary (first row) and transient (second row)   thermal driven cavity problem: graphs of velocity vector field, first and second componenets  of velocity, and pressure with $Ra=2000$.}}
\label{ldgbenchmark2} 
\end{figure}

\begin{figure}[hbtp]
\centering
\subfloat[ Heat--flux]
{{\includegraphics[width=3.5cm]{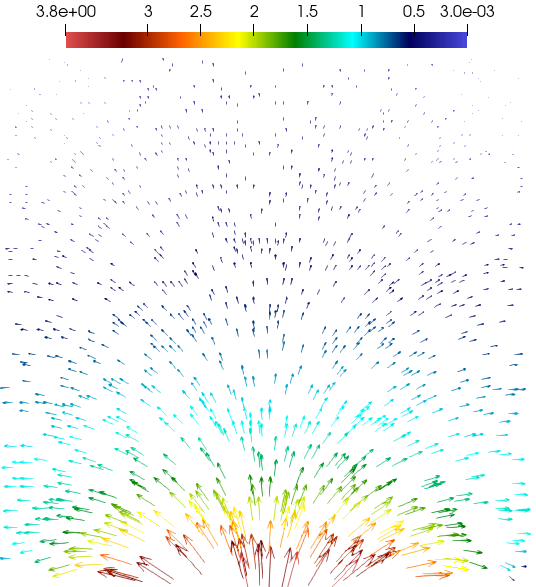}}}
\hspace{0.2cm}
\subfloat[First component of temperature gradient]
{{\includegraphics[width=3.5cm]{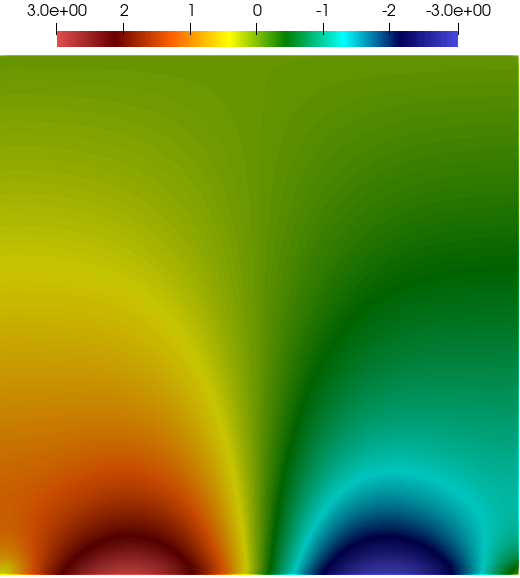}}}
\hspace{0.2cm}
\subfloat[Second component of temperature gradient]
{{\includegraphics[width=3.5cm]{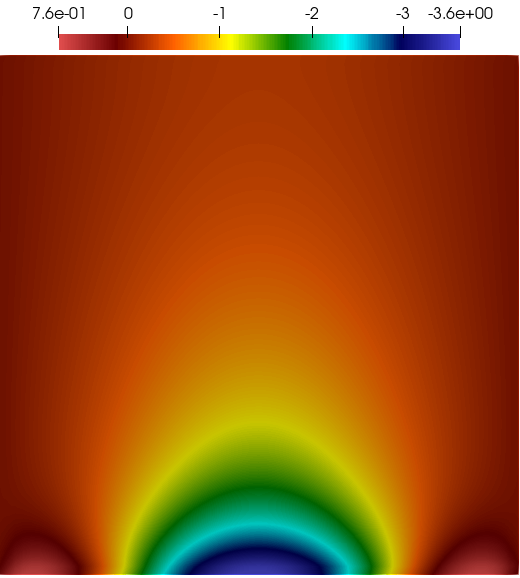 }}}
\hspace{0.2cm}
\subfloat[Temperature]
{{\includegraphics[width=3.5cm]{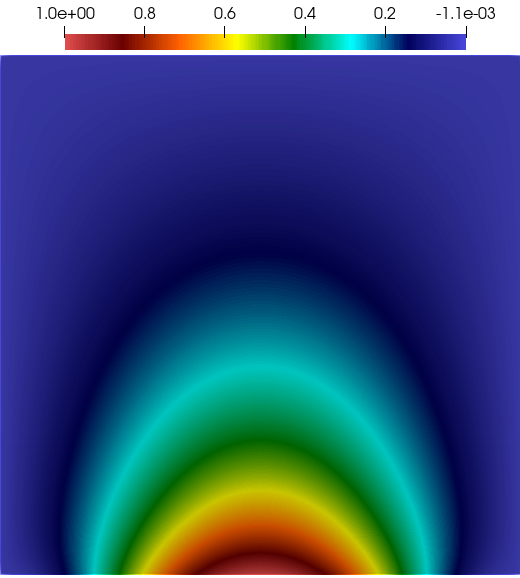}}}\\
\subfloat[heat--flux]
{{\includegraphics[width=3.5cm]{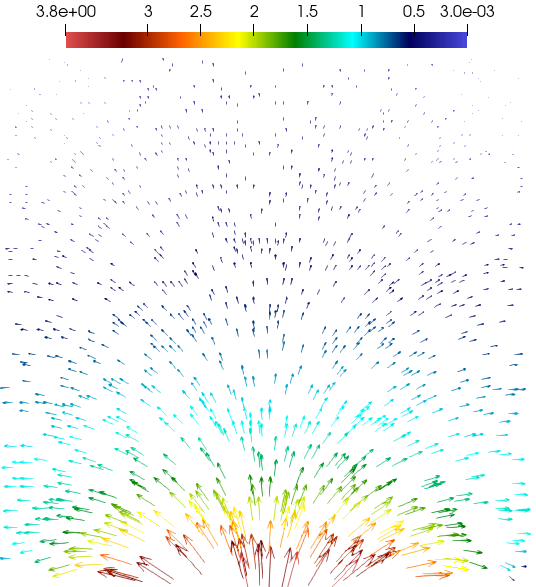}}}
\hspace{0.2cm}
\subfloat[First component of temperature gradient]
{{\includegraphics[width=3.5cm]{./Figure/s_ldgT1  }}}
\hspace{0.2cm}
\subfloat[Second component of temperature gradient]
{{\includegraphics[width=3.5cm]{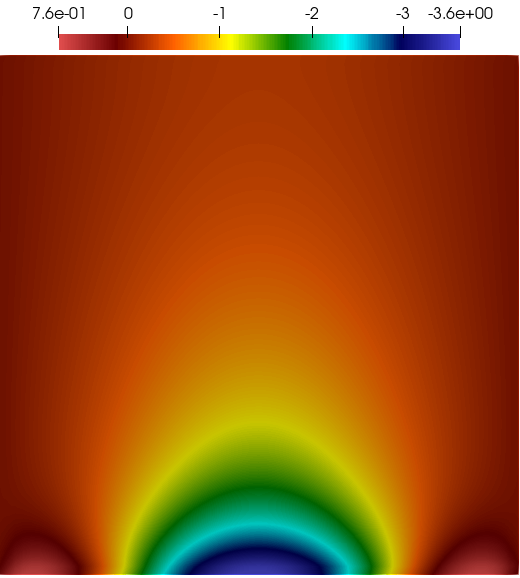 }}}
\hspace{0.2cm}
\subfloat[Temperature]
{{\includegraphics[width=3.5cm]{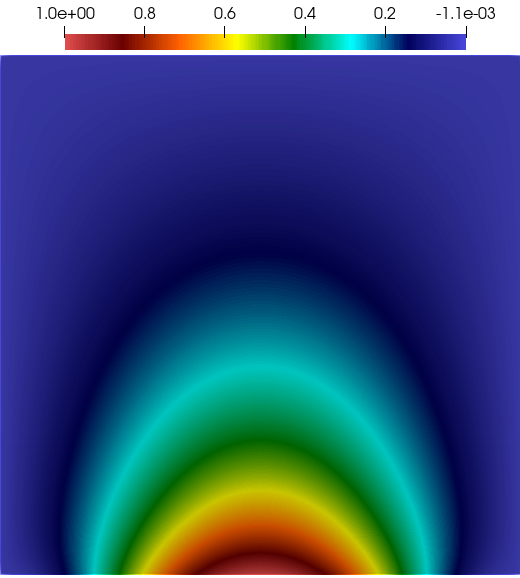}}}
\caption{{ LDG  method for stationary (first row) and transient (second row)  thermal driven cavity problem: graphs of heat--flux,  first and second components of temperature gradient and temperature with $Ra=2000$.}} 
\label{ldgbenchmark3}
\end{figure}

\vspace{1.5cm}
\begin{example}\label{ldgex3}
{\bf{Benchmark problem}}
\end{example}
\noindent
This work focuses on the simulation of thermally driven fluid flow within a square cavity, where the side walls are maintained at different temperatures. The system is governed by key dimensionless parameters: the viscosity \( \nu = 0.5 \), Rayleigh number \( \text{Ra} = 2000 \), and thermal conductivity \( \kappa = 1 \). The coupling parameter \( \boldsymbol{\mu} \) is defined as \( \boldsymbol{\mu} = -\nu \, \text{Ra} (0, 1) \). The velocity field \( \mathbf{u} \) is subject to homogeneous Dirichlet boundary conditions. For the temperature distribution \( \theta \), the top, left, and right walls are set to zero temperature, while the bottom wall is heated according to the profile
\(
\theta(\bx, 0) = 0.5(1 - \cos(2\pi x))(1 - y),
\)
which generates a sinusoidal heating with a peak at \( x = 0.5 \), where \( \theta = 1 \). Figures \ref{ldgbenchmark2}--\ref{ldgbenchmark3}  computed with a refined mesh of size \( h = 1/64 \) and penalty parameters \( c_{11} = c_{22} = 50 \).
The first row in Figures  \ref{ldgbenchmark2} present the steady--state numerical approximations of velocity vector field and pressure, components of velocity    for the stationary Boussinesq system where as second columns highlight the numerical  approximations  of velocity vector field and pressure, component of velocity, for  the  transient Boussinesq system of equations.  These visualizations reveal the formation of symmetric convection patterns within the cavity, consistent with established results in the literature \cite{dd2006}.
 Additionally, Figure \ref{ldgbenchmark3}  illustrates  the heat flux vectors,  the temperature, and  the temperature gradient components. The observed direction of the heat flux aligns with physical expectations flowing from hotter to cooler regions.  
Finally, the results demonstrate that as time advances toward \( T = 40 \), the transient Boussinesq solution gradually stabilizes, approaching the steady-state solution, in agreement with previously reported findings~\cite{COTD}.

\begin{figure}[hbtp]
\centering
{{\includegraphics[width=6cm]{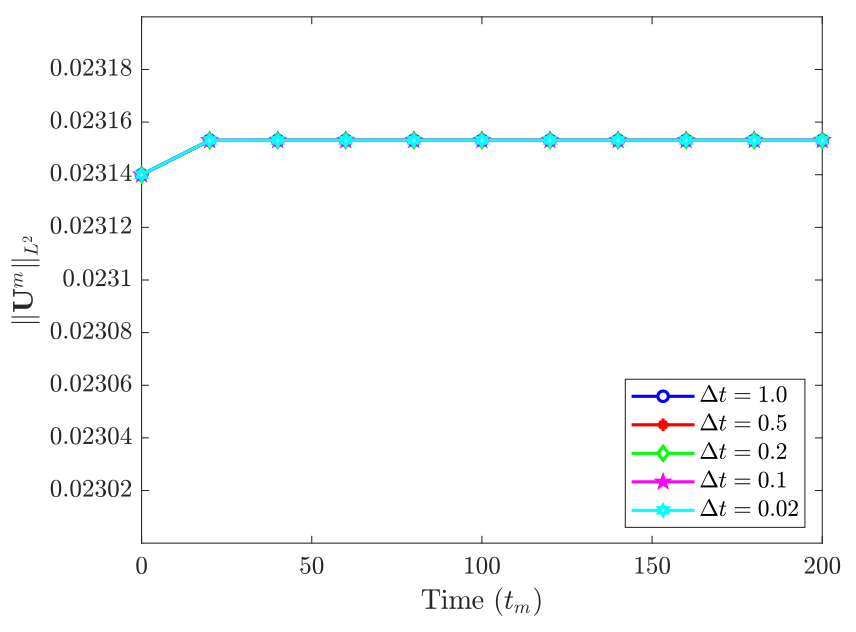}}}
\hspace{1cm}
{{\includegraphics[width=6cm]{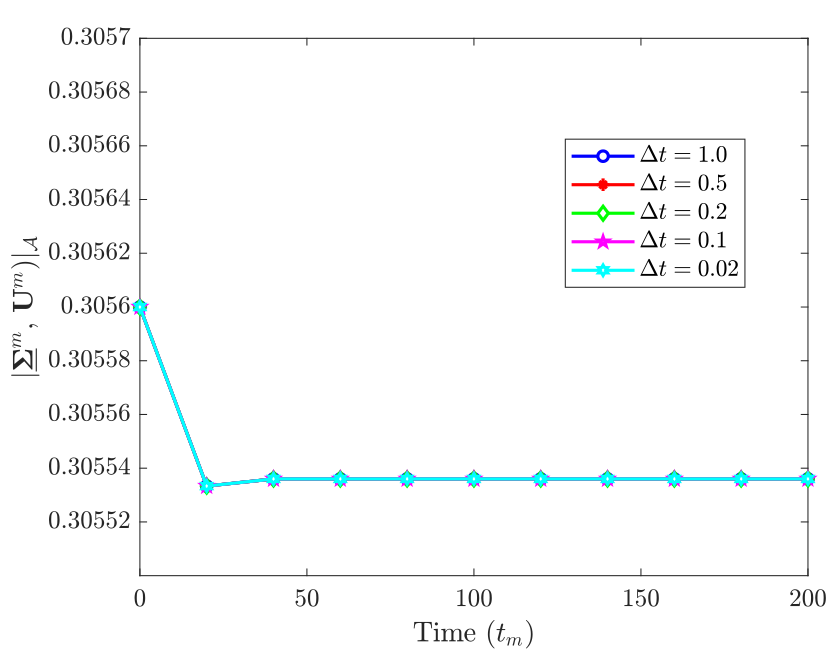  }}}\\
{{\includegraphics[width=6cm]{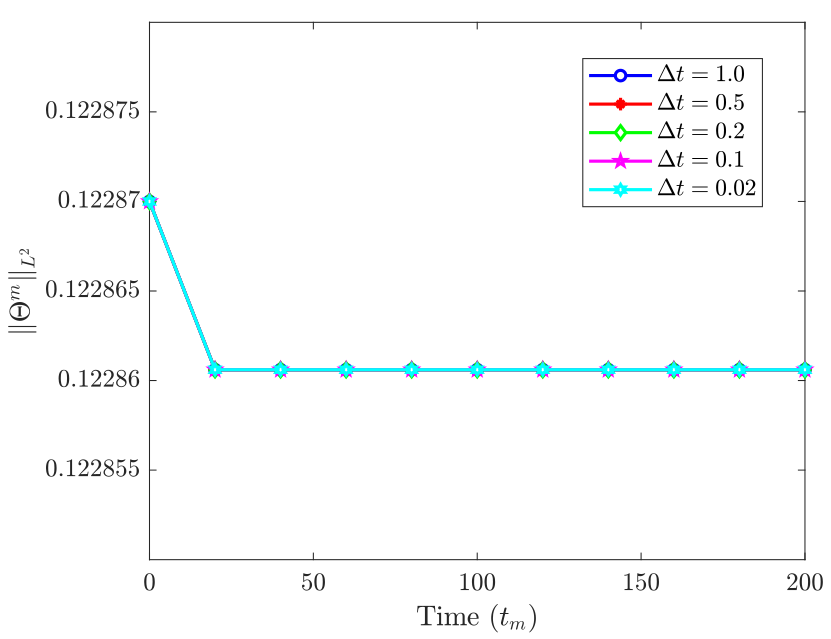 }}}
\hspace{1cm}
{{\includegraphics[width=6cm]{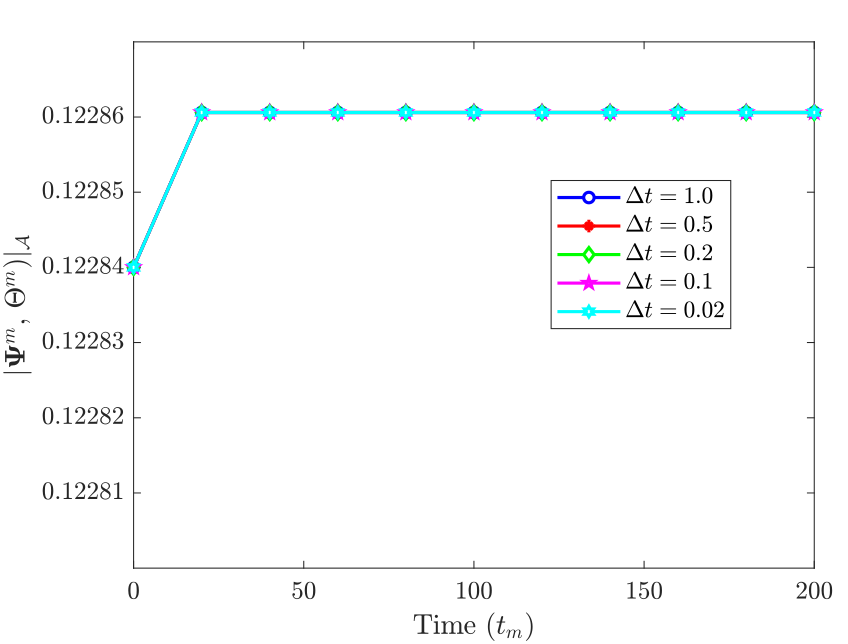 }}}
{{\includegraphics[width=6cm]{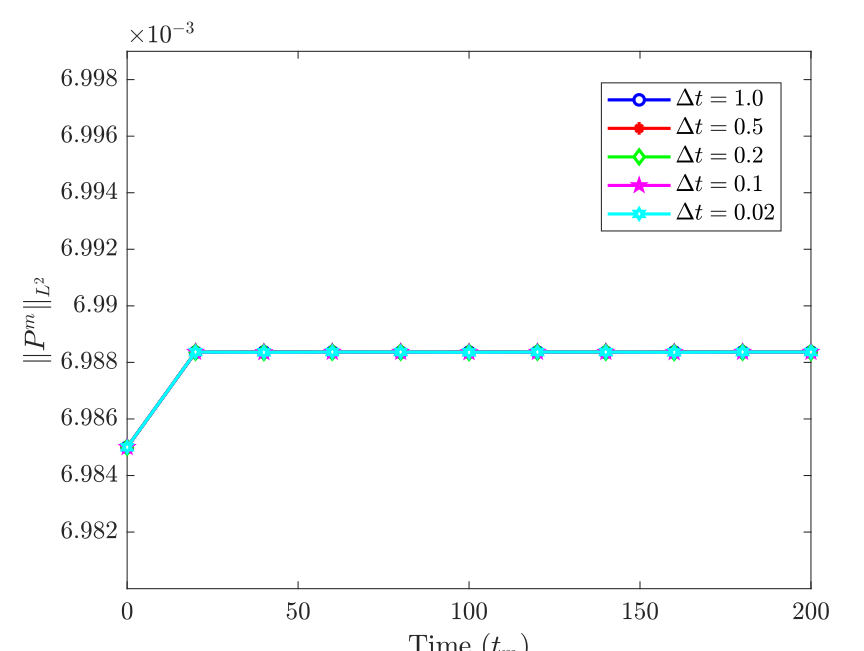 }}}
\caption{{Time evolution of fully discrete stress tensor,   velocity, pressure,  temperature gradient, temperature  for varying 
$\Delta t$ in the thermal driven cavity problem with $Ra=2000$, 
illustrating long--time stability.}}
\label{stabileLDG} 
\end{figure}

\noindent { Furthermore, we examine the stability of the fully discrete stress tensor, velocity, pressure, 
temperature gradient, and temperature, denoted by 
$(\Su^m,\,\bU^m,\,P^m,\,\boldsymbol{\Psi}^m,\,\Theta^m)_{m \geq 1}$ 
at time $t_m$ through the quantities 
\(
|(\Su^m,\,\bU^m)|_{\mathcal{A}}\), 
\(\|\bU^m\|, \quad 
|(\boldsymbol{\Psi}^m,\,\Theta^m)|_{\mathcal{A}}, \quad  
\|\Theta^m\|, \quad 
\|P^m\|.
\)
Under suitable smallness conditions, it has been established that these quantities remain 
uniformly bounded in time, independent of the timestep $\Delta t$ and the mesh size $h$. 
To verify this result numerically, we plot them  against time for different choices of $\Delta t$.  \\
For $\nu=1/100$, the results are presented in Figure~\ref{stabileLDG}. 
The plots show that the the quantities remain bounded for all $t_m$ and are insensitive to the choice of $\Delta t$, 
thereby confirming the long--time stability predicted by the theory.
}

\section{Conclusion}\label{ldgsection8}
The present work has applied the LDG method to the Boussinesq system (\ref{1.1})--(\ref{1.3}), providing novel contributions to its numerical analysis. We have formulated a semidiscrete LDG method for the proposed system and obtained {\it a priori} regularity estimates for the velocity and temperature variables.
To carry out the error analysis, we have introduced the $\bL^2$--projection operators 
$\underline{\bP}_h^\sib$, $\bP_h^\bu$, $\mathcal{P}_h^\bchi$, and $\mathcal{P}_h^\theta$, together with the modified Stokes operators 
$\underline{\boldsymbol{\bS}}_h^\sib$, $\boldsymbol{\bS}_h^\bu$, $\boldsymbol{\bS}_h^\bchi$, and $\mathcal{S}_h^\theta$ in appropriate DG spaces, and established their approximation properties. 
The estimates derived for the projection operators 
$\underline{\bP}_h^\sib$, $\bP_h^\bu$, $\mathcal{P}_h^\bchi$, and $\mathcal{P}_h^\theta$ play a crucial role in obtaining corresponding bounds for the modified Stokes operators 
$\underline{\boldsymbol{\bS}}_h^\sib$, $\boldsymbol{\bS}_h^\bu$, $\boldsymbol{\bS}_h^\bchi$, and $\mathcal{S}_h^\theta$.  Using weak solution regularity, duality arguments, and the approximation properties of \(\underline{\bP}^\sib_h\), \(\bP^\bu_h\), $\mathcal{P}_h^\bchi$, $\mathcal{P}_h^\theta$,  \(\underline{\boldsymbol{\bS}}^\sib_h\), \(\bS^\bu_h\), $\boldsymbol{\bS}_h^\bchi$, and $\mathcal{S}_h^\theta$, we have obtained optimal error estimates for the stress tensor, velocity, pressure, temperature, and its gradient. Furthermore, under a suitable smallness condition on the given data, we have shown that the semidiscrete error estimates are uniform in time. For full discretization, we have employed a backward Euler method, and achieved optimal convergence rates for the numerical approximations. The numerical results verify the theoretical analysis and demonstrate that the fully discrete LDG approach remains accurate and stable for small values of the viscosity $\nu$ and thermal conductivity $\kappa$.

\section*{Acknowledgement}
The first author expresses sincere gratitude to the Council of Scientific \& Industrial Research (CSIR) for financial support (09/1290(0001)/2019--EMR--I). The second author acknowledges support from the Anusandhan National Research Foundation (ANRF), Government of India, under the Core Research Grant project no. CRG/2021/008278. Additionally, the first, second, and third authors express their gratitude to ANRF, Government of India, for support under the Core Research Grant project no. CRG/2023/006366.\\
 
\vspace*{0.5cm}
\noindent
\textbf{Declarations of interest}: None.

\noindent
\textbf{Author Contributions}:
All authors contributed equally to this work. And all authors read and approved the final manuscript.

\noindent
\textbf{Data Availability}: Not applicable.
\section*{References}
	\begingroup
	\renewcommand{\section}[2]{}


\begin{thebibliography}{}
\bibitem{A75}
R.A. Adams, Sobolev Spaces, Academic Press, New York, 1975.
\bibitem{BD24B}
Bajpai, S. and Swain, D. K.
On error estimates of a discontinuous Galerkin method of the Boussinesq system of equations, {\it 
Comput. Methods Appl. Math.}, 25 (2), 313--348, 2024.
	\bibitem{Benard}
	    {Bénard, H. Rev. Gen. Sci. Pure Appl. 11, 1261 (1900); 1309 (1900); Ann. Chim. Phys. 23, 62 (1)}.
	    \bibitem{BHH} Burns, J. A.,  He, X. M. and Hu, W.{ \it Control of the Boussinesq equations and implications for sensor location in energy efficient buildings}, in: Proceedings of the 2012 American Control Conference, Montreal, CA,  2232--2237, 2012.
	    \bibitem{CCPS2000}
	    Castillo, P., Cockburn, B., Perugia, I. and Schötzau, D. An a priori error analysis of the local discontinuous Galerkin method for elliptic problems, {\it SIAM J. Numer. Anal.}, 38 (5), 1676--1706, 2000.
	\bibitem{CCQ17} Cesmelioglu, A., Cockburn, B. and Qiu, W.  { Analysis of a hybridizable discontinuous Galerkin method for the steady-state incompressible Navier--Stokes equations}, {\it Math. Comp.,} 86, 1643--1670, 2017.
	\bibitem{cockburn00} Cockburn, B., Kanschat, G., Perugia, I. and  Schötzau, D.  Superconvergence of the local discontinuous Galerkin method for elliptic problems on cartesian grids*, {\it SIAM J. Numer. Anal.}, 39 (1),  264–285, 2000.    
	    
\bibitem{cockburn02}
Cockburn, B., Kanschat, G.,   Schötzau, D. and  Schwab, C.  Local discontinuous Galerkin methods for the Stokes system, \emph{ SIAM J. Numer. Anal.}, 40 (1), 319--343, 2002.

\bibitem{cockburn03}
Cockburn, B., Kanschat, G. and Schötzau, D.  The local discontinuous Galerkin method for the Oseen equations, {\it Math. Comp.}, 73 (246),  569–593, 2003. 
\bibitem{CKS09}  Cockburn, B.,  Kanschat, G. and  Sch{\"o}tzau, D. {An equal--order DG Method for the incompressible Navier--Stokes equations},  {\it J. Sci. Comput.,}  40, 188--210, 2009.
\bibitem{cockburn04} 
Cockburn, B., Kanschat, G. and Schötzau, D. A locally conservative LDG method for the incompressible Navier--Stokes equations, {\it Math. Comp.}, 74 (251), 1067–1096, 2004.
\bibitem{CKS00} Cockburn, B., Karniadakis, G., and Shu, C.-W. Discontinuous Galerkin methods. theory, computation and
applications, {\it Lect. Notes Comput. Sci. Eng.}, 11, Springer--Verlag, Berlin, 2000.
\bibitem{dd2006}
		Dalal, A. and Das, M. Natural convection in a rectangular cavity heated from below and uniformly cooled from the top and both sides, \emph{ Numer. Heat Transf. Part A Appl.}, 49 (3), 301--322, 2006.
\bibitem{DE12} Di Pietro, D. A. and Ern, A. Mathematical aspects of discontinuous Galerkin methods,  Springer Science and  Business Media, 2012.
\bibitem{GR09} Girault, V. and  Rivi{\`e}re, B.  DG approximation of coupled Navier--Stokes and Darcy equations by Beaver--Joseph--Saffman interface condition. \emph {SIAM J. Numer. Anal.}, 47, 2052--2089, 2009.
\bibitem{GRW005} Girault, V.,  Rivi{\`e}re, B. and Wheeler, M. F. A discontinuous Galerkin method with non--overlapping domain decomposition for the Stokes and Navier--Stokes problems. \emph{ Math. Comp.,} 74, 53--84, 2005.
\bibitem{GRW05} Girault, V.,  Rivi{\`e}re, B. and Wheeler, M. F.  A splitting method using discontinuous Galerkin for the transient incompressible Navier--Stokes equations. \emph{ ESAIM: Math. Model. Numer. Anal.}, 39, 1115--1147, 2005.

\bibitem{guzman08} 
 Guzmán, J.  Local and pointwise error estimates of the local discontinuous Galerkin method applied to the Stokes problem, {\it Math. Comp.}, 77 (263),  1293--1322, 2008.

  \bibitem{OK98}
Karakashian, O. A. and  Jureidini, W. N. A nonconforming ﬁnite element method for the stationary Navier--Stokes equations, {\it SIAM Numer. Anal.}, 35 (1), 93--120, 1998

\bibitem{KR05} Kaya, S. and Rivi{\`e}re, B.  {A discontinuous subgrid eddy viscosity method for the time--dependent Navier--Stokes equations}, {\it SIAM J. Numer. Anal.,} 43, 1572--1595, 2005.
\bibitem{KR19} Kirk, K. L. A. and Rhebergen, S.  {Analysis of a pressure-robust hybridized discontinuous Galerkin method for the stationary Navier--Stokes equations}, {\it J. Sci. Comput.}, 81, 881--897, 2019.
\bibitem{LR74}  Lesaint, P. and  Raviart, P. A.  On a finite element method for solving the neutron transport equation, in mathematical aspects of finite element methods in partial differential equations., \emph{C.A. de Boor Ed.,}  Academic Press, 89--123, 1974.
\bibitem{MLR22} Masri, R., Liu, C. and Riviere, B., {A discontinuous Galerkin pressure correction scheme for the incompressible Navier--Stokes equations: Stability and convergence},  {\it Math. Comp.}, 91, 1625--1654, 2022.
\bibitem{KO}
        \'{O}eda, K. {\it On the initial value problem for the heat convection equation of Boussinesq approximation in a time--dependent domain}, Proc. Japan Acad. 64, Ser. A, 143--146, 1988.

 \bibitem{COTD}
Oyarzua, R., Qin, T. and Schotzau, D.  An exactly divergence free finite element method for a generalized Boussinesq problem, \emph{ IMA J. Numer. Anal.}, 34 (3),  1104--1135, 2013.
\bibitem{R08}   Rivi{\`e}re, B. Discontinuous Galerkin methods for solving elliptic and parabolic equations: theory and implementation, \emph{SIAM Philadelphia}, 2008.

\bibitem{T77}
 Temam, R. Navier--Stokes equations, theory and numerical analysis, studies in mathematics and its applications, North--Holland Publishing Co., Amsterdam, 1976.
 \bibitem{T95}
 Temam, R.  Navier--Stokes equations and nonlinear functional analysis. second ed., CBMS--NSF regional conference series in applied mathematics, 66, SIAM, Philadelphia, PA, 1987.
\bibitem{WDYW17} Wang, S., Deng, W., Yuan, J. and Wu, Y., {Characteristic local discontinuous Galerkin methods for incompressible Navier-Stokes equations}, {\it Commun. Comput. Phys.,} 22, 202--227, 2017.
\end{thebibliography}
\end{document}